\UseRawInputEncoding
\documentclass[12pt]{article}
\usepackage{bm}
\usepackage{mathptmx}
\usepackage{amsfonts}
\usepackage{latexsym}
\usepackage{graphicx}
\usepackage{color}
\usepackage{ifpdf}
\usepackage[titletoc]{appendix}
\newtheorem{theorem}{Theorem}[section]
\newtheorem{corollary}[theorem]{Corollary}
\newtheorem{lemma}[theorem]{Lemma}

\newtheorem{proposition}{Proposition}[section]
\newtheorem{definition}{Definition}[section]
\newtheorem{remark}{Remark}[section]
\newtheorem{problem}{Problem}[section]
\usepackage{amsfonts,amssymb,eucal,amsmath}
\title{\bf On the Problem of  Infinite Spin in  Total
Collisions of the Planar $N$-Body Problem} 

\author{{ Xiang Yu\footnote{Email:yuxiang@swufe.edu.cn, xiang.zhiy@gmail.com}} \\
\small \it School of Economic and Mathematics, Southwestern
University of Finance and Economics, \\
\small \it Chengdu 611130, China}

\date{}

\begin{document}
\maketitle

\begin{abstract}
 For the planar $N$-body problem, we first introduce a class of moving frame suitable for orbits near central configurations, especially for total collision orbits, which is the main new ingredient of this paper. The  moving frame allows us to reduce the degeneracy of the problem according to intrinsic symmetrical characteristic of the $N$-Body problem. First, we give a full answer to the infinite spin or $Painlev\acute{e}$-$Wintner$ problem in the case corresponding to nondegenerate central configurations. Then following some original ideas of C.L. Siegel, especially the idea of normal forms, and applying the theory of  central manifolds, we give a partial answer to the problem  in the case corresponding to degenerate central configurations.  We completely answer the problem  in the case corresponding to  central configurations with one degree of degeneracy. Combining some results on the planar nonhyperbolic equilibrium point, we give a criterion for the case corresponding to   central configurations with two  degrees of degeneracy.  We further answer the problem  in the case corresponding to all  known degenerate central configurations of four-body. Therefore, we solve the problem for almost every choice
of the masses of the four-body problem.  Finally, we give a measure of the set of initial conditions leading to total collisions.\\



\end{abstract}

\begin{center}
\tableofcontents
\end{center}

\section{Introduction}
\ \ \ \  We consider $N $ particles with positive masses moving in an Euclidean plane ${\mathbb{R}}^2$ interacting under
the law of universal gravitation. Let the $k$-th particle have mass $m_k$ and position
$\mathbf{r}_k \in {\mathbb{R}}^2$ ($k=1, 2, \cdots, N$), then the equations of motion of the $N$-body problem are written as
\begin{equation}\label{eq:Newton's equation1}
m_k \ddot{\mathbf{r}}_k =\sum_{1 \leq j \leq N, j \neq k} \frac{m_km_j(\mathbf{r}_j-\mathbf{r}_k)}{|\mathbf{r}_j-\mathbf{r}_k|^3}, ~~~~~~~~~~~~~~~k=1, 2, \cdots, N.
\end{equation}
where $|\cdot|$  denotes the Euclidean norm in ${\mathbb{R}}^2$. 

For the $N$-body problem, collision singularities are inevitable focuses, and also the main difficulties.

It is relatively simple to understand binary collisions of the $N$-body problem. Indeed,  one can change variables so that a binary
collision transforms to a regular point of the equations for the two-body problem \cite{Moser2010Regularization}. Such a
transformation is called a regularization of the binary collision. The
solution can then be extended through the singularity. Sundman \cite{Sundman} showed that binary collisions can also be regularized
in the three-body problem. That is, one can transform the variables in
such a way that the solution can be continued through the binary
collision as an analytic function of a new time variable. This is also  true for several binary collisions occurring simultaneously in the $N$-body problem \cite{Sperling1970On}.

Collisions involving more than two particles are more
complicated, only some partial results  are known. Consider
the \emph{normalized configuration} of the particles to be the configuration divided by a norm
which corresponds physically to the moment of inertia. Sundman \cite{Sundman}
 showed that, for triple collision in the three-body problem, the normalized
configuration approaches the set of central configurations
(cf. \cite{wintner1941analytical} and Sections 2 and 3 below). Wintner \cite{wintner1941analytical}  observed that
Sundman's techniques can be used to show that the normalized
configurations of solutions ending in total
collision in the $N$-body problem also approach the set of central configurations. This is also true for general collision singularity of the $N$-body problem in which several clusters of particles collapse simultaneously \cite{Sperling1970On}.

It is natural to ask whether this implies that
the normalized
configuration  of the particles must approach a certain central configuration, or, may the  normalized
configuration  of the particles make
an infinite number of revolutions before arriving at a collision. This is a long standing open problem on the collision singularity of the Newtonian $N$-body problem. This
problem was posed by Painlev\'{e} and discussed by Wintner \cite[p.283]{wintner1941analytical} in the total
collision of the $N$-body problem. So this problem is usually called the $Painlev\acute{e}$-$Wintner$ problem or the problem of infinite spin. For simplicity, the abbreviation ``$PISPW$" will be used to mean ``the problem of infinite spin or $Painlev\acute{e}$-$Wintner$" in this paper.

Although there has been tremendous interest in the problem, so far, only  few research progress has been made. Indeed,  ones knew that $PISPW$ could be solved in the case corresponding to nondegenerate central configurations for a long time \cite{Chazy,C1967Lectures,Saari1984The}; for example,  one of the ideas is to apply  the theory of normally hyperbolic invariant manifolds, 
 and there are several papers  mentioned this \cite[etc]{Saari1984The}. On the other hand, because little is known about central configurations for $N > 3$ \cite[etc]{zbMATH01516323,Palmore1976Measure,Hampton2006Finiteness,zbMATH06074021}, especially on the degeneracy of central configurations,  $PISPW$ is completely solved only for the three-body problem.

All in all, though several claims of even stronger results were published \cite[etc]{Saari1984The}, however, as  pointed out in \cite{Albouy2012Some},  the proofs in the corresponding papers could not be found up till now. Thus recently Chenciner and Venturelli  \cite{Albouy2012Some} asked for an  solution  to $PISPW$ even in the basic
case: \emph{in the total
collision of the planar $N$-body problem.}

The main goal of this article is to study $PISPW$ in the total
collision of the planar $N$-body problem. To this end, we first introduce a class of moving frame suitable for  orbits near central configurations, especially for total collision orbits. 

In the coordinate system originated from the moving frame, the degeneracy of the equations of motion according to intrinsic symmetrical characteristic of $N$-body problem can easily be reduced. As a result, $PISPW$ can be well described.  In fact, once the moving frame is successfully set, one can describe the motion of collision orbit effectively, and give the practical equations of motion.

As  by-products, the moving frame and the concomitant coordinate system are found useful investigating other questions of the planar $N$-body problem. Indeed, in addition to $PISPW$, we have found that the moving frame and the concomitant coordinate system are also useful to investigate the stability  of  relative equilibrium solutions,  the degenerate central configurations and periodic orbits of the planar $N$-body problem so far.

It is noteworthy that the coordinates  from the moving frame have many remarkable differences from the well-known McGehee's coordinate system \cite{Mcgehee1974Triple}. For instance, the coordinates from the moving frame  primarily focus on the configuration space  of the $N$-body problem, but McGehee's coordinates focus on the phase space of the $N$-body problem. 

It is shown that collision orbits belong to unstable manifolds of origin with regard to a subsystems of equations. Unfortunately, results on stable manifolds and unstable manifolds  cannot apply to $PISPW$ directly.
However we find that some original ideas of Siegel  \cite{C1967Lectures} are applicable. The ideas are related to normal forms, which is especially important for us. Since the original results on normal forms in \cite{C1967Lectures} can only be applied  to the case corresponding to nondegenerate central configurations, it is necessary to generalize the results of normal forms in \cite{C1967Lectures}  to the case corresponding to degenerate central configurations. Thus the theory of  central manifolds is also introduced to explore the case corresponding to degenerate central configurations.

First we give a full answer to $PISPW$ in the case corresponding to nondegenerate central configurations: the normalized configuration of the particles must approach a certain central configuration without undergoing infinite spin for collision orbits. This result is an  immediate  application of the theory of  hyperbolic dynamics, or equivalently, the theory of normal forms.   Therefore, as a separate method, we give a new rigorous and simple proof of the above result in this paper.

However, for $PISPW$ in the case corresponding to degenerate central configurations, the problem is  unexpectedly difficult.

It is well known that central configurations are connected with significant dynamic phenomenons of the $N$-body problem, for instance,  central configurations are important for finding critical points of integral manifolds of the $N$-body problem. Similarly, the work of this paper further shows that degenerate central configurations are connected with important  nonhyperbolic dynamic phenomenons of the $N$-body problem.
Indeed, degenerate central configurations  will lead to a degenerate system (i.e., there are  zeros in the eigenvalues of linear part of the system at an equilibrium point). Theoretically, degenerate systems are an especially difficult class of nonhyperbolic systems.

It may be an effective strategy that ones make an intensive study on central configurations to reduce the theoretical difficulties from nonhyperbolic dynamic as far as possible. Unfortunately, the problem of central configurations is also very difficult,  as indicated by the researches on this topic in the past decades. Up to now, only central configurations with special conditions are discussed.

Therefore,  in the paper, we mainly  study $PISPW$ corresponding to  central configurations with two or less degrees of degeneracy. We completely solve the problem  in the case corresponding to  central configurations with zero or one degree of degeneracy: the configuration of the particles must approach a certain central configuration without undergoing infinite spin. Combining some results on the planar nonhyperbolic equilibrium point, we give a criterion for the case corresponding to  central configurations with  degree of degeneracy two.

Because it can be shown that the exceptional masses corresponding to degenerate central configurations form a proper algebraic subset of the mass space for the  four-body problem \cite{moeckel2014lectures}, we conclude that,  for almost every choice
of the masses of the four-body problem,  the configuration of the particles must approach a certain central configuration without undergoing infinite spin for collision orbits.
Furthermore, based upon our investigation of  a kind of  symmetrical degenerate central configurations of four-body,  we  answer the problem  in the case corresponding to all  known degenerate central configurations of four-body.

After $PISPW$ is investigated, we naturally study  the manifold of all the collision orbits (i.e., the set of initial conditions leading to total collisions). It is showed that this set is locally a finite union of  real  submanifold in the neighbourhood of the collision instant, and  the
dimensions of the submanifolds depend upon the index of the limiting central configuration (i.e., the number of positive eigenvalues of the limiting  central configuration).

Finally, we examine the question of whether orbits can be extended through total
collision from the viewpoint of Sundman and Siegel, that is, whether a single
solution can be extended as an analytic function of time. We only consider the case corresponding to nondegenerate central configurations.

The paper is structured as follows. In \textbf{Section 2}, we introduce some notations,  and some preliminary results of central configurations; in particular, we introduce  a moving frame and  write equations of motion.
In \textbf{Section 3}, it will be seen that collision orbits are well described in the moving frame. In particular, $PISPW$ can be put in a form suitable   for deeper investigation. In \textbf{Section 4}, we investigate $PISPW$. In \textbf{Section 5}, we investigate the set of initial conditions leading to total collisions locally and examines the question  whether a single
solution can be extended as an analytic function of time. In \textbf{Section 6}, we summarize the main results and give some interesting  open questions. Finally, in \textbf{Appendix \ref{DegeneracyC.C.}}, we give a  criterion  for the degeneracy of central configurations by using the cartesian coordinates; in \textbf{Appendix \ref{C.C.ofFour-body}}, we investigate the  degenerate central configurations of the planar four-body problem  with an axis of symmetry; in \textbf{Appendix \ref{Lagrangiansystems}}, we sketchily give the theory of Lagrangian dynamical systems to deduce the general equations of motion; in \textbf{Appendix \ref{MotionforCollisionOrbits}}, we mainly show that the flow of the system (\ref{equationzr3}) restricted to some invariant manifold
 is gradient-like; in \textbf{Appendix \ref{NormalFormsth}}, we give the theory of normal forms (or reduction theorem) to  simplify equations of the problem in this paper; in \textbf{Appendix \ref{Plane Equilibrium Points}}, we discuss some aspects of planar equilibrium points.

As it should be,  this work indicates the fact that  $PISPW$ is  links of many disciplines. So $PISPW$ is a very good problem of great value to study.

\section{Preliminaries}
\label{Preliminaries}
\indent\par
In this section we fix notations and give some definitions, in particular, we will introduce a moving frame to  describe  effectively the orbits  of motion  near central configurations. 

\subsection{Central Configurations}
\indent\par
Let $(\mathbb{R}^2)^N $ denote the space of configurations for $N$ point particles in the Euclidean plane $\mathbb{R}^2$: $ (\mathbb{R}^2)^N = \{ \mathbf{r} = (\mathbf{r}_1,  \cdots, \mathbf{r}_N)|\mathbf{r}_j \in \mathbb{R}^2, j = 1,  \cdots, N \}$. In this paper, unless otherwise specified, the cartesian space $(\mathbb{R}^2)^N $ is considered as a column space.

For each pair of indices $j, k \in\{ 1,\ldots,N\}$, let $\Delta_{(j,k)} = \{ \mathbf{r} \in (\mathbb{R}^2)^N | \mathbf{r}_j = \mathbf{r}_k \}$ denote the collision set of the j-th and k-th particles. Let $\Delta = \bigcup_{j,k} \Delta_{(j,k)}$ be the collision set in $(\mathbb{R}^2)^N $. Then $(\mathbb{R}^2)^N  \backslash \Delta$ is the space of collision-free configurations.

The \emph{mass scalar product} in the space $(\mathbb{R}^2)^N $ is defined as:
\begin{displaymath}
\langle\mathbf{r},\mathbf{s}\rangle = \sum_{j=1}^{N} {{m_j (\mathbf{r}_j,\mathbf{s}_j)}},
\end{displaymath}
where $\mathbf{r} = (\mathbf{r}_1,  \cdots, \mathbf{r}_N)$ and $\mathbf{s} = (\mathbf{s}_1,  \cdots, \mathbf{s}_N)$ are two  configurations in $(\mathbb{R}^2)^N$, and $(\cdot,\cdot)$ denotes the standard scalar product in $\mathbb{R}^2$. Certainly,  the norm $\|\cdot\|$ for the mass scalar product can be introduced as
\begin{displaymath}
\|\mathbf{r}\|=\sqrt{\langle\mathbf{r},\mathbf{r}\rangle}.
\end{displaymath}
As a result, the cartesian space $(\mathbb{R}^2)^N $ is a new Euclidean space.

Given a configuration $\mathbf{r}$, let $\widehat{\mathbf{r}}:= \frac{\mathbf{r}}{\|\mathbf{r}\|}$  be the unit vector corresponding to $\mathbf{r}$ henceforth. In particular, the unit vector $\widehat{\mathbf{r}}$ is called the  \textbf{normalized configuration} of the configuration $\mathbf{r}$.

Let $\mathbf{r}_c = \frac{\sum_{k = 1}^{N} {m_k \mathbf{r}_k}}{\mathfrak{m}}$ be the center of mass, where $\mathfrak{m}=\sum_{k = 1}^{N} m_k$ is the total mass. Observe that the equations (\ref{eq:Newton's equation1})  of motion
are invariant by translation, so it is always  assumed that the center of mass $\mathbf{r}_c$ is at the origin.
Let $\mathcal{X}$ denote the space of configurations whose center of mass is at the origin; that is,
\begin{displaymath}
\mathcal{X} = \{ \mathbf{r} = (\mathbf{r}_1,\cdots, \mathbf{r}_N)\in (\mathbb{R}^2)^N| \sum_{k = 1}^{N} {m_k \mathbf{r}_k} = 0  \}.
\end{displaymath}
Then  $\mathcal{X}$ is a $2(N-1)$-dimensional subspace of the Euclidean space $(\mathbb{R}^2)^N $.
The open subset $\mathcal{X} \backslash \Delta$ is the space of collision-free configurations with the center of mass at the origin.\\

Let us recall the important concept of  central configurations \cite{wintner1941analytical}:
\begin{definition}
A configuration $\mathbf{r} \in \mathcal{X} \backslash \Delta$ is called a central configuration if there exists a constant $\lambda\in {\mathbb{R}}$ such that
\begin{equation}\label{centralconfiguration}
\sum_{j=1,j\neq k}^N \frac{m_jm_k}{|\mathbf{r}_j-\mathbf{r}_k|^3}(\mathbf{r}_j-\mathbf{r}_k)=-\lambda m_k\mathbf{r}_k,1\leq k\leq N.
\end{equation}
\end{definition}

The value of $\lambda$ in (\ref{centralconfiguration}) is uniquely determined by
\begin{equation}\label{lambda}
\lambda=\frac{\mathcal{U}(\mathbf{r})}{I(\mathbf{r})},
\end{equation}
where the opposite of the potential energy (force function) $\mathcal{U}$ and  the moment
of inertia $I$ are respectively defined as
\begin{displaymath}
\mathcal{U}(\mathbf{r}) = \sum_{1\leq k<j\leq N} {\frac{m_k m_j }{|\mathbf{r}_k - \mathbf{r}_j|}},
\end{displaymath}
\begin{displaymath}
I(\mathbf{r}) = \sum_{j=1}^{N} {{m_j |{\mathbf{r}}_j-{\mathbf{r}}_c|^2}}.
\end{displaymath}
Hereafter,
for given $m_j$($j=1, 2, \cdots, N$) and a fixed $\lambda$, let $\textbf{CC}_{\lambda}$ be the set of central configurations  satisfying equations (\ref{centralconfiguration}).

Central configurations are
important for the  N-body problem,  Saari \cite{zbMATH03657180} even said that central configurations play a central role in the $N$-body Problem. Indeed, in the $N$-body Problem, central configurations can produce homographic solutions; central configurations play a critical role in studying the topology of integral manifolds of energy  and angular momentum \cite{smale1970topology}; central configurations are limiting configurations in studying the asymptotic behavior of every motion starting and ending in a collision and every parabolic motion \cite{Sundman,wintner1941analytical,Sperling1970On,Saari1981On} and so on.

Although there were a lot of works on  central configurations, many significant problems of central configurations are open
up to now, in these open problems, the most famous one is the conjecture on the Finiteness of Central Configurations \cite{zbMATH01516323}:  for any given masses $m_1,\cdots,m_N$, the number of central configurations in the $N$-body Problem should be finite. Note that if $\mathbf{r} \in \mathcal{X} \backslash \Delta$ is  a central configuration,  so is $\rho A\mathbf{r}$ for any $A \in S\mathbb{O}(2)$ and $\rho>0$ via the transformation $\rho A\mathbf{r}=(\rho A\mathbf{r}_1, \rho A\mathbf{r}_2, \cdots, \rho A\mathbf{r}_N)$, where  $S\mathbb{O}(2)$ is the special orthogonal group of the plane.  The transformation $\rho A$ defines an equivalence relation of $\mathcal{X}$. Thus when counting central configurations, we actually count the equivalence classes of them.

\indent\par

There are several  equivalent definitions  of  central configurations. One of the equivalent definitions considers  a central configuration  as a critical point of the function  $\widetilde{\mathcal{U}}:=I^{\frac{1}{2}}\mathcal{U}$ refereed as \emph{normalized potential}, since
\begin{displaymath}
\nabla\widetilde{\mathcal{U}} (\mathbf{r})= I^{-\frac{1}{2}} \mathcal{U}  \mathbf{r} + I^{\frac{1}{2}} \nabla \mathcal{U},
\end{displaymath}
and the equations (\ref{centralconfiguration}) are just
\begin{equation}\label{centralconfiguration1}
 \nabla \mathcal{U} = -\lambda{\mathbf{r}},
\end{equation}
where $\nabla f$ is the gradient of a differential function $f$ on $\mathcal{X}$ (or $\mathcal{X} \backslash \Delta$)  with respect to the mass scalar product, i.e., given $\mathbf{r}\in \mathcal{X},$
\begin{equation*}
df({\mathbf{r}})(\mathbf{v})=\langle\nabla f({\mathbf{r}}),\mathbf{v}\rangle, ~~~~~~~~~~~~~~~~~~~~ \forall \mathbf{v}\in \mathrm{T}_{\mathbf{r}}\mathcal{X},
\end{equation*}
where $\mathrm{T}_{\mathbf{r}}\mathcal{X}$ is the tangent space of $\mathcal{X}$ at the point $\mathbf{r}$.

It is well known that the critical points of  $\widetilde{\mathcal{U}}$ are not isolated but rather occur as manifolds of critical points. Thus these critical points are always degenerate in the ordinary sence. More specifically,
let  $\mathbf{r} \in \mathcal{X} \backslash \Delta$ be  a central configuration, i.e., a critical point $\mathbf{r}$ of $\widetilde{\mathcal{U}}$.  By the invariance of  $\widetilde{\mathcal{U}}$ under  the action of the transformation $\rho A$, it follows that the manifold $\{\rho A\mathbf{r}|\rho>0,A \in S\mathbb{O}(2)\}$ is consisting of critical points  of $\widetilde{\mathcal{U}}$;
furthermore, it follows that the Hessian of $\widetilde{\mathcal{U}}$ evaluated at a central configuration $\mathbf{r}$ must
contain
\begin{center}
 $\mathcal{P}_{\mathbf{r}}:=\{\rho A\mathbf{r}|\rho>0,A \in S\mathbb{O}(2)\}\bigcup\{\mathbf{0}\}$
\end{center}
in its kernel. One naturally removes the trivial degeneracy according to this intrinsic symmetrical characteristic of the $N$-body problem.  Taking into account these facts, a central configuration $\mathbf{r}$ will be called \textbf{nondegenerate}, if  the kernel of the Hessian of $\widetilde{\mathcal{U}}$ evaluated at $\mathbf{r}$ is exactly $\mathcal{P}_{\mathbf{r}}$. Obviously, this definition of nondegeneracy is
equivalent to the one used by Palmore in his study of planar central
configurations \cite{Palmore1976Measure}.

If  all the central configurations are nondegenerate  for any choice of positive masses, then  the famous conjecture on the Finiteness of Central Configurations is correct. However,  it has been known that  degenerate central configurations exist in the  $N$-body problem for any $N>3$. Thus one naturally conjectures that all the central configurations are nondegenerate for almost every choice of positive masses. Unfortunately,  no practical progress has been made for this open problem so far, except for $N=4$.
For the four-body problem it has been shown that the exceptional masses corresponding to degenerate central configurations form a proper algebraic subset of the mass space by Moeckel \cite{moeckel2014lectures}.

In this paper,  the following concept of  \textbf{degrees of degeneracy}  for central configurations is important.
 \begin{definition}
Given a central configuration $\mathbf{r}\in \mathcal{X} \backslash \Delta$, consider the Hessian $D^2\widetilde{\mathcal{U}}(\mathbf{r})$ of $\widetilde{\mathcal{U}}$  at the point $\mathbf{r}$, i.e., the linear operator $D^2\widetilde{\mathcal{U}}(\mathbf{r}):\mathrm{T}_{\mathbf{r}}\mathcal{X}\rightarrow \mathrm{T}_{\mathbf{r}}\mathcal{X}$ characterized by
\begin{equation*}
    d^2\widetilde{\mathcal{U}}(\mathbf{r})(\mathbf{u},\mathbf{v})=\langle D^2\widetilde{\mathcal{U}}(\mathbf{r})\mathbf{u},\mathbf{v}\rangle,~~~~~~~~~~~~~~~~~~~~ \forall \mathbf{u},\mathbf{v}\in \mathrm{T}_{\mathbf{r}}\mathcal{X}.
\end{equation*}
Assume the nullity (i.e., the dimension of the kernel) of the Hessian $D^2\widetilde{\mathcal{U}}(\mathbf{r})$ is $n_0+2$,
then the number $n_0$   is called  the degree of degeneracy of  $\mathbf{r}$. In particular, a central configuration  with $0$ degree of degeneracy is nondegenerate.
\end{definition}

Given a central configuration $\mathbf{r}_0\in \mathcal{X} \backslash \Delta$, let $\mathcal{P}^{\bot}_{\mathbf{r}_0}$ be the orthogonal complement of $\mathcal{P}_{\mathbf{r}_0}$ in $\mathcal{X}$, i.e.,
\begin{equation*}
    \mathcal{X}=\mathcal{P}_{\mathbf{r}_0}\oplus\mathcal{P}^{\bot}_{\mathbf{r}_0}.
\end{equation*}
By the facts that  the Hessian $D^2\widetilde{\mathcal{U}}(\mathbf{r}_0)$ is  a symmetric linear operator and $\mathcal{P}_{\mathbf{r}_0}$ is an invariant subspace of  $D^2\widetilde{\mathcal{U}}(\mathbf{r}_0)$, it follows that $\mathcal{P}^{\bot}_{\mathbf{r}_0}$  is also an invariant subspace of  $D^2\widetilde{\mathcal{U}}(\mathbf{r}_0)$. Therefore,  $D^2\widetilde{\mathcal{U}}(\mathbf{r}_0)$ can be diagonalized in an orthogonal basis of $\mathcal{X}$ consisting of eigenvectors of $D^2\widetilde{\mathcal{U}}(\mathbf{r}_0)$, and the first two vectors of the
basis can  be chosen as any two mutually orthogonal  central configurations in $\mathcal{P}_{\mathbf{r}_0}$.

It is noteworthy that both of the subspaces $\mathcal{P}_{\mathbf{r}_0}$ and $\mathcal{P}^{\bot}_{\mathbf{r}_0}$ are   invariant under the action of the transformation $\rho A$.

It is necessary to parameterize $S\mathbb{O}(2)$ in the following, but for the ease of   notations, we simply  set
\begin{center}
$S\mathbb{O}(2)=\{e^{\mathbf{i}\theta}|\theta\in \mathbb{R}\}$,
\end{center}
where $\mathbf{i}$ is the imaginary unit. In particular, $\mathbf{i}\mathbf{r}_0$ is just $\mathbf{r}_0$ rotated by an angle $\frac{\pi}{2}$. Thus, when necessary, one may identify $\mathbb{R}^2$ with $\mathbb{C}$ and $(\mathbb{R}^2)^N$ with $\mathbb{C}^N$ and so on.

Assume
$\{\mathcal{E}_3, \mathcal{E}_4, \mathcal{E}_5,\cdots, \mathcal{E}_{2N}\}$ is an orthogonal basis of $\mathcal{X}$ such that
\begin{equation}\label{eigenorthogonalbasis}
    D^2\widetilde{\mathcal{U}}(\mathbf{r}_0)\mathcal{E}_j=\mu_j\mathcal{E}_j,~~~~~~~~~~~~j=3,\cdots,2N,
\end{equation}
and
\begin{center}
$\mathcal{E}_3=\mathbf{r}_0$, $\mathcal{E}_4=\mathbf{i}\mathbf{r}_0$.
\end{center}
Then $\mathcal{P}_{\mathbf{r}_0}=span\{{\mathcal{E}}_3, {\mathcal{E}}_4\}$, $\mathcal{P}^{\bot}_{\mathbf{r}_0}=span\{{\mathcal{E}}_5, \cdots, {\mathcal{E}}_{2N}\}$ and $\mu_3=\mu_4=0$.
It is noteworthy that the values of $\mu_j$ ($j=5, \cdots,2N$) depend on (the scale $\|\mathcal{E}_{3}\|$ of) the central configuration $\mathcal{E}_3=\mathbf{r}_0$. Indeed, the Hessian $D^2\widetilde{\mathcal{U}}(\mathbf{r}_0)$ is minus quadratic about $\mathbf{r}_0$, i.e.,
\begin{equation*}
    D^2\widetilde{\mathcal{U}}(\rho\mathbf{r}_0)=\frac{1}{\rho^2} D^2\widetilde{\mathcal{U}}(\mathbf{r}_0),  ~~~~~~~~~~~~\forall \rho>0.
\end{equation*}
Furthermore, by the invariance of  $\widetilde{\mathcal{U}}$ under  the action of the transformation $\rho A$, it follows that
\begin{equation}\label{Hessiansimilarly}
    D^2\widetilde{\mathcal{U}}(\rho A\mathbf{r}_0)=\frac{1}{\rho^2}diag(A,\cdots,A) D^2\widetilde{\mathcal{U}}(\mathbf{r}_0)diag(A^{-1},\cdots,A^{-1}),  ~~~\forall \rho>0,\forall A \in S\mathbb{O}(2).
\end{equation}

Assume that there are $n_0$ zeros and $n_p$ positive numbers in $\mu_j$ ($j=5, 6, \cdots, 2N$), then $n_0$ is just  the degrees of degeneracy of the central configuration $\mathcal{E}_3$.  By (\ref{Hessiansimilarly}), it follows that $n_0$ and $n_p$ are  invariant under  the action of the transformation $\rho A$, i.e.,
\begin{equation*}
    n_0(\rho A\mathbf{r}_0)=n_0(\mathbf{r}_0),~~~~~~~~~~ n_p(\rho A\mathbf{r}_0)=n_p(\mathbf{r}_0).
\end{equation*}
Note that, by a classic result (see \cite{Palmore1976Measure,moeckel1990central}), it follows that $n_p \geq N-2$.

The degeneracy of central configurations  can also be described in a specific coordinate system of $(\mathbb{R}^2)^N $, in particular, this method offers some convenience at practical calculations of  degeneracy of central configurations, for more detail please refer to Appendix \ref{DegeneracyC.C.} and Appendix \ref{C.C.ofFour-body}.

Obviously, $\{\widehat{\mathcal{E}}_3, \widehat{\mathcal{E}}_4, \cdots, \widehat{\mathcal{E}}_{2N}\}$ is an orthonormal basis of $\mathcal{X}$, in this orthonormal basis, every configuration $\mathbf{r}\in \mathcal{X}$ can be written as $\mathbf{r} =  \sum_{j = 3}^{2N} {y^j \widehat{\mathcal{E}}_j}$ and $I(\mathbf{r})=\sum_{j = 3}^{2N} (y^j)^{2}$. It is remarkable that the equations (\ref{eq:Newton's equation1}) of motion can be expressed clearly in the $y^j$ coordinates, this is especially useful when we study  relative equilibrium solutions of the Newtonian $N$-body problem. However, we will adopt another coordinate system originating from a kind of moving frame, which is more suitable for collision orbits and relative equilibrium solutions.

In the next two subsections, we first introduce a moving frame to  describe effectively the  orbits near central configurations. Here we say a  configuration $\mathbf{r}$ is near the central configuration $\mathbf{r}_0$, if $\langle\mathbf{r},\mathbf{r}_0\rangle\neq 0$ or $\langle\mathbf{r},\mathbf{i}\mathbf{r}_0\rangle\neq 0$. Then we deduce the general equations of motion in the concomitant system of coordinates.\\

\subsection{Moving Frame}
\indent\par

The notations in this subsection  are inherited from the previous subsection.

For any configuration $\mathbf{r}\in \mathcal{X} \backslash  \mathcal{P}^{\bot}_{\mathbf{r}_0}$, it is easy to see that  there exists a unique point $e^{\mathbf{i}\theta(\mathbf{r})}\widehat{\mathcal{E}}_3$ on $\mathbf{S}$ such that
 \begin{equation}
\|e^{\mathbf{i}\theta(\mathbf{r})}\widehat{\mathcal{E}}_3-\mathbf{r}\|= min_{\theta \in \mathbb{R}} \|e^{\mathbf{i}\theta}\widehat{\mathcal{E}}_3 - \mathbf{r}\|,\nonumber
\end{equation}
where $\mathbf{S}= \{e^{\mathbf{i}\theta}\widehat{\mathcal{E}}_3| \theta \in \mathbb{R}\}$ is a unit circle in the space $(\mathbb{R}^2)^N $ with the origin as the center.

Note that the variable $\theta$ in  $\theta(\mathbf{r})$  can be continuously determined as a continuous function of the independent variable $\mathbf{r}$. Indeed, according to the relation
\begin{equation}
\frac{d\langle e^{\mathbf{i}\theta}\widehat{\mathcal{E}}_3 - \mathbf{r},e^{\mathbf{i}\theta}\widehat{\mathcal{E}}_3 - \mathbf{r}\rangle}{d\theta}= \langle e^{\mathbf{i}\theta}\widehat{\mathcal{E}}_3 - \mathbf{r},\mathbf{i}e^{\mathbf{i}\theta}\widehat{\mathcal{E}}_3 \rangle,\nonumber
\end{equation}
it follows that
$\theta(\mathbf{r})$ can  be determined by the following relations
 \begin{equation}\label{varepsilon4}
\langle\mathbf{r},\mathbf{i}e^{\mathbf{i}\theta(\mathbf{r})}\widehat{\mathcal{E}}_3 \rangle = 0
\end{equation}
\begin{equation}\label{varepsilon3}
\langle\mathbf{r},e^{\mathbf{i}\theta(\mathbf{r})}\widehat{\mathcal{E}}_3 \rangle  > 0.
\end{equation}

Set $\Xi_3 = e^{\mathbf{i}\theta(\mathbf{r})}\widehat{\mathcal{E}}_3, \Xi_4 = e^{\mathbf{i}\theta(\mathbf{r})}\widehat{\mathcal{E}}_4, \cdots, \Xi_{2N} = e^{\mathbf{i}\theta(\mathbf{r})}\widehat{\mathcal{E}}_{2N}$, then $\{\Xi_3, \Xi_{4}, \cdots, \Xi_{2N}\}$ is an orthonormal basis of $\mathcal{X}$, and
\begin{displaymath}
span\{\Xi_3, \Xi_4\}=span\{\widehat{\mathcal{E}}_3, \widehat{\mathcal{E}}_4\}=\mathcal{P}_{\mathbf{r}_0},~~~~~~span\{\Xi_5, \cdots, \Xi_{2N}\}=span\{\widehat{\mathcal{E}}_5, \cdots, \widehat{\mathcal{E}}_{2N}\}=\mathcal{P}^{\bot}_{\mathbf{r}_0}.
\end{displaymath}

The  set $\{\Xi_3, \Xi_{4}, \cdots, \Xi_{2N}\}$  is the \textbf{moving frame} for us.

Set $r = \|\mathbf{r}\|$, then $\mathbf{r}=r \widehat{\mathbf{r}}$. In the moving frame, $\widehat{\mathbf{r}}$ can be written as $\widehat{\mathbf{r}} =  \sum_{k = 3}^{2N} {z_k \Xi_k}$.
By (\ref{varepsilon4}), it is easy to see that $z_4=0$. It follows from $\|\widehat{\mathbf{r}}\| = 1$ and (\ref{varepsilon3}) that
\begin{equation}\label{z3}
z_3 = \sqrt{1 - \sum_{j = 5}^{2N} z^2_j}.
\end{equation}
Then the total set of the variables $r, \theta, z_5, \cdots, z_{2N}$ can be thought as the  coordinates of $\mathbf{r}\in \mathcal{X} \backslash  \mathcal{P}^{\bot}_{\mathbf{r}_0}$ in the moving frame. Indeed, we have defined a real analytic diffeomorphism:
\begin{displaymath}
(0, +\infty) \times \mathbf{S}^1 \times \mathbf{B}^{2N-4} \rightarrow \mathcal{X} \backslash  \mathcal{P}^{\bot}_{\mathbf{r}_0}:
\end{displaymath}
\begin{displaymath}
(r, (\cos\theta,\sin\theta), z_5, \cdots, z_{2N}) \mapsto r(z_3 \Xi_3 + \sum_{j = 5}^{2N} {z_j \Xi_j}),
\end{displaymath}
and a classic covering map:
\begin{displaymath}
\mathbb{R} \rightarrow \mathbf{S}^1: \theta \mapsto (\cos\theta,\sin\theta),
\end{displaymath}
where $\mathbf{S}^1$ is the unit circle in the plane $\mathbb{R}^2$ and $\mathbf{B}^{2N-4}=\{(z_5, \cdots, z_{2N})|\sum_{j = 5}^{2N} {z^2_j }<1\}$ is the $(2N-4)$-dimensional unit ball in  $\mathbb{R}^{2N-4}$.

\input{Template1.TpX}
Geometrically, to make the visual understanding of the system of coordinates $r, \theta, z_5, \cdots, z_{2N}$,  please see Figure \ref{coordinates}. Note that the $z$-axis in Figure \ref{coordinates} represents the space $\mathcal{P}^{\bot}_{\mathbf{r}_0}$.

\subsection{General Equations of Motion}
\indent\par

{Consider  the kinetic energy, the total energy, the angular momentum, and the Lagrangian  function, respectively, defined by}
\begin{displaymath}
\mathcal{K} (\dot{\mathbf{r}})= \sum_{j=1}^{N} {\frac{1}{2}{m_j |\dot{\mathbf{r}}_j|^2}},
\end{displaymath}
\begin{displaymath}
\mathcal{H }(\mathbf{r},\dot{\mathbf{r}})= \mathcal{K}(\dot{\mathbf{r}})- \mathcal{U}(\mathbf{r}),
\end{displaymath}
\begin{displaymath}
\mathcal{J}(\mathbf{r}) = \sum_{j=1}^{N} {m_j {\mathbf{r}}_j \times {\dot{\mathbf{r}}}_j},
\end{displaymath}
\begin{displaymath}
\mathcal{L}(\mathbf{r},\dot{\mathbf{r}}) = \mathcal{L} = \mathcal{K} + \mathcal{U} = \sum_j \frac{1}{2} m_j |\dot{\mathbf{r}}_j|^2  + \sum_{k<j}{\frac{m_k m_j}{|\mathbf{r}_k - \mathbf{r}_j|}},
\end{displaymath}
where $\times$ denotes the standard cross product in $\mathbb{R}^2$.

It is well known that the equations (\ref{eq:Newton's equation1}) of motion are just the Euler-Lagrange equations
of  Lagrangian
system with the configuration manifold $\mathcal{X}$ and the Lagrangian function $\mathcal{L}(\mathbf{r},\dot{\mathbf{r}})$. Then this yields a quick method
for writing equations of motion in various coordinate systems. In fact, to write the equations of motion in a new coordinate system, it is sufficient to
express the Lagrangian function in the new coordinates. Please refer to \ref{LagrangianDynamicalsystems} for more detail.

Thus to write the general equations of motion in the above system of coordinates
\begin{center}
$r, \theta, z_5, \cdots, z_{2N}$,
\end{center}
which is well defined in the submanifold $\mathcal{X}\backslash \mathcal{P}^{\bot}_{\mathbf{r}_0}$,
 it suffices to rewrite
 the kinetic energy and the force function  as
\begin{equation}\label{kineticenergy}
\mathcal{K}(\mathbf{r}) = \frac{\dot{r}^2}{2}+ \frac{r^2}{2} (\dot{z}^2_3+\sum_{j=5}^{2N}\dot{z}^2_j + 2\dot{\theta}\sum_{j,k=5}^{2N}\langle\widehat{\mathcal{E}}_j,\mathbf{i}\widehat{\mathcal{E}}_k\rangle\dot{z}_j {z}_k+ \dot{\theta}^2),\nonumber
\end{equation}
\begin{equation}\label{force function}
\mathcal{U}(\mathbf{r}) = \frac{\mathcal{U}(z_3 \widehat{\mathcal{E}}_3+\sum_{j = 5}^{2N} {z_j \widehat{\mathcal{E}}_j})}{r}.\nonumber
\end{equation}
Since $\mathcal{U}(z_3 \widehat{\mathcal{E}}_3+\sum_{j = 5}^{2N} {z_j \widehat{\mathcal{E}}_j})$ only contains the variables $z_j$ ($j=5, \cdots, 2N$), we will simply
write it as  $U(z)$ henceforth and we always think that $z= (z_5, \cdots, z_{2N})^\top$.  Set $q_{jk}=\langle\widehat{\mathcal{E}}_j,\mathbf{i}\widehat{\mathcal{E}}_k\rangle$, then the square matrix $Q:=(q_{jk})_{(2N-4)\times (2N-4)}$ is an anti-symmetric orthogonal matrix.
As a result, it follows that the Lagrangian  function $\mathcal{L}$ is
\begin{equation}\label{Lagrangianfunction1}
    \mathcal{L}(z,r,\dot{z},\dot{r},\dot{\theta})=\frac{\dot{r}^2}{2}+ \frac{r^2}{2} (\dot{z}^2_3+\sum_{j=5}^{2N}\dot{z}^2_j + 2\dot{\theta}\sum_{j,k=5}^{2N}q_{jk}\dot{z}_j {z}_k+ \dot{\theta}^2)+\frac{U(z)}{r}.
\end{equation}
The equations of motion in the coordinates $ r, \theta, z_5, \cdots, z_{2N}$ are just the Euler-Lagrange equations
of (\ref{Lagrangianfunction1}):
\begin{equation}\label{Euler-Lagrangeequations0}
\left\{
             \begin{array}{l}
           \frac{d}{dt}\frac{\partial \mathcal{L}}{\partial \dot{z}_i}-\frac{\partial \mathcal{L}}{\partial z_i}=0,~~~~~~~~~~~~~~~~~~i=5,\cdots, 2N;    \\
             \frac{d}{dt}\frac{\partial \mathcal{L}}{\partial \dot{r}}-\frac{\partial \mathcal{L}}{\partial r}=0,   \\
            \frac{d}{dt}\frac{\partial \mathcal{L}}{\partial \dot{\theta}}-\frac{\partial \mathcal{L}}{\partial \theta}=0.
             \end{array}
\right.
\end{equation}
We remark that the degeneracy of $z_3, z_4$ according to intrinsic symmetrical characteristic  of the $N$-body problem (i.e., the Newton equations (\ref{eq:Newton's equation1}) are invariant under the transformation $(\mathbf{r},t)\mapsto (\rho A \mathbf{r}, \rho^{\frac{3}{2}}t )$) has  been reduced  in the coordinates $r, \theta, z_5, \cdots, z_{2N}$.

It is noteworthy that the variable $\theta$ is not involved in the function $U(z)$, this is a main reason of introducing the moving frame.  In particular, the variable $\theta$ is not involved in the Lagrangian $\mathcal{L}$, that is, the variable $\theta$ is an ignorable coordinate. Then $\frac{\partial \mathcal{L}}{\partial \dot{\theta}}$ is conserved. In fact,   a straightforward computation shows that
\begin{equation}\label{angular momentum}
\frac{\partial \mathcal{L}}{\partial \dot{\theta}}= \mathcal{J} = \langle{\mathbf{i}\mathbf{r}},{\dot{\mathbf{r}}}\rangle
= r^2(\dot{\theta}+\sum_{j,k=5}^{2N}q_{jk}\dot{z}_j {z}_k).
\end{equation}

By \emph{Routh's method for eliminating
ignorable coordinates}, we introduce the function
\begin{equation}\label{reducedLagrangian}
\mathcal{L}_{\mathcal{J}}(z,r,\dot{z},\dot{r})=\mathcal{L}(z,r,\dot{z},\dot{r},\dot{\theta}) - {\mathcal{J}} \dot{\theta}|_{r,z,\dot{r},\dot{z},{\mathcal{J}}} \nonumber
\end{equation}
as the reduced Lagrangian function on the level set of $\mathcal{J}$, where $\mathcal{J}$ is certainly a constant and we represent
$\dot{\theta}$ as a function of $z,r,\dot{z},\dot{r}$ and $\mathcal{J}$ by using the equality (\ref{angular momentum}). A straightforward computation shows that
\begin{equation}\label{reducedLagrangian1}
\mathcal{L}_{\mathcal{J}}(z,r,\dot{z},\dot{r})= \frac{\dot{r}^2}{2}+ \frac{r^2}{2} \left[\dot{z}^2_3+\sum_{j=5}^{2N}\dot{z}^2_j -\left(\sum_{j,k=5}^{2N}q_{jk}\dot{z}_j {z}_k- \frac{\mathcal{J}}{r^2}\right)^2\right]+\frac{U(z)}{r}.
\end{equation}
In particular, if the angular momentum $\mathcal{J}$ is zero, then the reduced Lagrangian function becomes
\begin{equation}\label{reducedLagrangian2}
\mathcal{L}_{0}(z,r,\dot{z},\dot{r})= \frac{\dot{r}^2}{2}+ \frac{r^2}{2} \left[\dot{z}^2_3+\sum_{j=5}^{2N}\dot{z}^2_j -\left(\sum_{j,k=5}^{2N}q_{jk}\dot{z}_j {z}_k\right)^2\right]+\frac{U(z)}{r}.
\end{equation}
The equations of motion on the level set of $\mathcal{J}$ are just the Euler-Lagrange equations
of (\ref{reducedLagrangian1}):
\begin{equation}\label{Euler-Lagrangeequations1}
\left\{
             \begin{array}{l}
           \frac{d}{dt}\frac{\partial \mathcal{L}_{\mathcal{J}}}{\partial \dot{z}_i}-\frac{\partial \mathcal{L}_{\mathcal{J}}}{\partial z_i}=0,~~~~~~~~~~~~~~~~~~i=5,\cdots, 2N;    \\
             \frac{d}{dt}\frac{\partial \mathcal{L}_{\mathcal{J}}}{\partial \dot{r}}-\frac{\partial \mathcal{L}_{\mathcal{J}}}{\partial r}=0.
             \end{array}
\right.
\end{equation}
Especially, we give the equations of motion on the level set  $\mathcal{J}=0$ in detail here. For more detail of (\ref{Euler-Lagrangeequations0}) and (\ref{Euler-Lagrangeequations1}) please refer to Appendix \ref{EquationsMotionMovingCoordinateSystem}. On the level set  $\mathcal{J}=0$, the equations of motion reduce to
\begin{equation}\label{Euler-Lagrangeequations1detai2}
\left\{
             \begin{array}{l}
              \frac{d}{dt}\frac{\partial K(z,\dot{z})}{\partial \dot{z}_i}-\frac{\partial K(z,\dot{z})}{\partial z_i}
 +\frac{2\dot{r}}{r} \frac{\partial K(z,\dot{z})}{\partial \dot{z}_i} - \frac{1}{r^3}\frac{\partial U(z)}{\partial z_i}=0,~~~~~~~~~~~~~~~~~~i=5,\cdots, 2N; \\\\
             \ddot{r}- {2r} K(z,\dot{z})+ \frac{ U(z)}{r^2}=0.
             \end{array}
\right.
\end{equation}
where
\begin{equation*}
   \begin{array}{l}
      K(z,\dot{z})=\frac{1}{2} [\dot{z}^2_3+\sum_{j=5}^{2N}\dot{z}^2_j -(\sum_{j,k=5}^{2N}q_{jk}\dot{z}_j {z}_k)^2]=\frac{1}{2} \left[\frac{ (\sum_{j=5}^{2N}{z}_j \dot{z}_j)^2}{1 - \sum_{j = 5}^{2N} z^2_j}+\sum_{j=5}^{2N}\dot{z}^2_j - (\sum_{j,k=5}^{2N}q_{jk}\dot{z}_j {z}_k)^2\right],
   \end{array}
\end{equation*}
\begin{equation*}
  \begin{array}{l}
    \frac{\partial K(z,\dot{z})}{\partial \dot{z}_i}=\frac{z_i \sum_{j=5}^{2N}\dot{z}_j z_j}{1 - \sum_{j = 5}^{2N} z^2_j} + \dot{z}_i -  \sum_{j=5}^{2N}q_{ij} z_j\sum_{j,k=5}^{2N}q_{jk}\dot{z}_j {z}_k,
  \end{array}
\end{equation*}
and
\begin{equation*}
   \begin{array}{l}
      \frac{d}{dt}\frac{\partial K(z,\dot{z})}{\partial \dot{z}_i}-\frac{\partial K(z,\dot{z})}{\partial z_i} \\
      =\ddot{z}_i - \sum_{j=5}^{2N}q_{ij} z_j\sum_{j,k=5}^{2N}q_{jk}\ddot{z}_j {z}_k- 2  \sum_{j=5}^{2N}q_{ij} \dot{z}_j\sum_{j,k=5}^{2N}q_{jk}\dot{z}_j {z}_k
               +\frac{z_i \sum_{j=5}^{2N}(\ddot{z}_j z_j + \dot{z}^2_j)}{1 - \sum_{j = 5}^{2N} z^2_j} + \frac{{z}_i (\sum_{j=5}^{2N}{z}_j \dot{z}_j)^2}{(1 - \sum_{j = 5}^{2N} z^2_j)^2}.
   \end{array}
\end{equation*}
Note that by (\ref{angular momentum}) it follows that
\begin{equation}\label{angularmomentumequation1}
   \dot{\theta}=-\sum_{j,k=5}^{2N}q_{jk}\dot{z}_j {z}_k.
\end{equation}
We will see that total collision orbits satisfy  the equations (\ref{Euler-Lagrangeequations1detai2}) and (\ref{angularmomentumequation1}) near collision instants.

Incidentally, the total energy is conserved:
\begin{equation}
    \mathcal{H}=\frac{\dot{r}^2}{2}+ r^2K(z,\dot{z})-\frac{U(z)}{r}.\nonumber
\end{equation}
By (\ref{Euler-Lagrangeequations1detai2}), it follows that
\begin{equation*}\label{totalenergyconserved1}
    \mathcal{H}=\frac{\dot{r}^2}{2}+ \frac{r \ddot{r}}{2}-\frac{U(z)}{2r}.
\end{equation*}

We conclude this subsection with some discussions of $U(z)$ near the point $z=0$. First, we can expand $U(z)$ as
\begin{equation}\label{expandingU}
   \begin{array}{l}
     U(z)=\mathcal{U}(z_3 \widehat{\mathcal{E}}_3+\sum_{j = 5}^{2N} {z_j \widehat{\mathcal{E}}_j})=\widetilde{\mathcal{U}}(z_3 \widehat{\mathcal{E}}_3+\sum_{j = 5}^{2N} {z_j \widehat{\mathcal{E}}_j}) \\\\
     =\widetilde{\mathcal{U}}(\hat{\mathcal{E}}_3) + \sum_{k = 5}^{2N} d\widetilde{\mathcal{U}}|_{\hat{\mathcal{E}}_3}(\hat{\mathcal{E}}_k)z_k + d\widetilde{\mathcal{U}}|_{\hat{\mathcal{E}}_3}(\hat{\mathcal{E}}_3)(z_3-1) \\\\
     +\frac{1}{2}[\sum_{j,k = 5}^{2N} d^2 \widetilde{\mathcal{U}}|_{\hat{\mathcal{E}}_3}(\hat{\mathcal{E}}_j,\hat{\mathcal{E}}_k) z_j z_k + 2\sum_{k =5}^{2N} d^2 \widetilde{\mathcal{U}}|_{\hat{\mathcal{E}}_3}(\hat{\mathcal{E}}_3,\hat{\mathcal{E}}_k)(z_3-1)z_k] \\\\
     +\frac{1}{3!}\sum_{i,j,k = 5}^{2N} d^3 \widetilde{\mathcal{U}}|_{\hat{\mathcal{E}}_3}(\hat{\mathcal{E}}_i,\hat{\mathcal{E}}_j,\hat{\mathcal{E}}_k)z_i z_j z_k +\cdots,
   \end{array}
\end{equation}
where ``$\cdots$" denotes   power-series in  $z_j$ ($j=5, \cdots, 2N$)    starting with quartic terms,  and $d\widetilde{\mathcal{U}}|_{\hat{\mathcal{E}}_3}$,  $d^2\widetilde{\mathcal{U}}|_{\hat{\mathcal{E}}_3}$, $ d^3\widetilde{\mathcal{U}}|_{\hat{\mathcal{E}}_3}$ denote   the differential,  second order differential, third order differential of $\widetilde{\mathcal{U}}$ at $\hat{\mathcal{E}}_3$ respectively.

Without loss of generality, assume $\widehat{\mathcal{E}}_3={\mathcal{E}}_3=\mathbf{r}_0$. Then, by  (\ref{lambda}) (\ref{centralconfiguration1}) (\ref{eigenorthogonalbasis}) (\ref{z3}) and (\ref{expandingU}), it follows  that
\begin{equation}\label{expandforcefunction}
U(z)= \lambda + \frac{1}{2} \sum_{k = 5}^{2N} \lambda_k z^2_k  +\frac{1}{6}\sum_{i,j,k = 5}^{2N} a_{ijk}z_i z_j z_k +\cdots,
\end{equation}
where  $a_{ijk}=d^3 \widetilde{\mathcal{U}}|_{\widehat{\mathcal{E}}_3}(\widehat{\mathcal{E}}_i,\widehat{\mathcal{E}}_j,\widehat{\mathcal{E}}_k)$, thus $a_{ijk}$ is symmetric with respect to the subscripts $i,j,k$. We remark that 
\begin{equation*}
    a_{ijk}=d^3 \widetilde{\mathcal{U}}|_{\widehat{\mathcal{E}}_3}(\widehat{\mathcal{E}}_i,\widehat{\mathcal{E}}_j,\widehat{\mathcal{E}}_k)=d^3 \mathcal{U}|_{\widehat{\mathcal{E}}_3}(\widehat{\mathcal{E}}_i,\widehat{\mathcal{E}}_j,\widehat{\mathcal{E}}_k)=\frac{\partial^3 U(0)}{\partial z_i \partial z_j \partial z_k}.
\end{equation*}

By (\ref{expandforcefunction}), if the Hessian of   $U(z)$ is nondegenerate, that is, the central configuration  $\widehat{\mathcal{E}}_3$ is nondegenerate or $n_0=0$, then the function $U(z)$ has exactly one critical point $z=0$  in some small neighbourhood of the point $z=0$. However, even if  the central configuration  $\widehat{\mathcal{E}}_3$ is degenerate,   the function $U(z)$ still has  exactly one critical point $z=0$  in a small neighbourhood of $z=0$, provided that  the number of central configurations in the $N$-body problem is finite for  given masses $m_1,\cdots,m_N$ (in other words, as critical points of the function $\widetilde{\mathcal{U}}$, the equivalence classes of central configurations determined by the transformation $\rho A$ are isolated).

\begin{proposition}\label{isolatedcriticalpoint}
The function $U(z)$ has  exactly one critical point $z=0$  in a small neighbourhood of $z=0$, provided that  the  central configuration $\mathbf{r}_0$  is an isolated critical point of the function $\widetilde{\mathcal{U}}$ in the sence of the equivalence classes of central configurations determined by the transformation $\rho A$.
\end{proposition}
{\bf Proof.} 
Note that $U(z)=\mathcal{U}(z_3 \widehat{\mathcal{E}}_3+\sum_{j = 5}^{2N} {z_j \widehat{\mathcal{E}}_j})=\widetilde{\mathcal{U}}(z_3 \widehat{\mathcal{E}}_3+\sum_{j = 5}^{2N} {z_j \widehat{\mathcal{E}}_j})$, thus  $U$ is the composite function of the two differential functions $\widetilde{\mathcal{U}}$ and $\pi$, where the function
\begin{center}
$\pi:\mathbf{B}^{2N-4}\rightarrow  \mathbf{S}^{2N-4}_+$
\end{center}
is a diffeomorphism between $\mathbf{B}^{2N-4}$ and  upper hemisphere $\mathbf{S}^{2N-4}_+$ in $span\{\mathcal{E}_3,\mathcal{P}^{\bot}_{\mathbf{r}_0}\}$. Therefore, a critical point of $U$  in a small neighbourhood of $z=0$ is exactly corresponding to a central configuration
in a small neighbourhood of $\widehat{\mathcal{E}}_3$ confined to the upper hemisphere $\mathbf{S}^{2N-4}_+$. It is easy to see that in such a small neighbourhood only $\widehat{\mathcal{E}}_3$  is a central configuration.

So the proposition is proved.

$~~~~~~~~~~~~~~~~~~~~~~~~~~~~~~~~~~~~~~~~~~~~~~~~~~~~~~~~~~~~~~~~~~~~~~~~~~~~~~~~~~~~~~~~~~~~~~~~~~~~~~~~~~~~~~~~~~~~~~~~~~~~~~~~~~~~~~~~~~~~\Box$\\

An argument similar to the above proof, if we note that the space $span\{\mathcal{E}_3,\mathcal{P}^{\bot}_{\mathbf{r}_0}\}$ has removed the rotation freedom of $e^{\mathbf{i} \theta}$, shows  an extension of Proposition \ref{isolatedcriticalpoint} below.
\begin{proposition}\label{isolatedcriticalpoint1}
The function $U(z)$ has  finitely isolated critical points   in $\mathbf{B}^{2N-4}$, provided that  the number of central configurations in the $N$-body Problem is finite for  given masses $m_1,\cdots,m_N$.
\end{proposition}

\subsection{Invariant Set}
\indent\par
Let us finish the section by recalling some well known notions of differential equations \cite{Geometricalmethods,Hartman,Shilnikov}.

Given a differential system
\begin{equation}\label{differential system}
\dot{q} = v(q),
\end{equation}
where $ v: \Omega \rightarrow \mathbb{R}^n$ is a continuously differentiable vector field and $\Omega $ is an open set in $\mathbb{R}^n$.
For any $p \in \Omega $, let $\phi(t,p)$ be the solution  of (\ref{differential system}) passing through $p$ at $t=0$, i.e., if $q(t)=\phi(t,p)$, then $\dot{q}(t) = v(q(t))$ and $q(0)=p$. We also call $\phi$ the flow of (\ref{differential system}) if $\phi(t,p)$ is defined
for all $t \in \mathbb{R}$ and all $p \in \Omega $. The orbit $\mathcal{O}(p) $ of (\ref{differential system}) through $p$ is defined by $\mathcal{O}(p)=\{q=\phi(t,p)|t \in \mathbb{R}\} $, the positive semiorbit through   $p$ is $\mathcal{O}^+(p)=\{q=\phi(t,p)|t\geq 0 \} $ and the
negative semiorbit through $p$ is $\mathcal{O}^-(p)=\{q=\phi(t,p)|t\leq 0 \} $.

An \textbf{equilibrium point} of (\ref{differential system}) is a point $p$  such that $v(p)=0$.
A set $\Sigma$ in $\Omega$ is called an \textbf{invariant set} of (\ref{differential system}) if $\mathcal{O}(p) \subset\Sigma$ for any $p \in \Sigma$. Any orbit $\mathcal{O}$ of (\ref{differential system})
 is obviously an invariant set of (\ref{differential system}).
A set $\Sigma$ in $\Omega$ is called positively
(negatively) invariant if $\mathcal{O}^+(p) \subset\Sigma$ ($\mathcal{O}^-(p) \subset\Sigma$) for any $p \in \Sigma$.

\begin{definition}
The positive or $\omega$-limit set of an  orbit $\mathcal{O}$ is the set
\begin{equation}
\omega(\mathcal{O})= \bigcap_{p \in \mathcal{O}} \overline{\mathcal{O}^+(p) }.\nonumber
\end{equation}
where the bar denotes closure.
Similarly, The negative  or $\alpha$-limit set of a point $p$ is the set
\begin{equation}
\alpha( \mathcal{O})= \bigcap_{p \in \mathcal{O}} \overline{\mathcal{O}^-(p)}.\nonumber
\end{equation}
\end{definition}

Recall that $\omega(\mathcal{O})$ and $\alpha(\mathcal{O})$ are invariant and closed, and if $\phi(t,p)$ for $t \geq 0$ ($t \leq 0$) is bounded, then $\omega(\mathcal{O})$ ($\alpha(\mathcal{O})$) is
nonempty compact and connected, furthermore, $\phi(t,p) \rightarrow \omega(\mathcal{O})$ ($\alpha(\mathcal{O})$) as $t \rightarrow +\infty$ ($t \rightarrow -\infty$), that is, $dist(\phi(t,p) , \omega(\mathcal{O}))\rightarrow 0$ ($dist(\phi(t,p) , \alpha(\mathcal{O}))\rightarrow 0$) as $t \rightarrow +\infty$ ($t \rightarrow -\infty$), here $dist( p,q)$ denotes the distance of $p, q \in \mathbb{R}^n$.

\begin{definition}
The stable set of a positively
 invariant set $\Sigma$  is the set
\begin{equation}
\mathcal{W}^s (\Sigma) = \{p \in \Omega| dist(\phi(t,p), \Sigma) \rightarrow 0 ~~as ~~t\rightarrow +\infty\};\nonumber
\end{equation}
The unstable set of a negatively invariant set $\Sigma$  is the set
\begin{equation}
\mathcal{W}^u (\Sigma) = \{p \in \Omega| dist(\phi(t,p), \Sigma) \rightarrow 0 ~~as ~~t\rightarrow -\infty\};\nonumber
\end{equation}
\end{definition}
In particular, in case of $\Sigma$ consisting of one equilibrium point $p_0$, we have
\begin{definition}
The stable manifold of an {equilibrium point} $p_0$ is the set
\begin{equation}
\mathcal{W}^s (p_0) = \{p| \phi(t,p)\rightarrow  p_0 ~~as ~~t\rightarrow +\infty\};\nonumber
\end{equation}
The unstable manifold of an {equilibrium point} $p_0$ is the set
\begin{equation}
\mathcal{W}^u (p_0) = \{p| \phi(t,p)\rightarrow  p_0 ~~as ~~t\rightarrow -\infty\}.\nonumber
\end{equation}
\end{definition}

In a small neighbourhood of an equilibrium point $p_0$, we can expand $v$ in a Taylor series
\begin{equation}
v(p_0+q) = \frac{\partial v(p_0)}{\partial q}q + \cdots.\nonumber
\end{equation}

Let us consider the following linearized system of the system (\ref{differential system})
\begin{equation}
\dot{q} = \frac{\partial v(p_0)}{\partial q}q.\nonumber
\end{equation}

\begin{definition} The equilibrium
point $p_0$ of (\ref{differential system}) is hyperbolic if all of the eigenvalues of $\frac{\partial v(p_0)}{\partial q}$
have nonzero real parts.  Otherwise, i.e., if at least one  of the eigenvalues of $\frac{\partial v(p_0)}{\partial q}$ is on the imaginary axis, the equilibrium
point $p_0$ is nonhyperbolic. Furthermore, if at least one  of the eigenvalues of $\frac{\partial v(p_0)}{\partial q}$ is zero, the equilibrium
point $p_0$ is degenerate.
\end{definition}

It is well known that if the equilibrium
point $p_0$ of (\ref{differential system}) is hyperbolic, then the system (\ref{differential system}) is  topologically equivalent to its linearized system in a small neighbourhood of $p_0$. Unfortunately, the system (\ref{differential system}) is hard  to understand if the equilibrium
point $p_0$ is nonhyperbolic. In particular, if the equilibrium
point $p_0$ is more degenerate, that is, more zeros in  the eigenvalues of $\frac{\partial v(p_0)}{\partial q}$, it is more difficult to understand the system (\ref{differential system}).

\section{Moving Frame and Equations of Motion}
\indent\par

\subsection{Equations of Motion for Collision Orbits}
\indent\par

In order to make preparations for resolving $PISPW$, let us write the  equations of motion for collision orbits in a form as simple as possible.

Let us recall that,  an orbit $\mathbf{r}(t)$ of the $N$-body problem arrives at a \textbf{total collision} at some instant $t_0$
if and only if  $\mathbf{r}(t) \rightarrow 0$ as $t \rightarrow t_0$, that is to say,  $\|\mathbf{r}(t)\| \rightarrow 0$ as $t \rightarrow t_0$. Without loss of generality, assume the instant $t_0=0$ and we consider only that $\mathbf{r}(t) \rightarrow 0$ as $t \rightarrow 0+$ henceforth.
Then it is easy to see that a motion $\mathbf{r}(t)$  is a total collision orbit if and only if   $r \rightarrow 0, z_5\rightarrow 0, \cdots, z_{2N}\rightarrow 0$ as $t \rightarrow 0+$.

Some classical results concerning the total collision orbits can be found in \cite{wintner1941analytical}. We summarize the results as follows.
\begin{theorem}\label{asymptotic}
Suppose a solution $\mathbf{r}(t)$ of the $N$-body problem arrives at a total collision at the instant $0$, then there exists a constant $\kappa > 0$, such that
\begin{itemize}
  \item $I(\mathbf{r}(t)) \sim (\frac{3}{2})^{\frac{4}{3}} \kappa^{\frac{2}{3}} t^{\frac{4}{3}}, \dot{I}(\mathbf{r}(t)) \sim (12)^{\frac{1}{3}} \kappa^{\frac{2}{3}} t^{\frac{1}{3}}, \ddot{I}(\mathbf{r}(t)) \sim (\frac{2}{3})^{\frac{2}{3}} \kappa^{\frac{2}{3}} t^{-\frac{2}{3}}$ as $t \rightarrow 0+$.
  \item $\mathcal{U}(\mathbf{r}(t))\sim  (\frac{1}{18})^{\frac{1}{3}} \kappa^{\frac{2}{3}} t^{-\frac{2}{3}}, \mathcal{K}(\dot{\mathbf{r}}(t))\sim (\frac{1}{18})^{\frac{1}{3}} \kappa^{\frac{2}{3}} t^{-\frac{2}{3}}$ as $t \rightarrow 0+$.
  \item $\widehat{\mathbf{r}}(t) \rightarrow \textbf{CC}_{\lambda}$ as $t \rightarrow 0+$, where ${\lambda}=\frac{\kappa}{2}$.
  \item $\mathcal{J}(\mathbf{r}(t)) \equiv 0$.
\end{itemize}
\end{theorem}

Therefore, it is natural to ask whether there exists a certain central configuration $\mathbf{s}_0\in \textbf{CC}_{\lambda}$ such that
\begin{equation*}
   \widehat{\mathbf{r}}(t) \rightarrow \mathbf{s}_0,  ~~~~~~~~~~~~~~~~~~as ~~~t \rightarrow 0+.
\end{equation*}
The point is that:
\begin{enumerate}
  \item if  the number of central configurations  is infinite for  given masses $m_1,\cdots,m_N$, then it would be
possible that the  normalized configuration $\widehat{\mathbf{r}}(t)$ comes closer and closer to more than
one central configuration, which are not in the same equivalence classes, in such a way as to oscillate between these central configurations;
  \item if  the number of central configurations  is finite for  given masses $m_1,\cdots,m_N$, then it would be
possible that the  normalized configuration $\widehat{\mathbf{r}}(t)$ moves in spirals without asymptotes (or in other words,
 make
an infinite number of revolutions) as $t \rightarrow 0+$. This is the so-called \emph{ problem of infinite spin}.
\end{enumerate}
For more detail please refer to \cite[p.282--p.283]{wintner1941analytical}.

To investigate  the problem of infinite spin or $Painlev\acute{e}$-$Wintner$ (abbreviated ``$PISPW$"), from now on, unless otherwise stated, the conjecture on the Finiteness of Central Configurations is presumed to be correct: \emph{the number of central configurations  is finite for any given masses $m_1,\cdots,m_N$}.

Then it follows from Theorem \ref{asymptotic} that there is a central configuration  $\widehat{\mathcal{E}}_3(=\mathbf{r}_0) \in \textbf{CC}_{\lambda}$ such that $\widehat{\mathbf{r}}(t) \rightarrow \mathbf{S}$ $(= \{e^{\mathbf{i}\theta}\widehat{\mathcal{E}}_3| \theta \in \mathbb{R}\})$ as $t \rightarrow 0+$.  As a result, it is easy to see that $PISPW$ explores \emph{whether there exists a fixed central configuration $e^{\mathbf{i}\theta_0}\widehat{\mathcal{E}}_3 \in \mathbf{S}$ such that $\widehat{\mathbf{r}}(t) \rightarrow e^{\mathbf{i}\theta_0}\widehat{\mathcal{E}}_3$ as $t \rightarrow 0+$}.

Since $\widehat{\mathbf{r}}(t) \rightarrow \mathbf{S}$,  we have legitimate rights to use the coordinates $r, \theta, z_5, \cdots, z_{2N}$ for the total collision orbit $\mathbf{r}(t)$.  The coordinates $r(t), \theta(t), z_5(t), \cdots, z_{2N}(t)$ of the total collision orbit $\mathbf{r}(t)$ are real analytic functions for $t >0$ and satisfy the following relations according to Theorem \ref{asymptotic}:
\begin{equation}\label{asympticr}
\left\{
             \begin{array}{lr}
             r(t) \sim (\frac{3}{2})^{\frac{2}{3}} \kappa^{\frac{1}{3}} t^{\frac{2}{3}},  & as ~~~t \rightarrow 0+ \\
             \dot{r}(t) \sim (\frac{2}{3})^{\frac{1}{3}} \kappa^{\frac{1}{3}} t^{-\frac{1}{3}},  & as ~~~t \rightarrow 0+\\
             \ddot{r}(t) \sim -(\frac{2}{81})^{\frac{1}{3}} \kappa^{\frac{1}{3}} t^{-\frac{4}{3}},  & as ~~~t \rightarrow 0+\\
             z_k \rightarrow 0,~~~ k=5, \cdots, 2N,  & as ~~~t \rightarrow 0+
             \end{array}
\right.
\end{equation}
Therefore, $r(t), \theta(t), z_5(t), \cdots, z_{2N}(t)$ of the total collision orbit $\mathbf{r}(t)$ satisfy  the equations (\ref{Euler-Lagrangeequations1detai2}) and (\ref{angularmomentumequation1}) as $t \rightarrow 0+$,
in particular, $PISPW$ exactly explores whether $\theta(t)$ approaches a fixed limit as $t \rightarrow 0+$ in the coordinates $ r, \theta, z_5, \cdots, z_{2N}$.\\

For a total collision orbit, let us introduce the time transformation:
\begin{equation}\label{timetransformation}
dt= r^{\frac{3}{2}}d\tau.
\end{equation}
According to (\ref{asympticr}) (\ref{timetransformation}),  it follows that
 \begin{equation}
\ln t \sim \frac{3\kappa^{\frac{1}{2}}}{2}\tau,  ~~~\ln r \sim \kappa^{\frac{1}{2}}\tau\nonumber
\end{equation}
So we think that $\tau \rightarrow -\infty$  according to $t \rightarrow 0+$ henceforth, provided  discussing  total collision orbits.

Differentiation with respect to time $t$ is denoted by $\dot{ }: \frac{df}{dt}=\dot{f}$ in the previous pages. Similarly,
differentiation with respect to the new variable $\tau$ will be denoted by $': \frac{df}{d\tau}=f'$ henceforth. Thus,
\begin{displaymath}
\dot{f}= r^{-\frac{3}{2}} f', ~~~~~\ddot{f}= r^{-3} f''-\frac{3}{2}r^{-4} r' f'.
\end{displaymath}\\

Then it is easy to see that the equations of motion become:
\begin{equation}\label{equation1}
\left\{
             \begin{array}{l}
             \frac{d}{d\tau}\frac{\partial K(z,z')}{\partial z'_i}-\frac{\partial K(z,z')}{\partial z_i}
 +\frac{{r'}}{2r}\frac{\partial K(z,z')}{\partial z'_i} - \frac{\partial U(z)}{\partial z_i}=0,~~~~~~~~~~~~~~~~~~~~~i=5,\cdots, 2N; \\\\
             {\frac{r''}{r}}-\frac{3}{2}(\frac{{r'}}{r})^2-  2K(z,z')+ U(z)=0,\\\\
             {\theta}'=-\sum_{j,k=5}^{2N}q_{jk}{z}'_j {z}_k.
             \end{array}
\right.\nonumber
\end{equation}
where
\begin{equation*}
   \begin{array}{l}
      K(z,z')=\frac{1}{2} [\frac{ (\sum_{j=5}^{2N}{z}_j z'_j)^2}{z^2_3}+\sum_{j=5}^{2N}z'^2_j -(\sum_{j,k=5}^{2N}q_{jk}z'_j {z}_k)^2]\\
      =\frac{1}{2} \left[\frac{ (\sum_{j=5}^{2N}{z}_j z'_j)^2}{1 - \sum_{j = 5}^{2N} z^2_j}+\sum_{j=5}^{2N}z'^2_j - (\sum_{j,k=5}^{2N}q_{jk}z'_j {z}_k)^2\right],
   \end{array}
\end{equation*}
\begin{equation*}
  \begin{array}{l}
    \frac{\partial K(z,z')}{\partial z'_i}= \frac{z_i \sum_{j=5}^{2N}z'_j z_j}{1 - \sum_{j = 5}^{2N} z^2_j} + z'_i -  \sum_{j=5}^{2N}q_{ij} z_j\sum_{j,k=5}^{2N}q_{jk}z'_j {z}_k,
  \end{array}
\end{equation*}
and
\begin{equation*}
   \begin{array}{l}
      \frac{d}{d\tau}\frac{\partial K(z,z')}{\partial z'_i}-\frac{\partial K(z,z')}{\partial z_i} \\
      =z''_i - \sum_{j=5}^{2N}q_{ij} z_j\sum_{j,k=5}^{2N}q_{jk}z''_j {z}_k- 2  \sum_{j=5}^{2N}q_{ij} z'_j\sum_{j,k=5}^{2N}q_{jk}z'_j {z}_k
               +\frac{z_i \sum_{j=5}^{2N}(z''_j z_j + z'^2_j)}{1 - \sum_{j = 5}^{2N} z^2_j} + \frac{{z}_i (\sum_{j=5}^{2N}{z}_j z'_j)^2}{(1 - \sum_{j = 5}^{2N} z^2_j)^2}.
   \end{array}
\end{equation*}

We introduce new variables:
\begin{displaymath}
Z_k=z'_k, ~~~ k=5,\cdots, 2N, ~~~~~~~~~~\Upsilon = \frac{r'}{r}, ~~~~~~~~~~\Theta = \theta'.
\end{displaymath}
Then  the equations of motion above become:
\begin{equation}\label{equationzr3}
\left\{
             \begin{array}{lr}
             z'_i  =  Z_i, &  \\\\
             Z'_i  =  \frac{\partial U(z)}{\partial z_i} - \frac{z_i \sum_{j=5}^{2N}({Z_j'} z_j + {Z_j}^2)}{z^2_3} - \frac{{z}_i (\sum_{j=5}^{2N}{z}_j  {Z_j})^2}{z^4_3}-  {\sum_{j,k=5}^{2N}q_{kj}{Z_j}' {z}_k} \sum_{j=5}^{2N}q_{ij} z_j\\
- 2 \sum_{j,k=5}^{2N}q_{kj}{Z_j} {z}_k \sum_{j=5}^{2N}q_{ij} {Z_j}-\frac{\Upsilon}{2} [ {Z_i} +\frac{z_i \sum_{j=5}^{2N}{Z_j} z_j}{z^2_3} +  \sum_{j,k=5}^{2N}q_{kj}{Z_j} {z}_k \sum_{j=5}^{2N}q_{ij} z_j], &\\\\
             r'  =  r\Upsilon, &  \\\\
             \Upsilon'  =  \frac{1}{2}\Upsilon^2+ \frac{ (\sum_{j=5}^{2N}{z}_j  {Z_j})^2}{z^2_3}+\sum_{j=5}^{2N}{Z_j}^2 - (\sum_{j,k=5}^{2N}q_{jk}{Z_j} {z}_k)^2 -   U(z); &
             \end{array}
\right.
\end{equation}
and
\begin{equation}\label{equationtheta3}
\left\{
             \begin{array}{lr}
             \theta'  =  \Theta, &  \\
             \Theta  =  \sum_{j,k=5}^{2N}q_{kj}{Z_j} {z}_k. &
             \end{array}
\right.
\end{equation}

It is noteworthy that, the system of the first, second and last equations in (\ref{equationzr3}) is autonomous in $\Upsilon, z, Z$. Once $\Upsilon, z, Z$ are solved by them, the variables  $r, \theta$ can also be solved by quadrature.

Incidentally, the relation of total energy becomes
\begin{equation}\label{totalenergyconserved2}
   2r \mathcal{H}={\Upsilon}^2+ [\frac{ (\sum_{j=5}^{2N}{z}_j Z_j)^2}{z^2_3}+\sum_{j=5}^{2N}Z^2_j - (\sum_{j,k=5}^{2N}q_{jk}Z_j {z}_k)^2]- 2U(z),
\end{equation}
or\begin{equation}
   2r \mathcal{H}=\Upsilon'+\frac{1}{2}\Upsilon^2- U(z).\nonumber
\end{equation}

\begin{remark}
Consider the system (\ref{equationzr3})
in
$ (z, Z, r, \Upsilon)\in \mathbf{B}^{2N-4}\times \mathbb{R}^{2N-4}\times [0,+\infty)\times \mathbb{R}$.
It is easy to see that the system (\ref{equationzr3}) does not have singularities at $r=0$ and has an invariant manifold $\{r=0\}$.
Set
\begin{equation}
    \mathcal{N}_{\mathcal{H}}=\{(z, Z, r, \Upsilon)\in \mathbf{B}^{2N-4}\times \mathbb{R}^{2N-4}\times [0,+\infty)\times \mathbb{R}| (\ref{totalenergyconserved2})~~~ holds\}\nonumber
\end{equation}
for a fixed real number $\mathcal{H}$. Then $\mathcal{N}_{\mathcal{H}}$ is an invariant manifold of the system (\ref{equationzr3}).
 It is similar to results  of McGehee in \cite{Mcgehee1974Triple}, one can show that the flow of the system (\ref{equationzr3}) restricted to the invariant manifold
 \begin{center}
 $\{r=0\}\bigcap\mathcal{N}_{\mathcal{H}}$\\$=\{(z, Z, 0, \Upsilon)\in \mathcal{N}_{\mathcal{H}}|{\Upsilon}^2+ [\frac{ (\sum_{j=5}^{2N}{z}_j Z_j)^2}{z^2_3}+\sum_{j=5}^{2N}Z^2_j - (\sum_{j,k=5}^{2N}q_{jk}Z_j {z}_k)^2]- 2U(z)=0\}$
 \end{center}
 is gradient-like with respect
to
\begin{equation*}
    \pi_{\Upsilon}: \{r=0\}\bigcap\mathcal{N}_{\mathcal{H}}\rightarrow \mathbb{R}: (z, Z, 0, \Upsilon)\mapsto -\Upsilon.
\end{equation*}
For more detail please refer to Appendix \ref{Gradient-likeFlow}.
\end{remark}

It is easy to see that the  total collision solution $\mathbf{r}(t)$ of equations (\ref{eq:Newton's equation1}) corresponds to  a solution
$( z, Z,r, \Upsilon, \theta)$ of the  equations (\ref{equationzr3})  and (\ref{equationtheta3}) such that $r \rightarrow 0, z \rightarrow 0$, where  $Z=(Z_5, \cdots, Z_{2N})^\top$.

It follows from Proposition \ref{isolatedcriticalpoint}  that all the equilibrium points of equations (\ref{equationzr3}) and (\ref{equationtheta3}), which geometrically constitute a circle, satisfy $  z=0, Z=0,r=0, \Upsilon = \pm \kappa^{\frac{1}{2}},\theta=const $ so long as $z$ is small. Furthermore, we have
\begin{proposition}\label{asymptic2}
The total collision orbit $\mathbf{r}(t)$ of equations (\ref{eq:Newton's equation1})  corresponds exactly to a solution $( z(\tau), Z(\tau),r(\tau), \Upsilon(\tau), \theta(\tau))$ of equations (\ref{equationzr3}) and (\ref{equationtheta3}) such that
\begin{equation}
z \rightarrow 0, ~~~Z \rightarrow 0, ~~~r \rightarrow 0, ~~~\Upsilon \rightarrow \kappa^{\frac{1}{2}}, ~~~\theta'=\Theta \rightarrow 0, ~~~as~ \tau \rightarrow -\infty.\nonumber
\end{equation}
\end{proposition}
{\bf Proof.} 

We only need to prove a solution $(z(\tau), Z(\tau),r(\tau), \Upsilon(\tau), \theta(\tau))$ of equations (\ref{equationzr3}) and (\ref{equationtheta3}) corresponding  to the total collision orbit $\mathbf{r}(t)$ satisfies
\begin{equation}
  Z \rightarrow 0, ~~~\Upsilon \rightarrow \kappa^{\frac{1}{2}}.\nonumber
\end{equation}

First, by (\ref{asympticr}), it is easy to show that $\Upsilon \rightarrow \kappa^{\frac{1}{2}}$.

Let $\alpha(\mathbf{r})$  denote the $\alpha$-limit set of the solution $( z(\tau), Z(\tau),r(\tau), \Upsilon(\tau), \theta(\tau))$.
Then
\begin{equation}
\alpha(\mathbf{r}) \subset \{(z, Z,r, \Upsilon, \theta)| r=0, \Upsilon = \kappa^{\frac{1}{2}},  z=0\}.\nonumber
\end{equation}
Let us  investigate the maximum invariant set included in
\begin{equation}
\{(z, Z,r, \Upsilon, \theta)| r=0, \Upsilon = \kappa^{\frac{1}{2}},  z=0\},\nonumber
\end{equation}
it is easy to show that this set is precisely
\begin{equation}
\{(z, Z,r, \Upsilon, \theta)| r=0, \Upsilon = \kappa^{\frac{1}{2}}, z=0, Z=0\}.\nonumber
\end{equation}
The relation $Z \rightarrow 0$ therefore follows from
\begin{equation}
\alpha(\mathbf{r}) \subset \{(z, Z,r, \Upsilon, \theta)| r=0, \Upsilon = \kappa^{\frac{1}{2}}, z=0, Z=0\}.\nonumber
\end{equation}
$~~~~~~~~~~~~~~~~~~~~~~~~~~~~~~~~~~~~~~~~~~~~~~~~~~~~~~~~~~~~~~~~~~~~~~~~~~~~~~~~~~~~~~~~~~~~~~~~~~~~~~~~~~~~~~~~~~~~~~~~~~~~~~~~~~~~~~~~~~~~~~~~~~~~~~~~~~~~~~~~~~~\Box$\\

Let us  investigate all the solutions of equations (\ref{equationzr3}) and (\ref{equationtheta3}) satisfying the above asymptotic conditions.

Now define a new variable $\gamma = \Upsilon - \kappa^{\frac{1}{2}}$ and substituting $\gamma$ into equations (\ref{equationzr3}) and (\ref{equationtheta3}), we have
\begin{equation}\label{equationz4}
\begin{pmatrix}z'\\Z'\\\gamma'\end{pmatrix}=\mathfrak{A}\begin{pmatrix}z\\Z\\\gamma\end{pmatrix}+\begin{pmatrix}0\\\chi_Z(z,Z,\gamma)\\\chi_0(z,Z,\gamma)\end{pmatrix}
\end{equation}
\begin{equation}\label{equationr4}
r'  =  r (\kappa^{\frac{1}{2}} + \gamma);
\end{equation}
\begin{equation}\label{equationtheta4}
\theta'  =  \sum_{j,k=5}^{2N}q_{kj}{Z_j} {z}_k;
\end{equation}
where 
\begin{equation}\label{}
\mathfrak{A}= \left(
           \begin{array}{ccc}
              0 & \mathbb{I} &   \\
             \Lambda & - \frac{\kappa^{\frac{1}{2}}}{2}\mathbb{I} &   \\
               &   & \kappa^{\frac{1}{2}} \\
           \end{array}
         \right)\nonumber
\end{equation}
denotes the square matrix of the  coefficients  of the
linear terms,   $\Lambda = diag( \mu_5, \cdots,  \mu_{2N})$;
$\chi_Z=(\chi_5,\cdots,\chi_{2N})^\top$ and
the functions $\chi_k,\chi_0$ ($k=5, \cdots, 2N$) are power-series in the $4N-7$ real variables $z,Z,\gamma$  starting with quadratic terms and  all converge for sufficiently small  $z, Z, \gamma$:
\begin{equation}\label{chik}
\begin{array}{lr}
             \chi_k (z,Z,\gamma)=\frac{1}{2} [ a_{kkk}z^2_k + 2 \sum_{j=5,j\neq k}^{2N} a_{jkk} {z_k} z_j +  \sum_{i,j=5, i,j\neq k}^{2N}a_{ijk}{z}_i {z_j} - \gamma Z_k] + \cdots  & \\
               =  \frac{1}{2} [ \sum_{i,j=5, }^{2N}a_{ijk}{z}_i {z_j} - \gamma Z_k] + \cdots,  &
             \end{array}
\end{equation}
\begin{equation}\label{chi}
\chi_0 (z,Z,\gamma)=\frac{1}{2}\gamma^2 + \sum_{j=5}^{2N}{Z_j}^2 - \frac{1}{2}\sum_{j=5}^{2N}{\mu_j}{z_j}^2 + \cdots.
\end{equation}

As a result, total collision orbits of the planar $N$-body problem have the following interesting characterization.
\begin{theorem}\label{asymptic3}
A total collision orbit $\mathbf{r}(t)$ of equations (\ref{eq:Newton's equation1})  reduces exactly to a solution $(z(\tau), Z(\tau),\gamma(\tau))$ of equations (\ref{equationz4}) corresponding to some central configuration $\mathbf{r}_0$ such that
\begin{equation}
z \rightarrow 0, ~~~Z \rightarrow 0, ~~~\gamma \rightarrow 0 ~~~~~~as~ \tau \rightarrow -\infty.\nonumber
\end{equation}
Therefore, the manifold of all the total collision orbits corresponds exactly to the finite union of the unstable manifold  of the systems like type (\ref{equationz4}) at the origin.
\end{theorem}

We conclude this subsection by reviewing what we have accomplished so far. For a fixed central configuration $\mathbf{r}_0$, we introduced a moving frame to describe equations of motion for the $N$-body problem in the phase space $\mathrm{T}(\mathcal{X}\backslash \mathcal{P}^{\bot}_{\mathbf{r}_0})\subset\mathrm{T}\mathcal{X}$. Given masses $m_1, ..., m_N$, suppose that there are $n$ central configurations,   $\mathbf{c}_i$ $(i=1,\cdots,n)$.
	For a fixed total collision orbit $\mathbf{r}(t)$, suppose that 
\begin{equation*}
   \widehat{\mathbf{r}}(t) \rightarrow \mathbf{S} (= \{e^{\mathbf{i}\theta}\mathbf{c}_1| \theta \in \mathbb{R}\}) ~~~~~~as~~ t \rightarrow 0+.
\end{equation*}
Then when choosing $\mathbf{r}_0=\mathbf{c}_1$, the total collision orbit $\mathbf{r}(t)$ reduces to a solution $(z(\tau), Z(\tau),\gamma(\tau))$ of equations (\ref{equationz4}) such that
\begin{equation}
z \rightarrow 0, ~~~Z \rightarrow 0, ~~~\gamma \rightarrow 0 ~~~~~~as~ \tau \rightarrow -\infty.\nonumber
\end{equation}
Conversely, a solution $(z(\tau), Z(\tau),\gamma(\tau))$ of equations (\ref{equationz4}) satisfying the relations above
is obviously corresponding to some total collision orbit. We remark that, when choosing $\mathbf{r}_0=\mathbf{c}_i$ for some $i\neq 1$, the moving frame and the concomitant coordinate system  can still describe the total collision orbit $\mathbf{r}(t)$ well, provided $\mathbf{c}_1 \notin \mathcal{P}^{\bot}_{\mathbf{r}_0}$ $(=\mathcal{P}^{\bot}_{\mathbf{c}_i})$; in this case, $\mathbf{r}(t)$ reduces to a solution $(z(\tau), Z(\tau),r(\tau),\Upsilon(\tau))$ of equations (\ref{equationzr3}) such that
\begin{equation}
z \rightarrow z_0, ~~~Z \rightarrow 0, ~~~r \rightarrow 0, ~~~\Upsilon \rightarrow \sqrt{2U(z_0)}  ~~~~~~as~ \tau \rightarrow -\infty,\nonumber
\end{equation}
where $z_0$ is a certain critical point of $U(z)$ in $\mathbf{B}^{2N-4}$.

\subsection{The Problem of Infinite Spin or Painlev\'{e}-Wintner ($PISPW$)}
\indent\par
Now $PISPW$ can be formulated as:
\begin{problem}\emph{(\textbf{$PISPW$})}
 for all the solutions of equations (\ref{equationz4}) satisfying
\begin{equation}\label{Asymptotic condition}
z \rightarrow 0, ~~~Z \rightarrow 0, ~~~\gamma \rightarrow 0
\end{equation}
whether does $\theta(\tau)$ in (\ref{equationtheta4}) approach a fixed limit as $\tau \rightarrow -\infty$?
\end{problem}

$PISPW$ can also be stated as following:
\begin{problem}
given a solution $(z, Z, \gamma)$ of equations (\ref{equationz4}), whether does the implication
\begin{displaymath}
\alpha(z, Z, \gamma) = \{0\} ~~~~~\Rightarrow ~~~~~\alpha(z, Z, \gamma, \theta)  ~is ~a ~single ~point,
\end{displaymath}
hold?
\end{problem}
Where $\alpha(z, Z, \gamma)$ denotes the $\alpha$-limit set of the solution $(z, Z, \gamma)$ of equations (\ref{equationz4}) and  $\alpha(z, Z, \gamma, \theta)$  denotes the $\alpha$-limit set of the solution $(z, Z, \gamma, \theta)$ of equations (\ref{equationz4}) (\ref{equationtheta4}).

We remark that there is exactly one equilibrium point $(z, Z, \gamma)=0$ of equations (\ref{equationz4}) in some small neighbourhood of the original point $(z, Z, \gamma)=0$ according to isolation of central configurations. So to solve $PISPW$ is equivalent to judge
\begin{equation}\label{}
\mathcal{W}^u(\Sigma)=\bigcup_{p \in \Sigma}\mathcal{W}^u(p)\nonumber
\end{equation}
where $\Sigma=\{(z, Z, \gamma, \theta)|z=0, Z=0, \gamma=0, \theta \in \mathbb{R}\}$ is the set of all the equilibrium points of equations (\ref{equationz4}) (\ref{equationtheta4}).\\

Before discussing formally the above problem, we wish to give some examples to illustrate some ideas  proper or not for $PISPW$.

\begin{figure}
  \center
  \includegraphics[width=7cm]{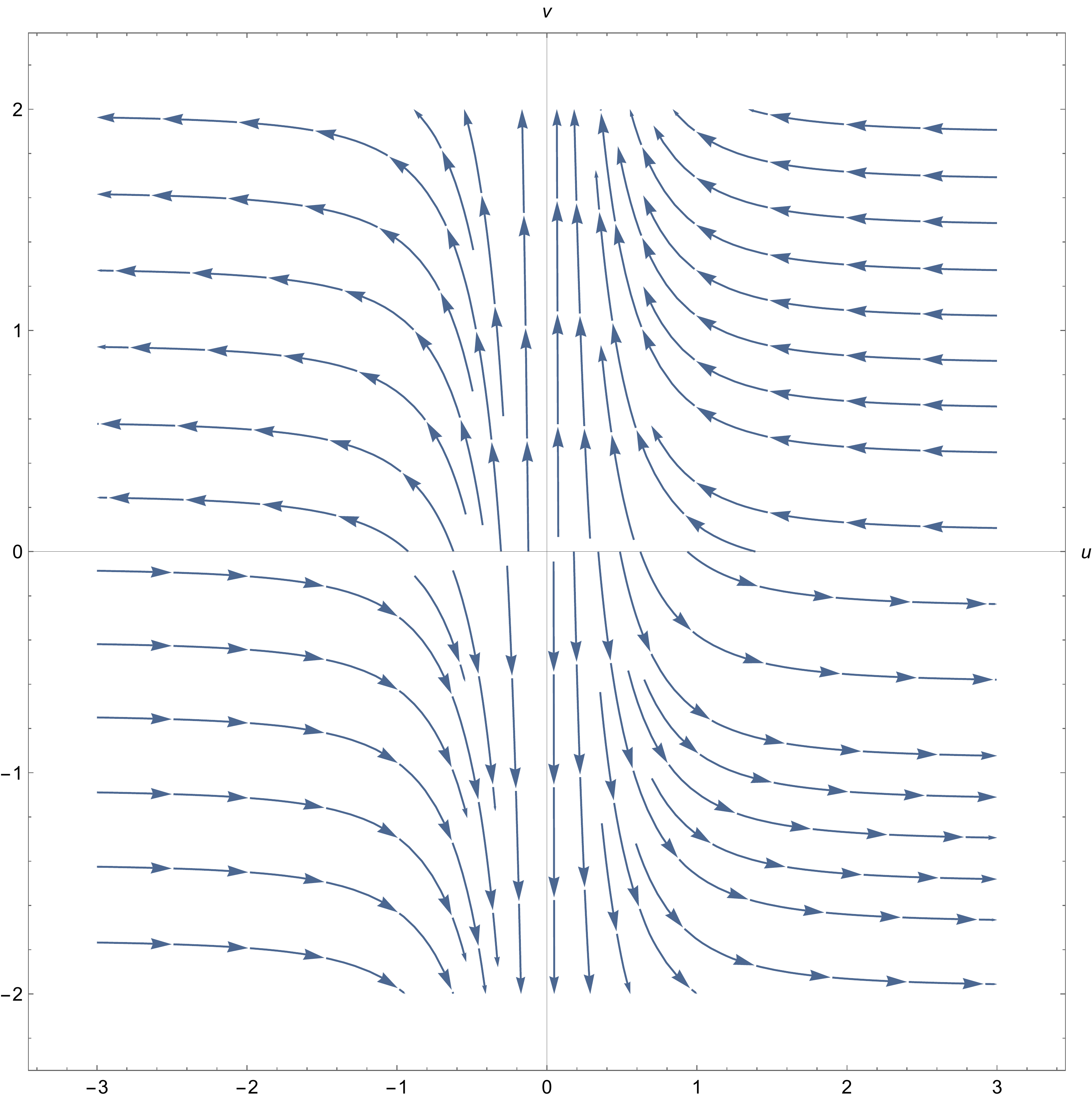}\\
  \caption{phase portrait of equations (\ref{Example 1}) in {Example 1} }\label{phase portrait1}
\end{figure}
\begin{figure}
  \center
  \includegraphics[width=6cm]{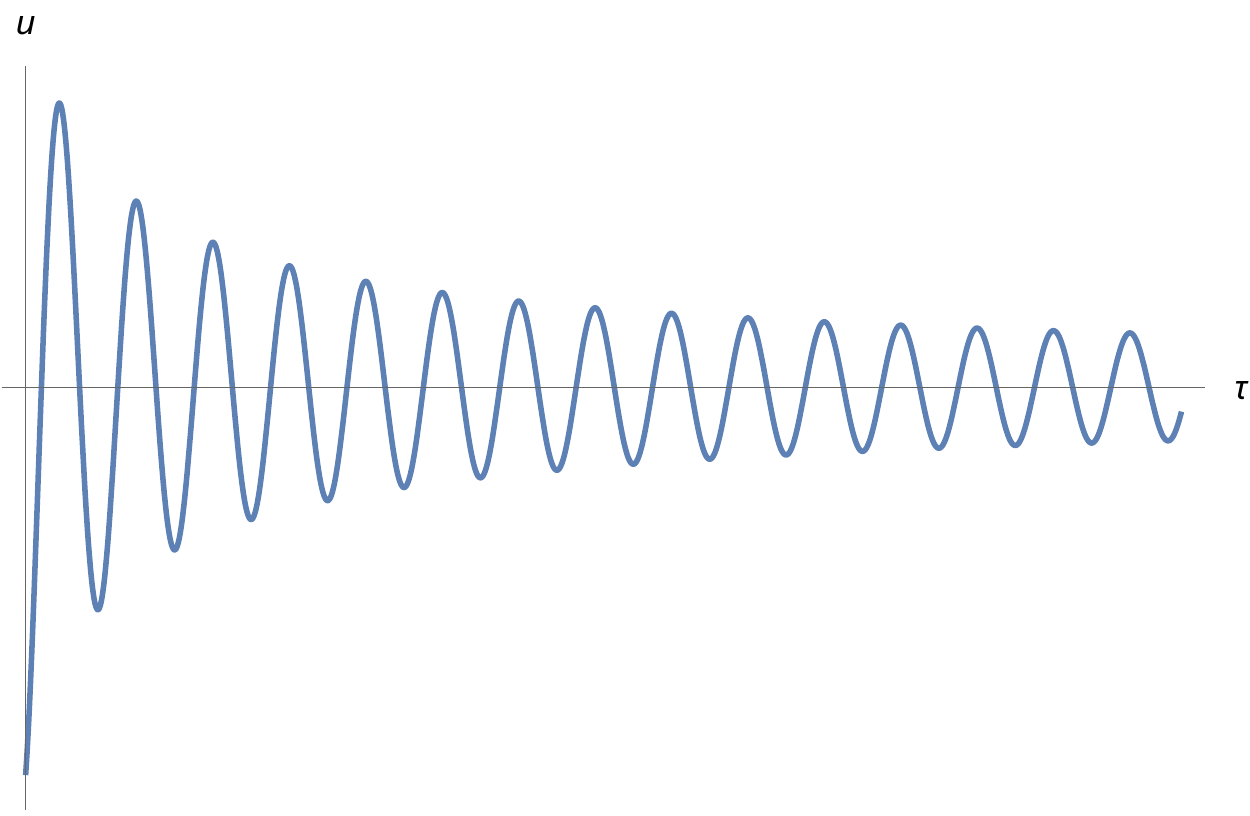}~~~~~~~~~~\includegraphics[width=6cm]{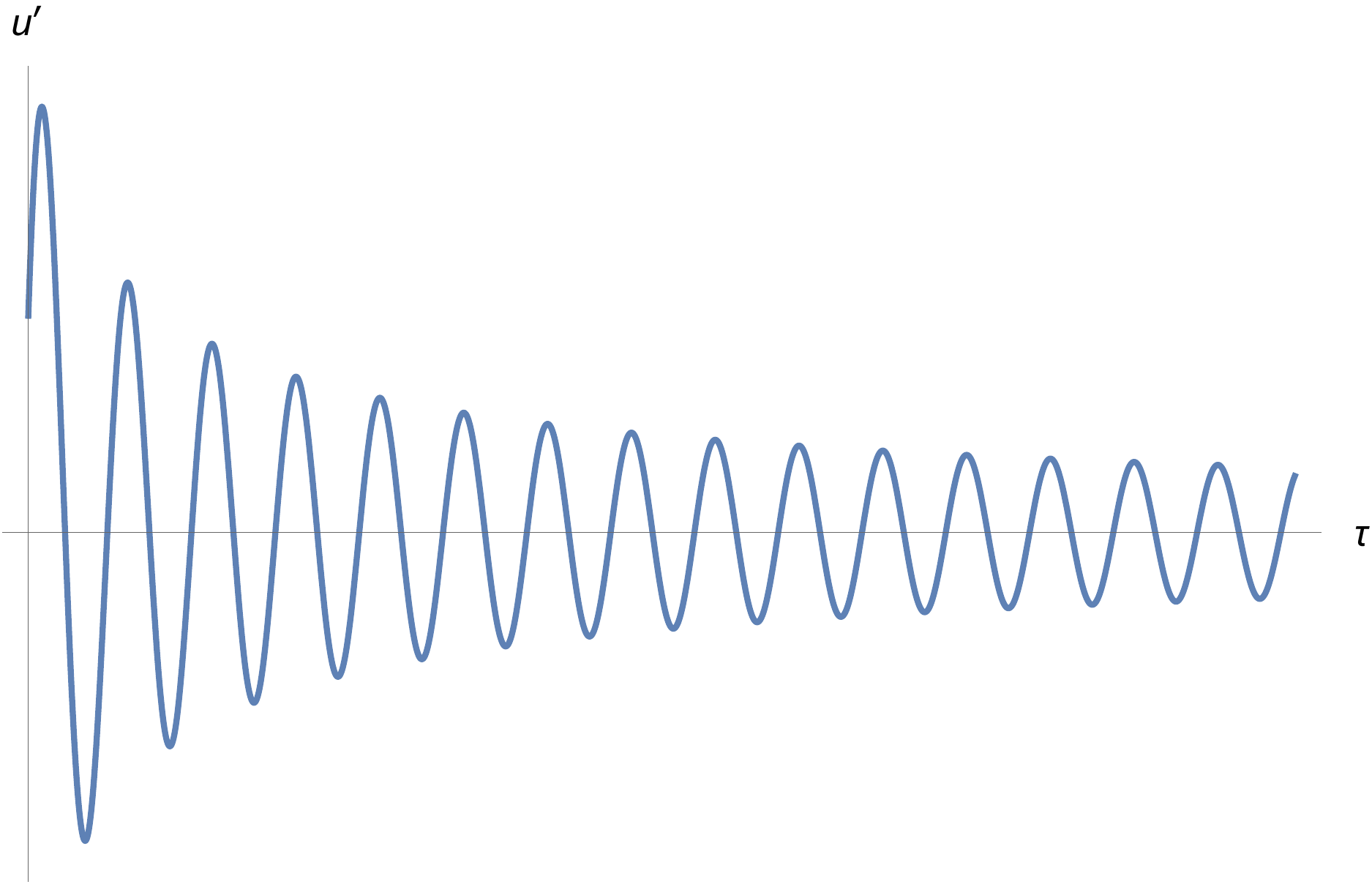}\\
  \includegraphics[width=6cm]{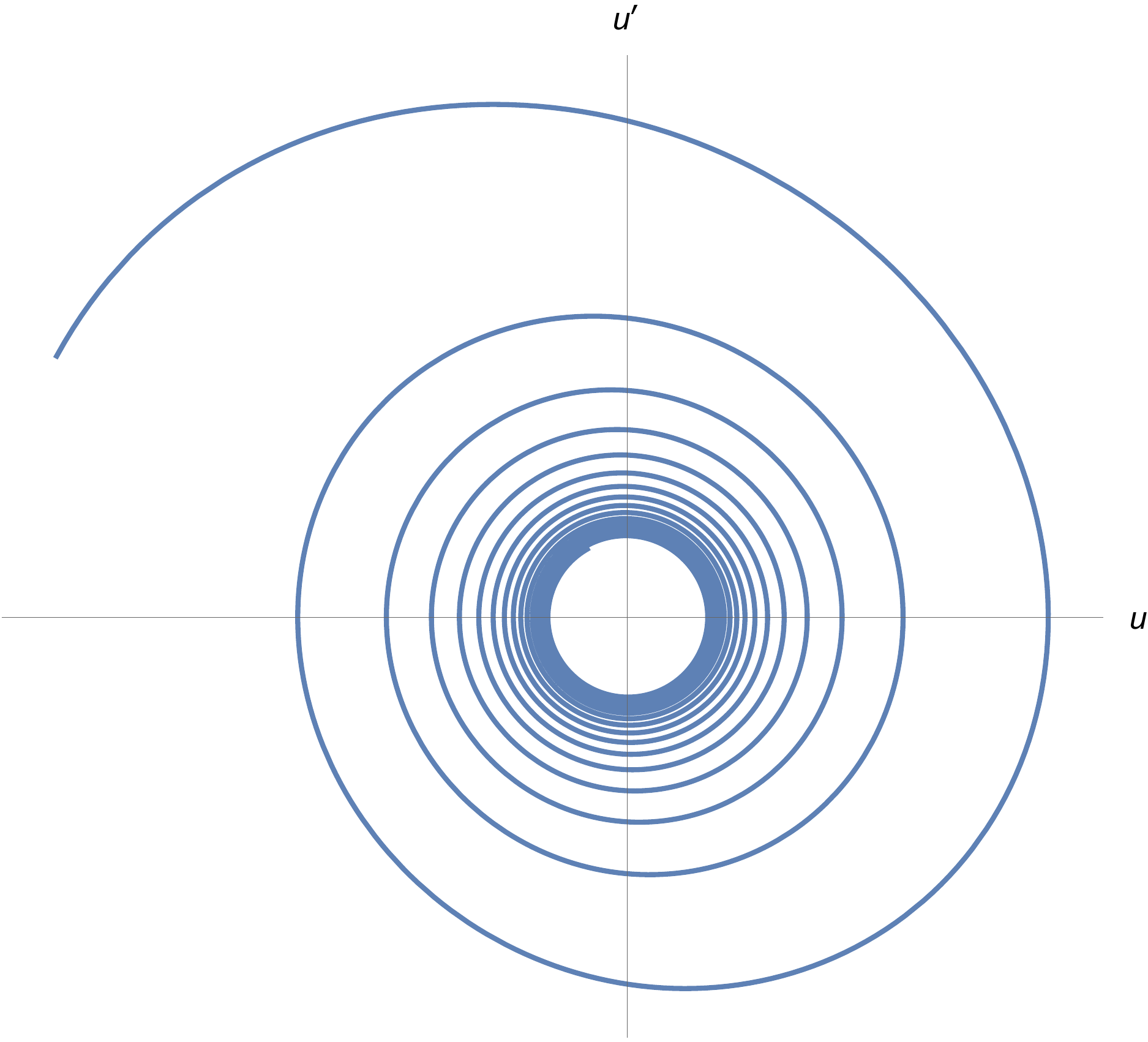}~~~~~~~~~~\includegraphics[width=6cm]{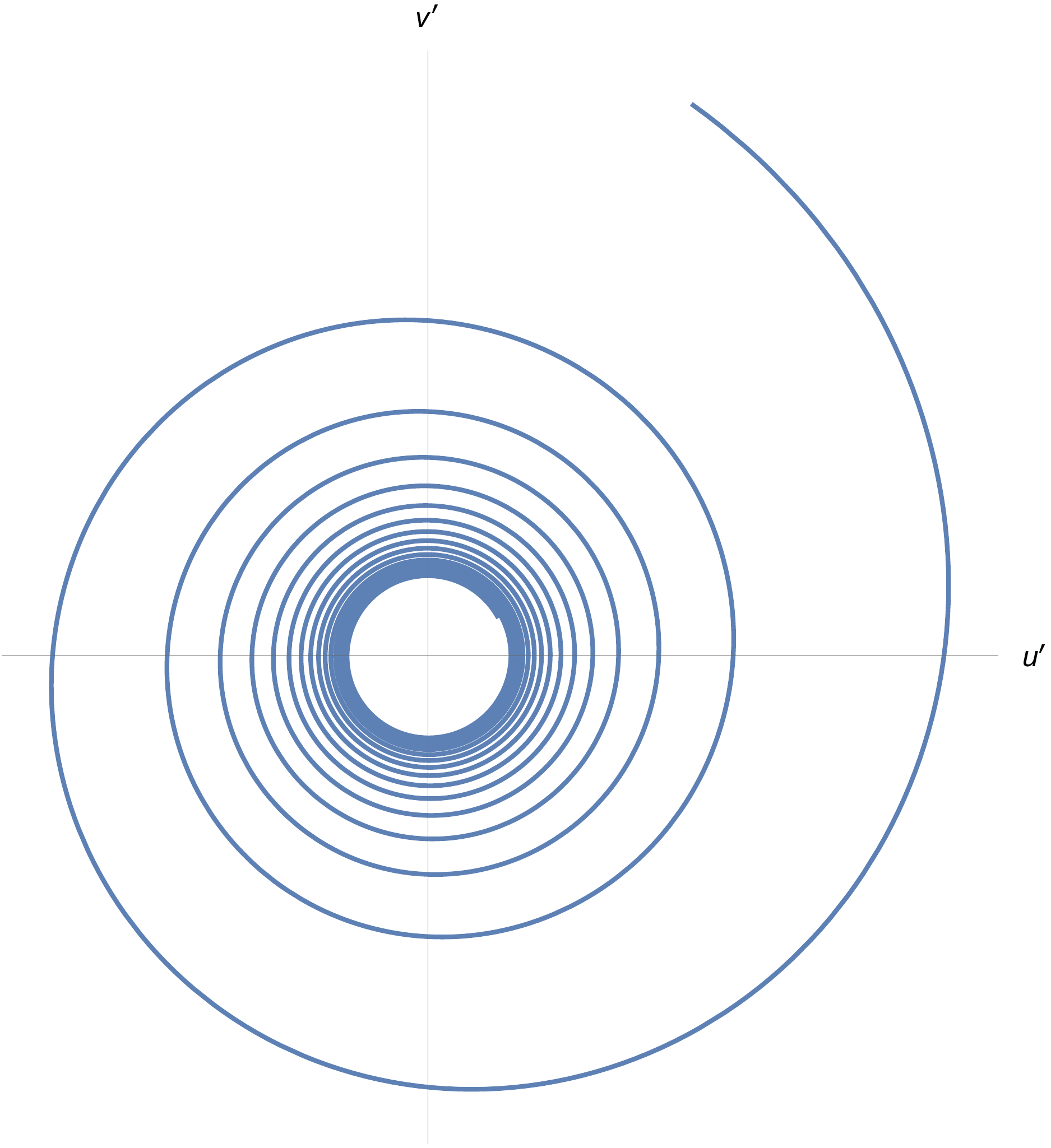}\\
  \caption{plots  in {Example 2}}\label{utao}
\end{figure}

{\bf Example 1.} Consider the system of differential equations
\begin{equation}\label{Example 1}
\left\{
             \begin{array}{lr}
             u'  = - u^2 (u^2+1)v &  \\
             v'  =  v &
             \end{array}
\right.
\end{equation}
The set of all the equilibrium points of the above equations is the $u$-axis $\Sigma =\{(u,v)|v=0\}$ and forms an invariant manifold, furthermore, $\Sigma$ is  a  central manifold (see \cite{Shilnikov,zbMATH03727247}).
We cannot simply utilize the theory of  central manifolds to prove that $\mathcal{W}^u(\Sigma)=\bigcup_{p \in \Sigma}\mathcal{W}^u(p)$, that is, it does not simply follow from the theory of central manifolds that $v(\tau) \rightarrow 0$ implies that $u(\tau)$ approaches a fixed limit as $\tau \rightarrow -\infty$. Indeed, the solution
\begin{equation}
\left\{
             \begin{array}{lr}
             \frac{1}{u} + \arctan u =  \exp \tau + \frac{\pi}{2} &  \\
             v  =  \exp \tau &
             \end{array}
\right.\nonumber
\end{equation}
or phase portrait of equations (\ref{Example 1}) as Figure \ref{phase portrait1} can show the statement.

\begin{remark}
If one can deduce that $\Sigma$ is a normally hyperbolic invariant manifold by the normal (to $\Sigma$) eigenvalue $1$ lies off the  imaginary axis, then  Example 1 shows that
$\mathcal{W}^u(\Sigma)=\bigcup_{p \in \Sigma}\mathcal{W}^u(p)$ may be wrong  for noncompact invariant manifold $\Sigma$.

On the other hand, for an invariant manifold $\Sigma$ of a vector field $\mathbf{v}$, except the case of $\Sigma = one ~point$, even if $\Sigma$ is a circle, the condition that the linearized vector field of $\mathbf{v}$ at $\Sigma$  have ``its normal (to $\Sigma$) eigenvalues off the  imaginary axis is neither necessary nor sufficient for the $\mathbf{v}$-flow. It thus remains an open, fuzzy question to formulate an integrated conditions on $\mathbf{v}$ at $\Sigma$ that guarantee normally hyperbolicity of the $\mathbf{v}$-flow" \cite[p.8]{Hirsch1970Invariant}.
\end{remark}

{\bf Example 2.}  Consider the following functions
\begin{equation}\label{Example 2}
\left\{
             \begin{array}{lr}
             u (\tau) = \frac{1}{\sqrt{\tau \ln \tau}} \sin \tau &  \\
             v (\tau) =  \frac{1}{\sqrt{\tau \ln \tau}} \cos \tau &
             \end{array}
\right.\nonumber
\end{equation}

Then $u (\tau), v (\tau), u' (\tau), v' (\tau)$ approach zero as $\tau \rightarrow +\infty$, however, it is easy to see that all the following improper integrals
\begin{displaymath}
\int^{+\infty}u (\tau)v' (\tau)d\tau, \int^{+\infty} v (\tau) u' (\tau) d\tau, \int^{+\infty}u (\tau)v' (\tau)- v (\tau) u' (\tau)d\tau
\end{displaymath}
are  not convergent.

Similarly, we cannot simply claim that $\theta(\tau)$  approaches a fixed limit as $\tau \rightarrow -\infty$, although
\begin{equation}
z \rightarrow 0, ~~~Z \rightarrow 0,\nonumber
\end{equation} and $\int_{-\infty}Z(\tau) d\tau$
is convergent.\\\\

Let us finish the section with further simplifying equations (\ref{equationz4}) to make preparations for resolving $PISPW$.

The aim  is to diagonalize the linear part of equations (\ref{equationz4}). More specifically,
 it is well known that the $4N-7$-rowed square matrix $\mathfrak{A}$  in equations (\ref{equationz4})  is generically diagonalizable, that is,  $\mathfrak{A}$ is similar to a diagonal matrix, provided that $\mathfrak{A}$ has $4N-7$ distinct eigenvalues; on the other hand,  $\mathfrak{A}$ is always ``almost" diagonalizable,  that is,  for any $\epsilon > 0$, there exists an upper triangular matrix, such that
 all the elements  above the diagonal are less than $\epsilon$ in absolute value,  is similar to $\mathfrak{A}$. Indeed, we find that the $4N-7$-rowed square matrix
 \begin{center}
 $\mathfrak{C}:=\begin{pmatrix}\mathfrak{P}^{-1}&\\&1\end{pmatrix} \mathfrak{A} \begin{pmatrix}\mathfrak{P}&\\&1\end{pmatrix}$
\end{center}
is an appropriative almost-diagonal-matrix  for any $\epsilon > 0$. Here and below,  please refer to Appendix \ref{DiagonalizationoftheLinearPart} for more detail.

Consequently,
 after applying the linear substitution
\begin{displaymath}
\left(
  \begin{array}{c}
    z \\
    Z \\
  \end{array}
\right)
=\mathfrak{P} q
\end{displaymath}
to the equations (\ref{equationz4}),  we arrive at the following equations which can be better handled:
\begin{equation}\label{equationzq}
\begin{pmatrix}q'\\
  \gamma'\end{pmatrix}=\mathfrak{C}\begin{pmatrix}q\\
  \gamma\end{pmatrix}-\varphi(q, \gamma),
\end{equation}
where
\begin{center}
$q=(q_5,\cdots,q_{2N}, q_{2N+5},\cdots,q_{4N})^\top$,
\end{center}
\begin{equation*}
    \varphi = ( \varphi_{5}, \cdots, \varphi_{2N},-\varphi_{5},\cdots, -\varphi_{2N},\varphi_{0})^\top,
\end{equation*}
the functions $\varphi_k, \varphi_0$ are power-series in the $4N-7$ variables $q,\gamma$  starting with quadratic terms:
\begin{equation}\label{varphi}
\left\{
             \begin{array}{lr}
             \varphi_k(q, \gamma)= \frac{1}{2 \sqrt{ \mu_{k} + \frac{\kappa}{16}}}\chi_k(z,Z,\gamma),  & k\in \{ 5, \cdots, 2N-n_3\} \\
             \varphi_k(q, \gamma)= \frac{1}{\epsilon}\chi_k(z,Z,\gamma),  & k\in \{2N-n_3 +1, \cdots, 2N\}\\
\varphi_0(q, \gamma)= -\chi_0(z,Z,\gamma); &
             \end{array}
\right.
\end{equation}
\begin{displaymath}
\tilde{\mu}_j = -\frac{\kappa^{\frac{1}{2}}}{4}+\sqrt{ \mu_j + \frac{\kappa}{16}} ~~~~for ~~~~j\in \{5, \cdots, 2N \};
\end{displaymath}
and $n_0, n_p, n_1, n_2, n_3$ are natural numbers such that
\begin{center}
$n_0+n_p+n_1+n_2+n_3=2N-4$,
\end{center}
\begin{displaymath}
\mu_j =0 ~~~~for ~~~~j\in \{5, \cdots, n_0+4 \},
\end{displaymath}
\begin{displaymath}
\mu_j >0 ~~~~for ~~~~j\in \{n_0+5, \cdots, n_0+n_p+4 \},
\end{displaymath}
\begin{displaymath}
0> \mu_j > - \frac{\kappa}{16} ~~~~for ~~~~j\in \{n_0+n_p+5, \cdots, n_0+n_p+n_1+4 \},
\end{displaymath}
\begin{displaymath}
\mu_j < - \frac{\kappa}{16} ~~~~for ~~~~j\in \{n_0+n_p+n_1+5, \cdots, n_0+n_p+n_1+n_2+4 \},
\end{displaymath}
\begin{displaymath}
\mu_j = - \frac{\kappa}{16} ~~~~for ~~~~j\in \{n_0+n_p+n_1+n_2+5, \cdots, 2N \}.
\end{displaymath}
Furthermore, note the following facts:
\begin{displaymath}
\tilde{\mu}_j =0 ~~~~for ~~~~j\in \{5, \cdots, n_0+4 \}
\end{displaymath}
\begin{displaymath}
\tilde{\mu}_j >0 ~~~~for ~~~~j\in \{ n_0+5, \cdots, n_0+n_p+4\}
\end{displaymath}
\begin{displaymath}
Re\tilde{\mu}_j <0 ~~~~for ~~~~j\in \{ n_0+n_p+5, \cdots, 2N\}
\end{displaymath}
and $n_p \geq N-2$.

After applying the above linear substitution
to equations (\ref{equationtheta4}), it follows that
\begin{equation}\label{equationtheta5}
\begin{array}{c}
  \theta'  =z^\top Q Z=  q^\top \mathfrak{P}^\top \begin{pmatrix}\mathbb{I}\\
  0\end{pmatrix}
Q \begin{pmatrix}0&\mathbb{I}\end{pmatrix}
\mathfrak{P}q \\
   = \sum_{5\leq k\leq n_0+4}\sum_{j=n_0+5}^{n_0+n_p+4}q_{kj}\tilde{{\mu}}_j {q}_k {q_j}+\sum_{j,k=n_0+5}^{n_0+n_p+4}q_{kj}\tilde{{\mu}}_j {q}_k {q_j}  + \cdots,
\end{array}
\end{equation}
the right hand side is a quadratic form of  of $q$, where
``$\cdots$" denotes all the quadratic terms which  contain at least one of $q_k$ $(k>n_0+n_p+4)$ as a factor.

\section{Resolving PISPW}
\indent\par
Let us turn to discuss $PISPW$ in this section. As  the work of Siegel  \cite{C1967Lectures} indicates, for discussing $PISPW$, it is important to make proper use of  results of normal forms. On the other hand, in the process of discussing $PISPW$ for the case corresponding to  central configurations with degree of degeneracy two, it has been natural to investigate  equilibrium points  in two-dimensional systems and   degenerate central configurations of the planar four-body problem. For the sake of readability, these materials are relegated to the Appendix.

In the following we divide the discussion into several cases according to the degree of degeneracy of central configurations.

\subsection{$n_0=0$}\label{subsectionn0=0}
\indent\par

First, let us discuss the case that $n_0=0$, i.e., the central configuration is nondegenerate.

Now the problem attributes to the following form:
\begin{problem}\label{problemn_0is0}
given a solution $(q(\tau),\gamma(\tau))$ of the system (\ref{equationzq}) such that
 \begin{center}
 $(q(\tau),\gamma(\tau))\rightarrow 0$ as $\tau \rightarrow -\infty$,
\end{center}
if $\theta(\tau)$ satisfies (\ref{equationtheta5}) then whether does $\theta(\tau)$  approach a fixed limit as $\tau \rightarrow -\infty$?
\end{problem}

By solving the above problem, we have
\begin{theorem}\label{nondegenerateproblem}
If all the central configurations are nondegenerate  for given  masses of the planar $N$-body problem, then the normalized configuration of any total collision orbit of the given   $N$-body problem
approaches a certain central configuration as time $t$ approaches the collision instant.
\end{theorem}

\textbf{Proof.}

It is well known that an orbit on the
unstable manifold of a hyperbolic equilibrium point  approaches the hyperbolic equilibrium point exponentially (for example, see Corollary \ref{hyperbolicapproximatetheorem}).
Since the equilibrium
point $0$ is hyperbolic for the system (\ref{equationzq}),
it follows that,
 for a fixed $\tau_0$, there are two positive constants $\varpi_1$ and $\sigma$ such that
\begin{displaymath}
\|u\| \leq \varpi_1 e^{\sigma \tau} ~~~~ ~~~~ \forall\tau \leq \tau_0,
\end{displaymath}
where
\begin{center}
$u=(\gamma, q_5,\cdots,q_{2N}, q_{2N+5},\cdots,q_{4N})^\top$
\end{center}
 and
\begin{displaymath}
 \|u\|= \sqrt{|\gamma|^2+|q_5|^2+\cdots+|q_{2N}|^2+ |q_{2N+5}|^2+\cdots+|q_{4N}|^2}.
\end{displaymath}

Thanks to the equation (\ref{equationtheta5}),   $\theta'$ is an analytic function of $u$ in the vicinity of the origin. For ease of notations, set
\begin{equation}\label{equationtheta6}
\theta'  =  \sum_{|\alpha|=2}c_\alpha u^\alpha,\nonumber
\end{equation}
where
\begin{center}
$\alpha=(\alpha_0,\alpha_5,   \cdots,\alpha_{2N}, \alpha_{2N+5},\cdots,  \alpha_{4+n_p})$
\end{center}
 is a multiindex.

Due to the Cauchy integral formula,  it is easy to see that there exist two positive constants $\rho$ and $\varpi_2$ such that
\begin{center}
$|c_\alpha| \leq \frac{\varpi_2}{\rho^{|\alpha|}}$.
\end{center}

Then
\begin{equation}\label{estimate1}
             \begin{array}{lr}
             |\theta' | \leq  \sum_{|\alpha|=2}\mid c_\alpha\mid~ \parallel u\parallel^{\mid \alpha\mid} \\
              \leq  \sum_{|\alpha|=2}\frac{\varpi_2 \varpi_1^{|\alpha|}}{\rho^{|\alpha|}} e^{\sigma |\alpha| \tau} \\
             =e^{\sigma \tau} \sum_{k=0}\frac{\varpi_2 \varpi_1^{k+2}}{\rho^{k+2}} e^{\sigma k \tau}.
             \end{array}\nonumber
\end{equation}
Obviously, $\sum_{k=0}\frac{\varpi_2 \varpi_1^{k+2}}{\rho^{k+2}} e^{\sigma k \tau}$ is bounded as $\tau \rightarrow -\infty$. Suppose that
\begin{displaymath}
\sum_{k=0}\frac{\varpi_2 \varpi_1^{k+2}}{\rho^{k+2}} e^{\sigma k \tau} \leq \varpi_3
\end{displaymath}
Then
\begin{displaymath}
|\theta' | \leq \varpi_3 e^{\sigma \tau}.
\end{displaymath}
It follows from Cauchy's test for convergence that $\theta(\tau)$  approaches a fixed limit as $\tau \rightarrow -\infty$.

In conclusion, the proof of Theorem \ref{nondegenerateproblem} is completed. That is, we  solve  $PISPW$ for nondegenerate central configurations.

$~~~~~~~~~~~~~~~~~~~~~~~~~~~~~~~~~~~~~~~~~~~~~~~~~~~~~~~~~~~~~~~~~~~~~~~~~~~~~~~~~~~~~~~~~~~~~~~~~~~~~~~~~~~~~~~~~~~~~~~~~~~~~~~~~~~~~~~~~~~~~~~~~~~~~~~~~~~~~~~~~~~\Box$

\subsection{$n_0>0$}\label{subsectionn0>0}
\indent\par

When $n_0>0$, i.e., the central configuration is degenerate,  the problem attributes to the following form:
\begin{problem}\label{problemn_0big0}
given a solution $(q(\tau),\gamma(\tau))$ of the system (\ref{equationzq}) such that
\begin{displaymath}
(q(\tau),\gamma(\tau))\rightarrow 0  ~~~as~~~  \tau \rightarrow -\infty,
\end{displaymath} if $\theta(\tau)$ satisfies (\ref{equationtheta5}) then whether does $\theta(\tau)$  approach a fixed limit as $\tau \rightarrow -\infty$?
\end{problem}

\subsubsection{Preliminaries of Proof}
\indent\par

\textbf{A.} The first trick  of discussing
the above problem is to simplify the system (\ref{equationzq}) by the theory of normal forms (or reduction theorem). In effect, it is to choose an orbit  with  actual calculations on the center manifold of the system (\ref{equationzq}) at the origin to approximate the orbit on the center-unstable manifold.

For ease of notations, set
\begin{displaymath}
  \mathfrak{C}^+=\left(
        \begin{array}{cccc}
          \kappa^{\frac{1}{2}} &   &   &   \\
            & \tilde{{\mu}}_{5+n_0} &  &  \\
          & & \ddots &  \\
          & &  & \tilde{{\mu}}_{4+n_0+n_p} \\
        \end{array}
      \right)
\end{displaymath}
\begin{equation}\label{C-}
\mathfrak{C}^-=
      \begin{pmatrix}
        \tilde{\Lambda}_1 &   &   &   &   &  & &\\
          & \tilde{\Lambda}_2 &   &   &   &  & &\\
          &   & - \frac{\kappa^{\frac{1}{2}}}{4}\mathbb{I}  &   &   &  & &\epsilon \mathbb{I} \\
           &   &   & - \frac{\kappa^{\frac{1}{2}}}{2}\mathbb{I} &0 & 0  & 0& 0 \\
          &   & & 0 & - \frac{\kappa^{\frac{1}{2}}}{2}\mathbb{I} - \tilde{\Lambda}_p & 0  & 0& 0 \\
          &   &  & 0& 0 & - \frac{\kappa^{\frac{1}{2}}}{2}\mathbb{I} - \tilde{\Lambda}_1 & 0  &  0 \\
          &   &  &0 & 0  & 0 & - \frac{\kappa^{\frac{1}{2}}}{2}\mathbb{I} - \tilde{\Lambda}_2 & 0  \\
          &   &  &0 &  0 & 0  & 0 &  - \frac{\kappa^{\frac{1}{2}}}{4}\mathbb{I}
      \end{pmatrix}
\end{equation}
\begin{displaymath}
\begin{array}{c}
  q^0 = ( q_{5}, \cdots, q_{4+n_0})^\top \\
  q^+ = (\gamma, q_{5+n_0}, \cdots, q_{4+n_0+n_p})^\top \\
  q^- = (q_{5+n_0+n_p}, \cdots, q_{2N}, q_{5+2N}, \cdots, q_{4N})^\top\\
   \varphi^0 = ( \varphi_{5}, \cdots, \varphi_{4+n_0})^\top \\
  \varphi^+ = (\varphi_{0}, \varphi_{5+n_0}, \cdots, \varphi_{4+n_0+n_p})^\top\\
  \varphi^- = (\varphi_{5+n_0+n_p},\cdots, \varphi_{2N}, -\varphi_{5},\cdots, -\varphi_{2N})^\top.
\end{array}
\end{displaymath}
Consequently, the system (\ref{equationzq}) can be rewritten  as the following system of equations:
\begin{equation}\label{equationzsim}
\left\{
             \begin{array}{lr}
             {q'^0} =  - \varphi^0(q^0,q^+,q^-), \\
             {q'^+} = \mathfrak{C}^+ q^+  - \varphi^+(q^0,q^+,q^-), \\
             {q'^-} = \mathfrak{C}^- q^-  - \varphi^-(q^0,q^+,q^-).
             \end{array}
\right.
\end{equation}

It follows from  Corollary \ref{normalform5} that we can
introduce a nonlinear substitution of the form
\begin{equation}\label{nonlinear substitution0>0}
\left\{
\begin{array}{lr}
             u_{k} =  q_{k}  ,   &k\in \{5, \cdots, n_0+n_p+4\} \\
             u_{0} =  \gamma ,  & \\
             u_{k} =  q_{k}  - F^{cu}_{k} ( q^0,q^+),   &  k\in \{n_0+n_p+5, \cdots,2N, 2N+5, \cdots,4N\}
             \end{array}
             \right.
\end{equation}
so that the system (\ref{equationzsim})  can be written  as the simpler form below
\begin{equation}\label{simple formsn0>0}
\left\{
             \begin{array}{lr}
             {u'^0} =  - \varphi^0(u^0,u^+,  F^{cu}(u^0,u^+)) + \psi^0(u)u^-, \\
             {u'^+} = \mathfrak{C}^+ u^+  - \varphi^+(u^0,u^+,  F^{cu}(u^0,u^+)) + \psi^+(u)u^-, \\
             {u'^-} = \mathfrak{C}^- u^-  + \psi^-(u)u^-,
             \end{array}
\right.
\end{equation}
where
\begin{displaymath}
\begin{array}{c}
  u^0 = ( u_{5}, \cdots, u_{4+n_0})^\top \\
  u^+ = (u_{0}, u_{5+n_0}, \cdots, u_{4+n_0+n_p})^\top \\
  u^- = (u_{5+n_0+n_p}, \cdots, u_{2N}, u_{5+2N}, \cdots, u_{4N})^\top,
\end{array}
\end{displaymath}
the function
\begin{displaymath}
q^{-}=F^{cu}(q^{0},q^{+})
\end{displaymath}
is a center-unstable manifold  of class $C^l$ such that $F^{cu}$ and all the partial derivative of $F^{cu}$ are vanishing at $(q^{0},q^{+})=0$;
and the functions $\psi^0, \psi^+$ are $C^{l}$-smooth and $ \psi^-$ is  $C^{l-1}$-smooth; in addition, all the functions $\psi^0, \psi^+, \psi^-$ are vanishing at the origin, i.e., $\psi(0)=0$.\\

It is easy to see that $u(\tau)\rightarrow 0$ as $\tau \rightarrow -\infty$, if $(q(\tau),\gamma(\tau))\rightarrow 0$ as $\tau \rightarrow -\infty$. Then we claim that
\begin{displaymath}
u_k(\tau) \equiv 0 ~~~~for ~~~~k\in \{n_0+n_p+5, \cdots, 2N, 2N+5, \cdots, 4N \}.
\end{displaymath}
The proof of the claim is quite standard  and  so is omitted. As a matter of fact, the claim is simply saying that the center-unstable manifold of the system (\ref{simple formsn0>0}) at the origin is  characterized as $u^-=0$.

Therefore, for the Problem \ref{problemn_0big0} the  system (\ref{simple formsn0>0}) reduces to
\begin{equation}\label{simple formsn0>01}
\left\{
             \begin{array}{lr}
             {u'^0} =  - \varphi^0(u^0,u^+,  F^{cu}(u^0,u^+)), \\
             {u'^+} = \mathfrak{C}^+ u^+  - \varphi^+(u^0,u^+,  F^{cu}(u^0,u^+)).
             \end{array}
\right.
\end{equation}

Then it follows from the Theorems  \ref{Center Manifold} and \ref{centerapproximatetheorem} that there exists a solution $\upsilon(\tau)=(\upsilon_5(\tau), \cdots,\upsilon_{n_0+4}(\tau))$ of the following system
\begin{equation}\label{center systemn_0>0}
{u'^0} = -\varphi^0 \left(u^0,F^{c}(u^0),  F^{cu}\left(u^0,F^{c}(u^0) \right) \right),
\end{equation}
 such that
\begin{equation}\label{centerapproximaten_0>0}
\left\{
             \begin{array}{lr}
             {u^0}(\tau) = \upsilon(\tau) + O(e^{\sigma \tau}) \\
             {u^+}(\tau) = F^{c}(\upsilon(\tau))  + O(e^{\sigma \tau})
             \end{array}
\right.
~~~as~ \tau\rightarrow -\infty,
\end{equation}
where the function $u^{+}=F^{c}(u^0)$  is a center manifold  of class $C^l$ such that $F^{c}$ and all the partial derivative of $F^{c}$ are vanishing at $u^0=0$; the abusive $\sigma >0$ is a constant depending only on $\mathfrak{C}^+$.\\

\textbf{B.}
Due to the equation (\ref{equationtheta5}) and the nonlinear substitution (\ref{nonlinear substitution0>0}), now we have
\begin{equation}\label{equationthetan_0>0}
\begin{array}{lr}
             \theta'  =  \left(u^0,u^+, F^{cu}(u^0,u^+)\right)^\top  \mathfrak{P}^\top \begin{pmatrix}\mathbb{I}\\
  0\end{pmatrix}
Q \begin{pmatrix}0&\mathbb{I}\end{pmatrix}
\mathfrak{P}\left(u^0,u^+, F^{cu}(u^0,u^+)\right) \\
             = \sum_{5\leq k\leq n_0+4}\sum_{j=n_0+5}^{n_0+n_p+4}q_{kj}\tilde{{\mu}}_j {u}_k {u_j}+\sum_{j,k=n_0+5}^{n_0+n_p+4}q_{kj}\tilde{{\mu}}_j {u}_k {u_j}  + \cdots
             \end{array}
\end{equation}
where
``$\cdots$" denotes all the  terms which  contain at least one of $F^{cu}_k$ $(k>n_0+n_p+4)$ as a factor.

By Taylor's formula,  the equation (\ref{equationthetan_0>0}) can be rewritten as
\begin{displaymath}
 \theta'=\sum_{5\leq k\leq n_0+4}\sum_{j=n_0+5}^{n_0+n_p+4}q_{kj}\tilde{{\mu}}_j {u}_k {u_j}+\sum_{j,k=n_0+5}^{n_0+n_p+4}q_{kj}\tilde{{\mu}}_j {u}_k {u_j}  + \sum^l_{|\alpha|=3}b_\alpha (u^{0+})^\alpha +o_l(u^{0+}),
\end{displaymath}
where $o_l$ denotes the reminder term which vanishes at the origin along with
the first $l$ derivatives,
\begin{displaymath}
u^{0+}=(u_5,\cdots,u_{4+n_0}, u_{0}, u_{5+n_0}, \cdots,u_{4+n_0+n_p})^\top,
\end{displaymath}
 and
\begin{displaymath}
\alpha=(\alpha_5,   \cdots, \alpha_{4+n_0}, \alpha_{0}, \alpha_{5+n_0},  \alpha_{4+n_0+n_p})
\end{displaymath}
is a multiindex.

Furthermore, it follows from  Theorem \ref{Center-Unstable Manifold} that the coefficients $b_\alpha$  only depend upon
the system (\ref{equationzsim}). In fact, the coefficients $b_\alpha$ are determined by the condition of invariance of the manifolds:
\begin{equation}\label{condition of invariancen0>0}
\begin{array}{c}
\mathfrak{C}^- F^{cu}(q^0,q^+) - \varphi^-(q^0,q^+, F^{cu}) = \\
-\frac{\partial F^{cu}(q^0,q^+)}{\partial q^0} \varphi^0(q^0,q^+, F^{cu}) +\frac{\partial F^{cu}(q^0,q^+)}{\partial q^+}\left[\mathfrak{C}^+ q^+ - \varphi^+(q^0,q^+, F^{cu})\right].
\end{array}
\end{equation}
This is because of the collection of  the eigenvalues of the matrix $\mathfrak{C}^-$
belongs to the Poincar\'{e} domain, then one can uniquely determine a formal power series as the formal solution of above equation.

Similarly, the relation
\begin{equation}\label{invariancencenter systemn_0>01}
\begin{array}{c}
  \mathfrak{C}^+ F^{c}(u^0) - \varphi^+\left(u^0,F^{c}(u^0),F^{cu}(u^0,F^{c}(u^0))\right) = \\
-\frac{\partial F^{c}(u^0)}{\partial u^0} \varphi^0\left(u^0,F^{c}(u^0),F^{cu}(u^0,F^{c}(u^0))\right)
\end{array}
\end{equation}
gives an algorithm for computing the Taylor's coefficients of $F^{c}(u^0)$. For example, the quadratic form in Taylor's formula of $F^{c}_k(u^0)$ ($k=n_0+5, \cdots, n_0+n_p+4$) is
 \begin{equation}\label{quadratic form}
\frac{\sum_{i,j=5}^{4+n_0}a_{ijk}{u}_i {u_j}}{4 \tilde{\mu}_{k} \sqrt{ \mu_{k} + \frac{\kappa}{16}}}.
\end{equation}

Based on the Taylor expansion, the system (\ref{center systemn_0>0}) can be written as
\begin{equation}\label{center systemn_0>01}
 \begin{array}{lr}
             {u'_{k}} = -\frac{\sum_{i,j=5}^{4+n_0}a_{ijk}{u}_i {u_j}}{  \sqrt{  \kappa}} + \cdots+ o_l(u^{0}),  & k\in \{5, \cdots, n_0+4\}.
             \end{array}
\end{equation}

Taking into consideration of (\ref{centerapproximaten_0>0}),  the equation (\ref{equationthetan_0>0}) can further be rewritten as
\begin{displaymath}
 \theta'=\sum_{j=n_0+5}^{n_0+n_p+4}\sum_{k, i,\iota=5}^{4+n_0} \frac{q_{kj}\tilde{{\mu}}_j a_{i\iota j} {\upsilon}_k  {\upsilon}_i {\upsilon_\iota}}{4 \tilde{\mu}_{k} \sqrt{ \mu_{k} + \frac{\kappa}{16}}}+ \cdots +o_l(\upsilon)+O(e^{\sigma \tau}).
\end{displaymath}
Obviously, we can discard the term $O(e^{\sigma \tau})$, and consider
\begin{equation}\label{equationthetan_0>01}
\theta'=\sum_{j=n_0+5}^{n_0+n_p+4}\sum_{k, i,\iota=5}^{4+n_0} \frac{q_{kj}\tilde{{\mu}}_j a_{i\iota j} {\upsilon}_k  {\upsilon}_i {\upsilon_\iota}}{4 \tilde{\mu}_{k} \sqrt{ \mu_{k} + \frac{\kappa}{16}}}+ \cdots +o_l(\upsilon)
\end{equation}
to establish the convergence of $\theta(\tau)$  as $\tau \rightarrow -\infty$.\\

\textbf{C.}
In general, the quadratic forms in (\ref{center systemn_0>01}) may be zero. However, we claim that the Taylor's coefficients of the right side of (\ref{center systemn_0>01}) are not all zero, that is,
\begin{proposition}\label{coefficientisnot0}
For sufficiently large  $l$,  there exists some nonzero Taylor's coefficient in the Taylor expansion of the right side of the system (\ref{center systemn_0>0}):
\begin{equation}
\varphi^0 \left(u^0,  F^{cu}\left(u^0,F^{c}(u^0),F^{c}(u^0) \right) \right).\nonumber
\end{equation}
\end{proposition}
{\bf Proof.} 
It suffices to prove the proposition in the case of formal power series.

If the statement was false, then
\begin{equation}
\varphi^0 \left(u^0,F^{c}(u^0),  F^{cu}\left(u^0,F^{c}(u^0) \right) \right)\equiv 0\nonumber
\end{equation}
in the sense of formal power series.

It follows from (\ref{invariancencenter systemn_0>01}) that
\begin{equation}
\mathfrak{C}^+ F^{c}(u^0) - \varphi^+\left(u^0,F^{c}(u^0),F^{cu}(u^0,F^{c}(u^0))\right) \equiv 0\nonumber
\end{equation}
By (\ref{condition of invariancen0>0}), it turns out that
\begin{equation}
\mathfrak{C}^- F^{cu}(u^0,F^{c}(u^0)) - \varphi^-\left(u^0,F^{c}(u^0),F^{cu}(u^0,F^{c}(u^0))\right) \equiv 0\nonumber
\end{equation}
As a result, $u^0$ and two formal power series $f_1(u^0)=F^{c}(u^0), f_2(u^0)=F^{cu}(u^0,F^{c}(u^0))$ satisfy
\begin{equation}
\left\{
             \begin{array}{lr}
             \varphi^0 \left(u^0,f_1,  f_2 \right)=0, & \\
             \mathfrak{C}^+ f_1 - \varphi^+\left(u^0,f_1,  f_2\right) =0,  & \\
              \mathfrak{C}^- f_2 - \varphi^-\left(u^0,f_1,  f_2\right) =0.  &
             \end{array}
\right.\nonumber
\end{equation}
It is noteworthy that the two equations
\begin{equation}
\left\{
             \begin{array}{lr}
             \mathfrak{C}^+ f_1 - \varphi^+\left(u^0,f_1,  f_2\right) =0,  & \\
              \mathfrak{C}^- f_2 - \varphi^-\left(u^0,f_1,  f_2\right) =0,  &
             \end{array}
\right.\nonumber
\end{equation}
are enough to determine $f_1(u^0), f_2(u^0)$; furthermore, it follows from the analytic version of implicit function theorem that $f_1(u^0), f_2(u^0)$  are more than just formal power series, they are really analytic  functions of $u^0$.

Therefore, we have  infinitely many critical points of the system (\ref{equationzsim}) in a small neighbourhood of the origin. However, we know that the origin is an isolated equilibrium point of the system (\ref{equationzsim}). This leads to a contradiction.

$~~~~~~~~~~~~~~~~~~~~~~~~~~~~~~~~~~~~~~~~~~~~~~~~~~~~~~~~~~~~~~~~~~~~~~~~~~~~~~~~~~~~~~~~~~~~~~~~~~~~~~~~~~~~~~~~~~~~~~~~~~~~~~~~~~~~~~~~~~~~~~~~~~~~~~~~~~~~~~~~~~~\Box$

\subsubsection{$n_0=1$}
\indent\par

By solving the Problem \ref{problemn_0big0} for $n_0= 1$,  we completely solve $PISPW$ for $n_0 \leq 1$ in this subsection. That is
\begin{theorem}\label{degenerateproblemn_0=1}
For given  masses of the $N$-body problem, if all the central configurations have  degree of degeneracy less than or equal to one, then the normalized configuration of any total collision orbit of the given   $N$-body problem
approaches a certain central configuration as time $t$ approaches the collision instant.
\end{theorem}
{\bf Proof.} 
 Consider the Problem \ref{problemn_0big0} for $n_0= 1$.

As a matter of notational convenience, set $u^0=u_5 = w$. Then the system (\ref{center systemn_0>01}) becomes
\begin{equation}
 \begin{array}{lr}
             w' = c_{2}w^2 + \cdots+ c_{l}w^l+ o_l(w).
             \end{array}\nonumber
\end{equation}
By the Proposition \ref{coefficientisnot0},  suppose
\begin{displaymath}
c_{2}=\cdots=c_{m-1}=0, c_{m}\neq 0, 2\leq m<l,
\end{displaymath}
thus
\begin{equation}\label{center systemn_0=1}
 \begin{array}{lr}
             w' = c_{m}w^m + \cdots+ c_{l}w^l+ o_l(w).
             \end{array}
\end{equation}

If $w(\tau) =0$ for some $\tau$, then $w(\tau)\equiv0$, and (\ref{equationthetan_0>01}) becomes $\theta'(\tau)=0$, therefore $\theta(\tau)$ is obviously convergent as $\tau \rightarrow -\infty$.

So we consider $w(\tau) \neq 0$ as $\tau \rightarrow -\infty$. Without loss of generality, we assume $w(\tau) > 0$, i.e., $w(\tau)\rightarrow 0+$ as $\tau \rightarrow -\infty$.

According to L'H\'{o}pital's rule, it follows from (\ref{center systemn_0=1}) that
\begin{displaymath}
\lim_{\tau \rightarrow -\infty} \frac{1}{\tau w^{m-1}}= (1-m)c_{m},
\end{displaymath}
or
\begin{displaymath}
 w= \left(\frac{1}{ (1-m)c_{m}\tau}\right)^\frac{1}{m-1} +o((\frac{1}{-\tau})^\frac{1}{m-1}) ~~~~~~~~~as~  \tau \rightarrow -\infty.
\end{displaymath}

We cannot prove the convergence of $\theta(\tau)$ by (\ref{equationthetan_0>01}), but a proof based on original (\ref{equationtheta4}) will be given right now.

It is simple to show that for every $j=5, \cdots, 2N$,
\begin{equation}
 \begin{array}{c}
             Z_j = c_{j,2}w^2 + \cdots+ c_{j,l}w^l+ o_l(w).
             \end{array}\nonumber
\end{equation}
If $c_{j,2}, \cdots, c_{j,l}$ are all zero, then  it is easy to show that
\begin{equation*}
    |Z_j|=O((\frac{1}{ -\tau})^{1+\frac{2}{l-2}}),
\end{equation*}
it follows that the improper integral $\int^{\tau_0}_{-\infty}|Z_j|d\tau$ converges.
If at least one of $c_{j,2}, \cdots, c_{j,l}$ is not zero, then $Z_j$  takes the following form
\begin{displaymath}
 Z_j = \tilde{c}_j(\frac{1}{-\tau})^\frac{1}{m_j} +o((\frac{1}{-\tau})^\frac{1}{m_j}).
\end{displaymath}
So  $Z_j$    converges to zero monotonically. Then it follows from $z'_j=Z_j$ that the improper integral $\int^{\tau_0}_{-\infty}|Z_j|d\tau$ also converges.

As a result,  the improper integral $\int^{\tau_0}_{-\infty}|Z_j|d\tau$  converges for every $j=5, \cdots, 2N$.
 Consequently, the improper integral
\begin{equation}
|\int^{\tau_0}_{-\infty}\theta'(\tau)d\tau|  \leq \int^{\tau_0}_{-\infty}\sum_{j,k=5}^{2N}\mid q_{kj}\mid \mid {Z_j}\mid \mid{z}_k\mid d\tau \nonumber
\end{equation}
also converges.

Therefore $\theta(\tau)$ is obviously convergent as $\tau \rightarrow -\infty$.

In conclusion, the proof of Theorem \ref{degenerateproblemn_0=1} is completed.

$~~~~~~~~~~~~~~~~~~~~~~~~~~~~~~~~~~~~~~~~~~~~~~~~~~~~~~~~~~~~~~~~~~~~~~~~~~~~~~~~~~~~~~~~~~~~~~~~~~~~~~~~~~~~~~~~~~~~~~~~~~~~~~~~~~~~~~~~~~~~~~~~~~~~~~~~~~~~~~~~~~~\Box$

\subsubsection{$n_0=2$}
\indent\par

Due to the intrinsic difficulty of  degenerate or nonhyperbolic differential equations,  the difficulty of the problem  increases rapidly for $n_0 \geq 2$. Indeed, we cannot completely resolve $PISPW$ even in the case of $n_0=2$. The main difficulty is from  estimating the rate of an orbit on a center-unstable manifold tending to the  equilibrium point.

In particular, we could not apply a similar method as in the above subsection of $n_0=1$, because we can not prove that every $Z_j$ ($j=5, \cdots, 2N$)  monotonically  converges to zero generally. Indeed, in general, $Z_j$ can wavily converges to zero as in Example 2.

Naturally, one hopes that researches on central configurations would help to resolve  $PISPW$. Unfortunately, the problem of central configurations is also very difficult,  as indicated by the researches on this topic in the past decades. We can  only get and utilize some special results of central configurations.

We consider the Problem \ref{problemn_0big0} for $n_0= 2$ in this subsection.  First, we give a criterion for the case of $n_0=2$; then we  give an answer to $PISPW$  for all  known degenerate central configurations of four-body. Therefore, for almost every choice
of the masses of the four-body problem,  $PISPW$ is answered.
\begin{theorem}\label{degenerateproblemn_0=2}
For given  masses of the $N$-body problem, if all the central configurations have  degree of degeneracy less than or equal to two, and central configurations with    degree of degeneracy two satisfy the condition (\ref{discriminate}), then the normalized configuration of any total collision orbit of the given   $N$-body problem
approaches a certain central configuration as time $t$ approaches the collision instant.
\end{theorem}
\begin{remark}
For collision orbits of the general four-body problem,  $PISPW$ is waiting for an answer in very few cases corresponding to   central configurations with  degree of degeneracy two. Recall that  the  masses corresponding to degenerate central configurations form a proper algebraic subset of the mass space for the four-body problem  \cite{moeckel2014lectures}. Furthermore, following a crude dimension count,  almost all of degenerate central configurations have one degree of degeneracy, and  the  masses corresponding to   central configurations with  degree of degeneracy two should form a subset of the mass space consisting of finite points.
\end{remark}
As an example, we have the following result:
\begin{corollary}\label{lastcorollary}
If the degenerate equilateral central configurations   are  the only  central configurations with  degree of degeneracy two for the four-body problem, then the normalized configuration of any total collision orbit of the  four-body problem
approaches a certain central configuration as time $t$ approaches the collision instant.
\end{corollary}

{\bf Proof of Theorem \ref{degenerateproblemn_0=2}:} 

Consider the Problem \ref{problemn_0big0} for $n_0= 2$.

As a matter of notational convenience, set $u^0=(u_5,u_6) = (\zeta,\eta)$. Then the system (\ref{center systemn_0>01}) becomes
\begin{equation}\label{center systemn_0=2}
\left\{
             \begin{array}{lr}
             {\zeta'} = c_1 \zeta^2 + 2 c_2 \zeta\eta+ c_3 \eta^2 + o(\zeta^2+\eta^2) \\
             {\eta'} = c_2 \zeta^2 + 2 c_3 \zeta\eta+ c_4 \eta^2 + o(\zeta^2+\eta^2)
             \end{array}
\right.
\end{equation}
where \begin{displaymath}
c_1=-\frac{a_{555}}{  \sqrt{  \kappa}}, ~~~c_2=-\frac{a_{556}}{  \sqrt{  \kappa}}, ~~~c_3=-\frac{a_{566}}{  \sqrt{  \kappa}}, ~~~c_4=-\frac{a_{666}}{  \sqrt{  \kappa}}.
\end{displaymath}

Let us introduce polar coordinates as in Appendix \ref{Plane Equilibrium Points}
\begin{displaymath}
\zeta=\rho \cos \vartheta, \eta=\rho \sin \vartheta
\end{displaymath}
to transform the system (\ref{center systemn_0=2})
into
\begin{equation}
\left\{
             \begin{array}{lr}
             \rho' = \rho^{2}\Phi(\vartheta) + o (\rho^{2}), \\
             \vartheta' = \rho\Psi(\vartheta) + o (\rho),
             \end{array}
\right.\nonumber
\end{equation}
where
\begin{equation}
\left\{
             \begin{array}{lr}
             \Phi(\vartheta) = c_1 \cos^3\vartheta + 3 c_2 \cos^2\vartheta\sin \vartheta + 3c_3 \cos\vartheta\sin^2 \vartheta +c_4 \sin^3\vartheta, \\
             \Psi(\vartheta) = c_2 \cos^3\vartheta + (2c_3- c_1) \cos^2\vartheta\sin \vartheta + (c_4-2c_2) \cos\vartheta\sin^2 \vartheta - c_3 \sin^3\vartheta;
             \end{array}
\right.\nonumber
\end{equation}
thus $\Phi,\Psi$ are homogeneous polynomials  of degree $3$ in $\cos \vartheta, \sin \vartheta$.

To prove the theorem, we need two lemmas.
\begin{lemma}\label{varthetamainn_0=2}
Assume $c_1,c_2,c_3,c_4$ are not all zero.
Then
\begin{equation}
              \vartheta_0= \lim_{\tau\rightarrow -\infty} \vartheta(\tau)~~~~~~exists ~(and~ is~ finite)
             \nonumber
\end{equation}
and $\Psi(\vartheta_0)= 0$.
\end{lemma}
\begin{lemma}\label{rhomainn_0=2}
Assume
\begin{displaymath}
c^2_1 c^2_4-6 c_1 c_2 c_3 c_4+4 c_1 c^3_3+4 c^3_2 c_4-3 c^2_2 c^2_3\neq 0.
\end{displaymath}
Then
\begin{displaymath}
 \rho = -\frac{1}{\tau \Phi(\vartheta_0)} +o(\frac{1}{\tau}),
\end{displaymath}
where $\vartheta_0$ is just the one in Lemma \ref{varthetamainn_0=2}.
\end{lemma}

Note that  (\ref{equationthetan_0>01}) becomes
\begin{equation}\label{equationthetan_0=2}
\theta'=P_3(\zeta,\eta) +o(\rho^3), \nonumber
\end{equation}
where $P_3(\zeta,\eta)$ is homogeneous polynomials of degree $3$.

So $|\theta' | \leq  \sigma \rho^3$ for some positive number $\sigma$ and sufficiently small $\rho$.
Then it follows from  Lemma \ref{rhomainn_0=2} and  Cauchy's test for convergence that $\theta(\tau)$  approaches a fixed limit as $\tau \rightarrow -\infty$,
provided that
\begin{displaymath}
c^2_1 c^2_4-6 c_1 c_2 c_3 c_4+4 c_1 c^3_3+4 c^3_2 c_4-3 c^2_2 c^2_3\neq 0,
\end{displaymath}
i.e.,
\begin{equation}\label{discriminate}
a^2_{555} a^2_{666}-6 a_{555} a_{556} a_{566} a_{666}+4 a_{555} a^3_{566}+4 a^3_{556} a_{666}-3 a^2_{556} a^2_{566}\neq 0.
\end{equation}

In conclusion, the proof of Theorem \ref{degenerateproblemn_0=2} is completed.

$~~~~~~~~~~~~~~~~~~~~~~~~~~~~~~~~~~~~~~~~~~~~~~~~~~~~~~~~~~~~~~~~~~~~~~~~~~~~~~~~~~~~~~~~~~~~~~~~~~~~~~~~~~~~~~~~~~~~~~~~~~~~~~~~~~~~~~~~~~~~~~~~~~~~~~~~~~~~~~~~~~~\Box$

{\bf Proof of Lemma \ref{varthetamainn_0=2}:}

Since the origin is an isolated equilibrium point of the system (\ref{equationzsim}), and the system (\ref{center systemn_0>01}) is just the restriction of the system (\ref{equationzsim}) to the center manifold $u^{+}=F^{c}(u^0)$,  it follows that the origin is also an isolated equilibrium point of the system (\ref{center systemn_0>01}).

A routine computation gives rise to
\begin{equation}
\Psi(\vartheta) = \frac{c_2+c_4}{4} \cos\vartheta + \frac{-c_1- c_3}{4} \sin \vartheta + \frac{3c_2-c_4}{4} \cos{3\vartheta} + \frac{3c_3 -c_1}{4}\sin{3\vartheta}.\nonumber
\end{equation}
As a result, $\Psi(\vartheta) \equiv 0$ if and only if
\begin{equation}
\left\{
             \begin{array}{lr}
             c_2+c_4 =0,  & \\
             -c_1- c_3 =0,  & \\
             3c_2-c_4 =0,  & \\
             3c_3 -c_1 =0,  &
             \end{array}
\right.\nonumber
\end{equation}
or
\begin{equation}
c_1=c_2=c_3=c_4=0.\nonumber
\end{equation}

The lemma is now a direct consequence of the Theorem \ref{varthetamain}.

$~~~~~~~~~~~~~~~~~~~~~~~~~~~~~~~~~~~~~~~~~~~~~~~~~~~~~~~~~~~~~~~~~~~~~~~~~~~~~~~~~~~~~~~~~~~~~~~~~~~~~~~~~~~~~~~~~~~~~~~~~~~~~~~~~~~~~~~~~~~~~~~~~~~~~~~~~~~~~~~~~~~\Box$

{\bf Proof of Lemma \ref{rhomainn_0=2}:}

Obviously, $c_1,c_2,c_3,c_4$ are not all zero. Therefore
\begin{equation}\label{vartheta11n_0=2}
    \begin{array}{lr}
              \vartheta_0= \lim_{\tau\rightarrow -\infty} \vartheta(\tau)~~~~~~exists ~(and~ is~ finite)
             \end{array}\nonumber
\end{equation}
and
\begin{center}
$\Psi(\vartheta_0)= 0$.
\end{center}

Furthermore, we claim that $\Phi(\vartheta_0) \neq 0$. If otherwise, then it is easy to see that
\begin{equation}
\left\{
             \begin{array}{lr}
             c_1 \zeta^2 + 2 c_2 \zeta\eta+ c_3 \eta^2 =0 \\
             c_2 \zeta^2 + 2 c_3 \zeta\eta+ c_4 \eta^2 =0
             \end{array}
\right.\nonumber
\end{equation}
for $\zeta=\cos\vartheta_0, \eta=\sin\vartheta_0$.

If $\sin\vartheta_0\neq0$, we conclude that
\begin{equation}
\left\{
             \begin{array}{lr}
             c_1 w^2 + 2 c_2 w+ c_3  =0 \\
             c_2 w^2 + 2 c_3 w+ c_4  =0
             \end{array}
\right.\nonumber
\end{equation}
for $w=\frac{\cos\vartheta_0}{\sin\vartheta_0}$. Thus the resultant
\begin{displaymath}
\left|
  \begin{array}{cccc}
    c_1 & 2c_2 & c_3 & 0 \\
    0 & c_1 & 2c_2 & c_3 \\
    c_2 & 2c_3 & c_4 & 0 \\
    0 & c_2 & 2c_3 & c_4 \\
  \end{array}
\right|
\end{displaymath}
is zero, i.e.,
\begin{displaymath}
c^2_1 c^2_4-6 c_1 c_2 c_3 c_4+4 c_1 c^3_3+4 c^3_2 c_4-3 c^2_2 c^2_3 =0.
\end{displaymath}
This is contrary to the assumption of the lemma.

If $\sin\vartheta_0=0$, then $\cos\vartheta_0\neq0$. Using the same argument as above, we can always obtain $\Phi(\vartheta_0) \neq 0$.

The lemma is now a direct consequence of Theorem \ref{rhoamain}.

$~~~~~~~~~~~~~~~~~~~~~~~~~~~~~~~~~~~~~~~~~~~~~~~~~~~~~~~~~~~~~~~~~~~~~~~~~~~~~~~~~~~~~~~~~~~~~~~~~~~~~~~~~~~~~~~~~~~~~~~~~~~~~~~~~~~~~~~~~~~~~~~~~~~~~~~~~~~~~~~~~~~\Box$

{\bf Proof of Corollary \ref{lastcorollary}:}

We only need to test and verify the condition (\ref{discriminate}) for the degenerate equilateral central configurations.

As in Appendix \ref{Equilateral Central Configurations}, we consider
\begin{displaymath}
m_1=m_2=m_3=1,m_4=\frac{81+64\sqrt{3}}{249}.
\end{displaymath}
Recall that
\begin{displaymath}
a_{ijk}=d^3 \mathcal{U}|_{\widehat{\mathcal{E}}_3}(\widehat{\mathcal{E}}_i,\widehat{\mathcal{E}}_j,\widehat{\mathcal{E}}_k), ~~~~~~~~~~~~~ i,j,k\in\{5,6\},
\end{displaymath}
and
\begin{equation}
\begin{array}{c}
  {\mathcal{E}}_3= \mathbf{r}=(-\frac{\sqrt{3}}{2},-\frac{1}{2},\frac{\sqrt{3}}{2},-\frac{1}{2},0,1,0,0)^\top, \\
  {\mathcal{E}}_4= \mathbf{i}\mathbf{r}=(\frac{1}{2},-\frac{\sqrt{3}}{2},\frac{1}{2},\frac{\sqrt{3}}{2},-1,0,0,0)^\top,\\
  \mathcal{E}_5=\left(\frac{64 \sqrt{3}+81}{498} ,-\frac{741 \sqrt{3}+908}{1494},\frac{64 \sqrt{3}+81}{498} ,\frac{741 \sqrt{3}+908}{1494},0,0,-1,0\right)^\top, \\
  \mathcal{E}_6 = \left(\frac{165 \sqrt{3}+179}{747} ,-\frac{371 \sqrt{3}+738}{2241},-\frac{165 \sqrt{3}+179}{747} ,-\frac{371 \sqrt{3}+738}{2241},0,\frac{2 \sqrt{3}+9}{27} ,0,1\right)^\top.
\end{array}\nonumber
\end{equation}

Some tedious computation yields
\begin{displaymath}
a_{555}=0,
\end{displaymath}
\begin{displaymath}
a_{556}=-6630331032 \sqrt{\frac{2}{13129701006956661 \sqrt{3}+22740709543896356}},
\end{displaymath}
\begin{displaymath}
a_{566}=0,
\end{displaymath}
\begin{displaymath}
a_{666}=3269394 \sqrt{\frac{2}{6312834009 \sqrt{3}+10926270656}}.
\end{displaymath}

Obviously,
\begin{equation}
a^2_{555} a^2_{666}-6 a_{555} a_{556} a_{566} a_{666}+4 a_{555} a^3_{566}+4 a^3_{556} a_{666}-3 a^2_{556} a^2_{566}\neq 0.\nonumber
\end{equation}

The proof is completed.

$~~~~~~~~~~~~~~~~~~~~~~~~~~~~~~~~~~~~~~~~~~~~~~~~~~~~~~~~~~~~~~~~~~~~~~~~~~~~~~~~~~~~~~~~~~~~~~~~~~~~~~~~~~~~~~~~~~~~~~~~~~~~~~~~~~~~~~~~~~~~~~~~~~~~~~~~~~~~~~~~~~~\Box$

\section{Manifold of  Collision Orbits}
\indent\par
Based on the work about $PISPW$,  we can consider now the manifold of all the collision orbits or the set of initial conditions leading to total collisions locally. We also  have to divide the discussion into several cases according to the value of $n_0$.

\subsection{$n_0=0$}
\indent\par
First, let us consider the case of $n_0=0$, i.e., the central configuration is nondegenerate.
\begin{theorem}\label{manifoldofallthecollisionorbitsn_0=0}
For the planar $N$-body problem, the manifold of all collision orbits corresponding to a fixed nondegenerate central configuration  is a real
analytic manifold of  dimension $n_p+8$  in a neighbourhood of the collision instant. Therefore the set of initial conditions leading to total collisions is locally a finite union of  real analytic submanifold in the neighbourhood of collision instant, here  the
dimensions of the submanifolds depend upon the index of the limiting central configuration and are at most $2N+4$.
\end{theorem}
{\bf Proof.} 
The equations of motion now become:
\begin{equation}\label{Manifoldequation1}
\left\{
             \begin{array}{lr}
             q'_k= \tilde{{\mu}}_k {q}_k - \varphi_k(q, \gamma),  & k\in \{ 5, \cdots, 2N-n_3\} \\
             q'_{2N+k}= (- \frac{\kappa^{\frac{1}{2}}}{2}-\tilde{\mu}_k) {q}_{2N+k} + \varphi_k(q, \gamma),  & k\in \{ 5, \cdots, 2N-n_3\}\\
             q'_k= - \frac{\kappa^{\frac{1}{2}}}{4} {q}_k + \epsilon {q}_{2N+k} - \varphi_k(q, \gamma),  & k\in \{2N-n_3 +1, \cdots, 2N\} \\
             q'_{2N+k}= - \frac{\kappa^{\frac{1}{2}}}{4} {q}_{2N+k} + \varphi_k(q, \gamma),  & k\in \{2N-n_3 +1, \cdots, 2N\}\\
             \gamma'  = \kappa^{\frac{1}{2}}\gamma - \varphi_0(q, \gamma),&
             \end{array}
\right.
\end{equation}
and
\begin{equation}\label{Manifoldequation2}
\left\{
             \begin{array}{lr}
             \theta' = \sum_{5\leq k\leq n_0+4}\sum_{j=n_0+5}^{n_0+n_p+4}q_{kj}\tilde{{\mu}}_j {q}_k {q_j}+ \cdots,&\\
             r'  =  r (\kappa^{\frac{1}{2}} + \gamma).
             \end{array}
\right.
\end{equation}

Since  the origin  is a hyperbolic  equilibrium
point of the system  (\ref{Manifoldequation1}), it follows from  Corollary \ref{normalform3} that we can
introduce a nonlinear substitution of the form
\begin{equation}
\left\{
\begin{array}{lr}
             u_{0} =  \gamma  - F_{0} (q_{5+n_p}, \cdots,q_{2N},q_{2N+5}, \cdots,q_{4N}), &\\
             u_{k} =  q_{k}  - F_{k} (q_{5+n_p}, \cdots,q_{2N},q_{2N+5}, \cdots,q_{4N}),   ~~~~~~~~~~~~~~k\in \{5, \cdots, 4+n_p\} , &\\
             u_{k} =  q_{k}  - F_{k} (\gamma, q_{5}, \cdots,q_{4+n_p}),   ~~~~k\in \{5+n_p, \cdots, 2N\}\bigcup\{2N+5, \cdots, 4N\}, &
             \end{array}
             \right.\nonumber
\end{equation}
so that the system  (\ref{Manifoldequation1}) is transformed into a simpler form below
\begin{equation}
\left\{
             \begin{array}{lr}
             {u'^+} = \mathfrak{C}^+ u^+  +  \psi^+(u)u^+, \\
             {u'^-} = \mathfrak{C}^- u^-  +  \psi^-(u)u^-,
             \end{array}
\right.\nonumber
\end{equation}
where  the functions $F_{k}$  are power-series starting with quadratic terms and convergent
for small independent variables;
\begin{displaymath}
\begin{array}{c}
  u^+ = (u_{0}, u_{5}, \cdots, u_{4+n_p})^\top, \\
  u^- = (u_{5+n_p}, \cdots, u_{2N}, u_{5+2N}, \cdots, u_{4N})^\top,\\
  \mathfrak{C}^+=\left(
        \begin{array}{cccc}
          \kappa^{\frac{1}{2}} &   &   &   \\
            & \tilde{{\mu}}_5 &  &  \\
          & & \ddots &  \\
          & &  & \tilde{{\mu}}_{4+n_p} \\
        \end{array}
      \right),
\end{array}
\end{displaymath}
 $\mathfrak{C}^-$ is just the same as (\ref{C-}),
and $\psi^+, \psi^-$ are two matrix-valued functions, whose elements are power-series in the $4N-7$ independent
 variables $u_{0}, u_{5}, \cdots, u_{4N}$ starting with linear terms and convergent
for small $u_{0}, u_{5}, \cdots, u_{4N}$.

As a result,  every collision orbit $(\gamma,q,r,\theta)$ of the system of equations (\ref{Manifoldequation1}) and (\ref{Manifoldequation2}) satisfies
\begin{equation}
\left\{
\begin{array}{lr}
             \gamma = u_{0}   + G_{0} (u^+), &\\
             q_{k} =   u_{k} + G_{k} (u^+),   ~~~~~~~~~~~~~~k\in \{5, \cdots, 4+n_p\} , &\\
             q_{k} =  G_{k} (u^+),   ~~~~k\in \{5+n_p, \cdots, 2N\}\bigcup\{2N+5, \cdots, 4N\}, &
             \end{array}
             \right.
\end{equation}
and
\begin{equation}\label{Manifoldsimple formsn0=01}
{u'^+} = \mathfrak{C}^+ u^+  +  \psi^+_0(u^+)u^+
\end{equation}
where  the $G_{k}$  are power-series starting with quadratic terms and convergent
for small $u^+=(u_{0},u_{5},\cdots,u_{4+n_p})^\top$ and $\psi^+_0(u^+)$ is power-series starting with linear terms and convergent
for small $u^+$.

Due to the  facts (\ref{conjugation}),  an argument similar to the one used in  \cite{C1967Lectures} shows that  the coefficients in the power-series $\psi^+_0$ and $u_k $ ($k\in \{0,5, \cdots, 4+n_p\}$) are real numbers.

We claim that  a general solution of (\ref{Manifoldsimple formsn0=01}) satisfying  $u^+(\tau)\rightarrow 0$ as $\tau \rightarrow -\infty$ contains exactly $n_p+1$ independent real parameters.

It suffices to prove the claim in the case that
\begin{displaymath}
0<\kappa^{\frac{1}{2}}\leq \tilde{{\mu}}_5\leq \tilde{{\mu}}_6 \leq\cdots\leq \tilde{{\mu}}_{4+n_p}
\end{displaymath}
or even
\begin{equation}\label{tezhenzhishunxv}
0<\kappa^{\frac{1}{2}}< \tilde{{\mu}}_5< \tilde{{\mu}}_6 <\cdots< \tilde{{\mu}}_{4+n_p}.
\end{equation}

Then it follows from Theorem \ref{normalformpd}   that we can
introduce a nonlinear substitution of the form
\begin{equation}\label{nonlinear substitutionpd1}
\tilde{u}_k =  u_k - \tilde{F}(u^+),  ~~~k\in \{0, 5, \cdots, 4+n_p\}\nonumber
\end{equation}
so that the system (\ref{Manifoldsimple formsn0=01}) can be reduced to the following simple form
\begin{equation}\label{simple formspd1}
\tilde{u}' = \mathfrak{C}^+ \tilde{u}+ R(\tilde{u}) ,
\end{equation}
where  $\tilde{F}$  is power-series starting with quadratic terms and convergent
for small $u^+$, and
 $R(\tilde{u})$ is a finite-order polynomial composed by resonant monomials.

Note that it follows from  (\ref{tezhenzhishunxv}) that the component $R_k$ of $R(\tilde{u})$ is a polynomial only in the first $k-1$
variables in $\tilde{u}_{0},\tilde{u}_{5},\cdots,\tilde{u}_{4+n_p}$. So finally we obtain the following system of differential equations
for $\tilde{u}_{0},\tilde{u}_{5},\cdots,\tilde{u}_{4+n_p}$:
\begin{equation}
\left\{
             \begin{array}{lr}
             \tilde{u}'_0 = \kappa^{\frac{1}{2}} \tilde{u}_0,  \\
             \tilde{u}'_5 = \tilde{{\mu}}_5 \tilde{u}_5 +R_5(\tilde{u}_0),\\
             \cdots \\
             \tilde{u}'_{4+n_p} = \tilde{{\mu}}_{4+n_p} \tilde{u}_{4+n_p} +R_{4+n_p}(\tilde{u}_0, \tilde{u}_0, \cdots, \tilde{u}_{3+n_p}).
             \end{array}
\right.
\nonumber
\end{equation}

Then we can determine the general solution of (\ref{simple formspd1}) inductively. Indeed, an easy induction gives that
\begin{equation}\label{polinomal}
\begin{array}{lr}
             \tilde{u}_0 =  c_0 e^{\kappa^{\frac{1}{2}} \tau}, &  \\
             \tilde{u}_k =  \left(c_k+P_k(c_0,c_5,\cdots,c_{k-1}, \tau)\right)e^{\tilde{{\mu}}_k \tau},  &k\in \{ 5, \cdots, 4+n_p\},
             \end{array}
\end{equation}
where  $c_0,c_5,\cdots,c_{4+n_p}$ are  constants of integration and are uniquely determined by the initial value $\tilde{u}(0)$,
 $P_k$ is  a polynomial in $c_0,c_5,\cdots,c_{k-1}, \tau$ and contains a positive
power of $\tau$. Obviously, $\tilde{u}(\tau)\rightarrow 0$ as $\tau \rightarrow -\infty$.

This completes the proof of the above claim.

Therefore, every collision orbit $(r,\theta,\gamma,q)$ contains exactly $n_p+3$ independent real parameters $c_0,c_5,\cdots,c_{4+n_p}$ and $c_\theta,c_r$, here $c_\theta,c_r$ are
constants of integration corresponding to the system (\ref{Manifoldequation2}). As a result, the coordinates $r, \Upsilon, \theta, z, Z$ or the position vector and velocity vector of collision orbit $\mathbf{r}\in \mathcal{X}$ contains exactly $n_p+3$ independent real parameters too. Since we have
assumed that  the center of mass $\mathbf{r}_c$ remains at the origin, and the
center of mass integrals involve four real parameters; in addition, the points
on a collision orbit are parameterized by the real variable $t$. Thus it follows that the position vector and velocity vector of  a collision orbit $\mathbf{r}\in \mathbb{R}^{2N}$ have $n_p+8$ independent real parameters.  Because the coordinate functions are regular analytic functions of these
parameters, we conclude that the manifold of all collision orbits corresponding to a fixed nondegenerate central configuration  is a real
analytic manifold of  dimension $n_p+8$  in a neighbourhood of the collision instant.

Here, for a fixed nondegenerate central configuration, $n_p$  satisfies $N-2 \leq n_p \leq 2N-4$. Hence $n_p+8 \leq 2N+4 \leq 4N -2$ for $N\geq3$. As a result, the set of initial conditions leading to total collisions is locally a finite union of  real analytic submanifold in the neighbourhood of collision instant, and  the
dimensions of the submanifolds depend upon the index of the limiting central configuration and are at most $2N+4$ for the $N$-body problem.

The proof is completed.

$~~~~~~~~~~~~~~~~~~~~~~~~~~~~~~~~~~~~~~~~~~~~~~~~~~~~~~~~~~~~~~~~~~~~~~~~~~~~~~~~~~~~~~~~~~~~~~~~~~~~~~~~~~~~~~~~~~~~~~~~~~~~~~~~~~~~~~~~~~~~~~~~~~~~~~~~~~~~~~~~~~~\Box$

\subsection{$n_0>0$}
\indent\par

Let us consider the case of $n_0>0$, i.e., the central configuration is degenerate.
\begin{theorem}
The set of initial conditions leading to total collisions is locally a finite union of real  submanifold in the neighbourhood of collision instant, and  the
dimensions of the submanifolds depend upon the index of the limiting central configuration and are at most $2N+4$ for the $N$-body problem.
\end{theorem}
\begin{remark}
Since $2N+4 \leq 4N -2$ for $N\geq3$, we get the result that the set of initial conditions leading to total collisions  has zero measure in the neighbourhood of collision instant.  According to the invariant set $\mathcal{J} \equiv 0$ is of dimension $4N -1$,  it follows that  the set of initial conditions leading to total collisions  has zero measure  even when restricted to the invariant set $\mathcal{J} \equiv 0$. However, let us quote a remark by Siegel in \cite{C1967Lectures}:``We remark that since our solutions are described only near $t=0$, the above description of the collision orbits is purely local. It is not possible
to describe the manifold of collision orbits in the large, that is, for all $t$,
by our method."
\end{remark}
{\bf Proof.} 
The equations of motion now become:
\begin{equation}\label{Manifoldequationn_0>0}
\left\{
             \begin{array}{lr}
             q'_k=  - \varphi_k(q, \gamma),  & k\in \{5, \cdots, n_0+4\} \\
             q'_{2N+k}= - \frac{\kappa^{\frac{1}{2}}}{2} {q}_{2N+k} + \varphi_k(q, \gamma),  & k\in \{5, \cdots, n_0+4\}\\
             q'_k= \tilde{{\mu}}_k {q}_k - \varphi_k(q, \gamma),  & k\in \{ n_0+5, \cdots, 2N-n_3\} \\
             q'_{2N+k}= (- \frac{\kappa^{\frac{1}{2}}}{2}-\tilde{\mu}_k) {q}_{2N+k} + \varphi_k(q, \gamma),  & k\in \{ n_0+5, \cdots, 2N-n_3\}\\
             q'_k= - \frac{\kappa^{\frac{1}{2}}}{4} {q}_k + \epsilon {q}_{2N+k} - \varphi_k(q, \gamma),  & k\in \{2N-n_3 +1, \cdots, 2N\} \\
             q'_{2N+k}= - \frac{\kappa^{\frac{1}{2}}}{4} {q}_{2N+k} + \varphi_k(q, \gamma),  & k\in \{2N-n_3 +1, \cdots, 2N\}\\
             \gamma'  = \kappa^{\frac{1}{2}}\gamma - \varphi_0(q, \gamma),&\\
             \theta' = \sum_{5\leq k\leq n_0+4}\sum_{j=n_0+5}^{n_0+n_p+4}q_{kj}\tilde{{\mu}}_j {q}_k {q_j}+ \cdots,&\\
             r'  =  r (\kappa^{\frac{1}{2}} + \gamma).
             \end{array}
\right.\nonumber
\end{equation}

Although we cannot completely resolve $PISPW$, we shall give a measure of the set of initial conditions leading to total collisions.

 Recall that in Subsection \ref{subsectionn0>0} it is shown that every collision orbit $(\gamma,q,r,\theta)$ satisfies
\begin{equation}
\left\{
\begin{array}{lr}
             \gamma = u_{0}  , &\\
             q_{k} =   u_{k},   ~~~~~~~~~~~~~~k\in \{5, \cdots, n_0+n_p+4\} , &\\
             q^{-}=F^{cu}(u^{0},u^{+}), &
             \end{array}
             \right.\nonumber
\end{equation}
and
\begin{equation}\label{Manifoldsimple formsn0>0}
\left\{
             \begin{array}{lr}
             {u'^0} =  - \varphi^0(u^0,u^+,  F^{cu}(u^0,u^+)), \\
             {u'^+} = \mathfrak{C}^+ u^+  - \varphi^+(u^0,u^+,  F^{cu}(u^0,u^+)).
             \end{array}
\right.\nonumber
\end{equation}
where $F^{cu}$ is a center-unstable manifold  of class $C^l$.

Therefore, the coordinates $r,\gamma,q$  depend on at most $n_0+n_p+2$ independent real parameters smoothly. Taking into consideration of  the
center of mass integrals and the variables $\theta, t$, it follows that the position vector and velocity vector of  a collision orbit $\mathbf{r}\in \mathbb{R}^{2N}$ have at most $n_0+n_p+8$ independent real parameters.  As a result, we conclude that the manifold of all collision orbits corresponding to a fixed degenerate central configuration  is a real
smooth manifold of dimension no more than $n_0+n_p+8$  in a neighbourhood of the collision instant.

According to $n_0+n_p+8 \leq 2N+4 \leq 4N -2$ for $N\geq3$, it follows that the set of initial conditions leading to total collisions is locally a finite union of  real smooth submanifold in the neighbourhood of collision instant, and  the
dimensions of the submanifolds depend upon the index of the limiting central configuration and are at most $2N+4$ for the $N$-body problem.

The proof is completed.

$~~~~~~~~~~~~~~~~~~~~~~~~~~~~~~~~~~~~~~~~~~~~~~~~~~~~~~~~~~~~~~~~~~~~~~~~~~~~~~~~~~~~~~~~~~~~~~~~~~~~~~~~~~~~~~~~~~~~~~~~~~~~~~~~~~~~~~~~~~~~~~~~~~~~~~~~~~~~~~~~~~~\Box$

\subsection{Analytic Extension}
\indent\par

Finally, let us examine the question of whether orbits can be extended through total
collision from the viewpoint of Sundman and Siegel, that is, whether a single
solution can be extended as an analytic function of time.
The problem has been studied by Sundman and Siegel for the three-body problem.

It is easy to show that the nature of the singularity corresponding to a total
collision depends on the arithmetical
nature of the eigenvalues $\kappa^{\frac{1}{2}}, \tilde{{\mu}}_5, \cdots, \tilde{{\mu}}_{4+n_p}$.
We here discuss only the case corresponding to  the rhombic central configurations for the four-body problem as a demonstration.

Recall  the  facts (\ref{Rhombiceigenvalues1}) and (\ref{Rhombiceigenvalues}). For simplicity, we consider only the case $\sqrt{3}<\zeta<\zeta_1 \approx 1.7889580612081344$, then
\begin{equation}
0< \kappa^{\frac{1}{2}}<\mu_5<\mu_6 <\mu_8<\mu_7.\nonumber
\end{equation}

It is clear that  $\kappa^{\frac{1}{2}},\mu_5, \mu_6, \mu_7, \mu_8$ are nonresonant for almost all $\zeta \in(\sqrt{3},\zeta_1)$.
According to (\ref{polinomal}),  the coordinate functions $r,\theta,\gamma,q$ of a collision orbit are regular analytic functions of the following variables
\begin{equation}
\tilde{u}_0 =  c_0 e^{\kappa^{\frac{1}{2}} \tau},\tilde{u}_5 = c_5e^{\tilde{{\mu}}_5 \tau},\tilde{u}_6 = c_6e^{\tilde{{\mu}}_6 \tau},\tilde{u}_7 = c_7e^{\tilde{{\mu}}_7 \tau},\tilde{u}_8 = c_8e^{\tilde{{\mu}}_8 \tau}.\nonumber
\end{equation}

Furthermore, it follows from (\ref{timetransformation}) (\ref{Manifoldequation1}) and (\ref{Manifoldequation2}) that  $\gamma,q,r,\theta$  are regular analytic functions of the following variables
\begin{equation}
\check{u}_0 =  c_0 t^{\frac{2}{3} },\check{u}_5 =  c_5 t^{{(2\tilde{{\mu}}_5})/{(3\kappa^{\frac{1}{2}}} )},\check{u}_6 =  c_6 t^{{(2\tilde{{\mu}}_6})/{(3\kappa^{\frac{1}{2}}} )},\check{u}_7 =  c_7 t^{{(2\tilde{{\mu}}_7})/{(3\kappa^{\frac{1}{2}}} )},\check{u}_8 =  c_8 t^{{(2\tilde{{\mu}}_8})/{(3\kappa^{\frac{1}{2}}} )}.\nonumber
\end{equation}

Since all the $\kappa^{\frac{1}{2}},\mu_5, \mu_6, \mu_7, \mu_8$ and their ratios are irrational for almost all $\zeta \in(\sqrt{3},\zeta_1)$,
thus generically we  have an  essential singularity at collision instant $t=0$.  In this case it is not possible to
continue the solutions analytically beyond the  collision.

\begin{remark}
Although one can conclude that  the eigenvalues $\mu_j$ ($j=5,\cdots,2N$) depend continuously upon the value
of the masses $m_k$ ($k=1,\cdots,N$), it is obvious that one can not simply claim that there are values of the masses giving rise to some eigenvalue $\mu_j$
with irrational values, since we can not simply exclude the case that $\mu_j$ are constant value functions about $m_k$.
\end{remark}

\section{Conclusion and Questions}
\indent\par

For the planar $N$-body problem, based on a novel moving frame, which allows us to describe the motion of collision orbit effectively,  we discussed $PISPW$: whether the normalized configuration of the particles must approach a certain central configuration without undergoing infinite spin for a collision orbit.  In the cases corresponding to  central configurations with degree of degeneracy less than or equal to one, we completely solve the problem. We also give a criterion of the problem for the case corresponding to central configurations with degree of degeneracy two;  we further give an answer to the problem  in the case corresponding to all  known degenerate central configurations of four-body. Therefore, for almost every choice
of the masses of the four-body problem,   $PISPW$ is solved. For all the solved cases, we conclude that the normalized configuration of the particles must approach a certain central configuration without undergoing infinite spin for a collision orbit. Finally, we give a measure of the set of initial conditions leading to total collisions: the set has zero measure in the neighbourhood of collision instant.

This work  indicates the fact that $PISPW$ is the link of many interesting subjects in dynamical system. It is our hope that this work may spark similar  interest to the problems on degenerate central configurations and/or degenerate equilibrium points etc.

Many  questions remain to be answered. For example, some concrete questions are the following:
\begin{description}
  \item[i.]  In the cases corresponding to central configurations with three or higher degrees of degeneracy,  how  should one  study  $PISPW$?
  \item[ii.]  For the spatial $N$-body problem,  how  should one  study  $PISPW$?
  \item[iii.] If a solution of a strongly degenerate system (there are many zeros in the eigenvalues of linear part of the system at an equilibrium point $O$) approaches $O$, how should one estimate the rate of tending to the  equilibrium point $O$?
  \item[iv.]  Is it true that all the four-body central configurations except equilateral central configurations have zero or one degree of degeneracy?
\end{description}

We hope to explore some of these questions in future work.

By the way, this work hints the fact that the moving frame and the concomitant coordinate system will be useful when we study   the Newtonian $N$-body problem near relative equilibrium solutions. Inspired by this, we  will investigate the stability  of  relative equilibrium solutions in future work by  making use of the moving frame. Indeed, it is shown that
the theory of KAM is successfully applied to study  the  $N$-body problem near relative equilibrium solutions.
\newpage

\begin{appendices}
\begin{center}
\textbf{\Large Appendix}
\end{center}

\section{Degeneracy of Central Configurations}\label{DegeneracyC.C.}

\subsection{Degeneracy of  a Constrained Critical Point}\label{DegeneracyConstrainedCriticalPoint}
      \indent\par
Let $f$ and $g_k$ be  smooth functions  defined on $\Omega$, $1\leq k\leq m <n$, where $\Omega$ is an open set in $\mathbb{R}^n$. Assume  $\mathcal{M}$ is a
 smooth submanifold  given as the common zero-set of $m$ smooth functions $g_k$.

The classical approach to find the critical points of $\tilde{f}:=f|_{\mathcal{M}}$ ,  which is the
 restriction of  $f$ on $\mathcal{M}$,
  involves the well-known
\emph{method of Lagrange multipliers}. The method avoids looking for local coordinates on the submanifold $\mathcal{M}$, but consists of introducing $m$  ``undetermined multipliers" $\lambda_k$, defining the following auxiliary function
\begin{equation*}
    L=f+ \sum_{k=1}^m\lambda_kg_k
\end{equation*}
and finding its critical points; i.e., by regarding $L$ as a function of
\begin{center}
$(x,\bar{\lambda})=(x_1,\cdots,x_n,\lambda_1,\cdots,\lambda_m)\in \Omega\times \mathbb{R}^m$,
\end{center}
 the critical
points are obtained by solving the equations
\begin{equation*}
    \nabla L=0.
\end{equation*}

Similarly, one can investigate the degeneracy of a critical
point  in a straightforward manner: assume $\underline{a}$ is  a critical point  of $\tilde{f}$, let $(\underline{a},\underline{\lambda})$ be  an auxiliary critical point  of $L$; by regarding $L$ as a function of $(x,\bar{\lambda})$, the degeneracy of   $\tilde{f}$ at $\underline{a}$ is the  same as that of $L$ at $(\underline{a},\underline{\lambda})$, i.e., we have the following result (please refer to \cite{zbMATH00579429} for more detail)
\begin{theorem}\label{DegeneracyConstrainedCriticalPoint}\emph{(\cite{zbMATH00579429})}
 The nullity of the Hessian $D^2\tilde{f}(\underline{a})$ equals the nullity of the Hessian $D^2L(\underline{a},\underline{\lambda})$.
\end{theorem}
Note that  the Hessian $D^2\tilde{f}(\underline{a})$ is the $n-m$ by $n-m$
 symmetric matrix, and the Hessian $D^2L(\underline{a},\underline{\lambda})$ which is called  the bordered Hessian is the $n+m$ by $n+m$
 symmetric matrix.

\subsection{Degeneracy of Central Configurations}
\indent\par

When investigating the degeneracy of central configurations, the method in \ref{DegeneracyConstrainedCriticalPoint} is certainly  useful. For example, $\mathcal{X}$ is the submanifold of $(\mathbb{R}^2)^N $ such that the center of mass  is at the origin; then the nullity of the Hessian $D^2\widetilde{U}$ in $\mathcal{X}$ is naturally reduced to that in $(\mathbb{R}^2)^N\times \mathbb{R}^2$. The benefits of this is that  practical calculations do not  look for local coordinates on $\mathcal{X}$.
As it happens often, it is difficult to find appropriate local coordinates for concrete problems.

In the following, we use the cartesian coordinates of $(\mathbb{R}^2)^N$ to give a  slightly better criterion  than the bordered Hessian in $(\mathbb{R}^2)^N\times \mathbb{R}^2$ for investigating the degeneracy of central configurations.

Let
\begin{displaymath}
\{\mathbf{e}_1, \cdots,\mathbf{e}_{2N}\}
\end{displaymath}
 be the standard basis of $(\mathbb{R}^2)^N $, where $\mathbf{e}_j\in (\mathbb{R}^2)^N $ has
unity at the $j$-th component and zero at all others. Then every N-body configuration $\mathbf{r}\in (\mathbb{R}^2)^N $ can be written as
\begin{displaymath}
\mathbf{r} =  \sum_{j = 1}^{2N} {x^j \mathbf{e}_j},
\end{displaymath}
and
\begin{displaymath}
 (x^1, x^2, \cdots, x^{2N})^\top
\end{displaymath}
are the coordinates of $\mathbf{r}$ in the standard basis,  here ``$\top$" denotes transposition of matrix. It is also true that $\mathbf{r}_j = (x^{2j-1},x^{2j})^\top$ for $j = 1, 2, \cdots, N$.
Then \begin{displaymath}
\mathcal{X} = \{ \mathbf{r} \in (\mathbb{R}^2)^N| \langle\mathbf{r},\mathcal{E}_1\rangle = 0, \langle\mathbf{r},\mathcal{E}_2\rangle = 0   \},
\end{displaymath}
where
\begin{displaymath}
\mathcal{E}_1 = \sum_{j = 1}^{N} {\mathbf{e}_{2j-1}}=(1,0,\cdots,1,0)^\top, \mathcal{E}_2 = \sum_{j = 1}^{N} {\mathbf{e}_{2j}}=(0,1,\cdots,0,1)^\top.
\end{displaymath}

Let $\mathfrak{M}$ be the matrix
\begin{displaymath}
diag(m_1,m_1,m_2,m_2, \cdots, m_N, m_N),
\end{displaymath} where ``diag" means diagonalmatrix.
Then
\begin{displaymath}
\langle\mathbf{r},\mathbf{r}\rangle =   \mathbf{r}^\top \mathfrak{M} \mathbf{r}.
\end{displaymath}

Given a central configuration $\mathbf{r}_0 = (\mathbf{r}_1,\cdots, \mathbf{r}_N) \in \mathcal{X}\backslash \Delta$,
a straightforward computation shows that the Hessian $D^2\widetilde{\mathcal{U}}(\mathbf{r}_0)$ in $(\mathbb{R}^2)^N $ is
\begin{displaymath}
I^{\frac{1}{2}} (\lambda \mathbb{I} + \mathfrak{M}^{-1} \mathfrak{B}) - 3 I^{-\frac{1}{2}} \lambda \mathbf{r}_0\mathbf{r}_0^\top \mathfrak{M},
\end{displaymath}
where $\mathfrak{B}$ is the Hessian of $\mathcal{U}$ evaluated at $\mathbf{r}_0$ and can be viewed as an $N\times N$ array of $2 \times 2$ blocks:
\begin{center}
$\mathfrak{B} = \left(
       \begin{array}{ccc}
         B_{11} & \cdots & B_{1N} \\
         \vdots & \ddots & \vdots \\
         B_{N1} & \cdots & B_{NN}\\
       \end{array}
     \right)
$
\end{center}
The off-diagonal blocks are given by:
\begin{center}
$B_{jk} = \frac{m_j m_k}{r^3_{jk}}[\mathbb{I}-\frac{3(\mathbf{r}_k - \mathbf{r}_j)(\mathbf{r}_k - \mathbf{r}_j)^\top}{r^2_{jk}}],
$
\end{center}
where $r_{jk}=|\mathbf{r}_k - \mathbf{r}_j|$, and $\mathbb{I}$ is the identity matrix of order 2. However, as a matter of notational convenience, the \textbf{identity matrix} of any order will always be denoted by $\mathbb{I}$, and the order of $\mathbb{I}$ can be determined according to the context.
The diagonal blocks are given by:
\begin{displaymath}
B_{kk} = -\sum_{1\leq j\leq N, j\neq k} B_{jk}.
\end{displaymath}

Thus $\{\mathcal{E}_1, \mathcal{E}_2,\mathcal{E}_3, \mathcal{E}_4, \mathcal{E}_5,\cdots, \mathcal{E}_{2N}\}$ is just an orthogonal basis of $(\mathbb{R}^2)^N $ such that
\begin{equation*}
    D^2\widetilde{\mathcal{U}}(\mathbf{r}_0)\mathcal{E}_j=\mu_j\mathcal{E}_j,~~~~~~~~~~~~j=1,\cdots,2N,
\end{equation*}
and
\begin{equation*}
    \begin{array}{c}
    \mu_1=\mu_2=I^{\frac{1}{2}}\lambda=I^{-\frac{1}{2}} U, ~~~~~~~~~~\mu_3=\mu_4=0,\\
       \mathcal{E}_3=\mathbf{r}_0, ~~~~~~~~~~\mathcal{E}_4=\mathbf{i}\mathbf{r}_0, \\
       \mathcal{P}_{\mathbf{r}_0}=span\{{\mathcal{E}}_3, {\mathcal{E}}_4\}, \mathcal{P}^{\bot}_{\mathbf{r}_0}=span\{{\mathcal{E}}_5, \cdots, {\mathcal{E}}_{2N}\}.
     \end{array}
\end{equation*}

Then
\begin{displaymath}
\{\widehat{\mathcal{E}}_1, \widehat{\mathcal{E}}_2, \widehat{\mathcal{E}}_3, \widehat{\mathcal{E}}_4, \cdots, \widehat{\mathcal{E}}_{2N}\}
\end{displaymath}
consisting of eigenvectors of $\mathfrak{D}$ is an orthonormal basis of the space $(\mathbb{R}^2)^N $ with respect to the scalar product $\langle,\rangle$, that is,
\begin{displaymath}
(\widehat{\mathcal{E}}_1, \widehat{\mathcal{E}}_2, \widehat{\mathcal{E}}_3, \widehat{\mathcal{E}}_4, \cdots, \widehat{\mathcal{E}}_{2N})^\top \mathfrak{M} (\widehat{\mathcal{E}}_1, \widehat{\mathcal{E}}_2, \widehat{\mathcal{E}}_3, \widehat{\mathcal{E}}_4, \cdots, \widehat{\mathcal{E}}_{2N}) = \mathbb{I}.
\end{displaymath}
By
\begin{displaymath}
D^2\widetilde{\mathcal{U}}(\mathbf{r}_0) (\widehat{\mathcal{E}}_1,  \cdots, \widehat{\mathcal{E}}_{2N}) =diag(\mu_1, \mu_2, \cdots, \mu_{2N})(\widehat{\mathcal{E}}_1,  \cdots, \widehat{\mathcal{E}}_{2N}),
\end{displaymath}
it follows that
\begin{equation}\label{Hessian0}
(\widehat{\mathcal{E}}_1,  \cdots, \widehat{\mathcal{E}}_{2N})^{\top} (I^{\frac{1}{2}} (\lambda \mathfrak{M} + \mathfrak{B}) - 3 I^{-\frac{1}{2}} \lambda \mathfrak{M} \mathcal{E}_3\mathcal{E}_3^\top \mathfrak{M}) (\widehat{\mathcal{E}}_1,  \cdots, \widehat{\mathcal{E}}_{2N}) = diag(\mu_1, \cdots, \mu_{2N}).
\end{equation}

Thanks to
\begin{displaymath}(\widehat{\mathcal{E}}_1, \cdots, \widehat{\mathcal{E}}_{2N})^{\top} ( \mathfrak{M} \mathcal{E}_1\mathcal{E}_1^\top \mathfrak{M}) (\widehat{\mathcal{E}}_1,  \cdots, \widehat{\mathcal{E}}_{2N}) = diag(\mathfrak{m}^2, 0, 0, 0, \cdots, 0),\end{displaymath}
\begin{displaymath}(\widehat{\mathcal{E}}_1, \cdots, \widehat{\mathcal{E}}_{2N})^{\top} ( \mathfrak{M} \mathcal{E}_2\mathcal{E}_2^\top \mathfrak{M}) (\widehat{\mathcal{E}}_1,  \cdots, \widehat{\mathcal{E}}_{2N}) = diag(0,\mathfrak{m}^2, 0, 0,  \cdots, 0),\end{displaymath}
\begin{displaymath}(\widehat{\mathcal{E}}_1, \cdots, \widehat{\mathcal{E}}_{2N})^{\top} (  \mathfrak{M} \mathcal{E}_3\mathcal{E}_3^\top \mathfrak{M}) (\widehat{\mathcal{E}}_1,  \cdots, \widehat{\mathcal{E}}_{2N}) = diag(0, 0, I, 0, \cdots, 0),\end{displaymath}
and \begin{displaymath}(\widehat{\mathcal{E}}_1, \cdots, \widehat{\mathcal{E}}_{2N})^{\top} (  \mathfrak{M} \mathcal{E}_4\mathcal{E}_4^\top \mathfrak{M}) (\widehat{\mathcal{E}}_1,  \cdots, \widehat{\mathcal{E}}_{2N}) = diag(0, 0,0,  I, \cdots, 0),\end{displaymath}
 we have
\begin{equation}\label{Hessian}
(\widehat{\mathcal{E}}_1,  \cdots, \widehat{\mathcal{E}}_{2N})^{\top}  (\lambda \mathfrak{M} + \mathfrak{B})  (\widehat{\mathcal{E}}_1,  \cdots, \widehat{\mathcal{E}}_{2N})
= diag( \lambda, \lambda,3  \lambda, 0,  \frac{\mu_5}{\|\mathcal{E}_{3}\|}, \cdots,  \frac{\mu_{2N}}{\|\mathcal{E}_{3}\|}),
\end{equation}
\begin{equation}\label{Hessiandeg1}
(\widehat{\mathcal{E}}_1,  \cdots, \widehat{\mathcal{E}}_{2N})^{\top}  (\lambda \mathfrak{M} + \mathfrak{B}+\mathfrak{M} \mathcal{E}_4\mathcal{E}_4^\top \mathfrak{M})  (\widehat{\mathcal{E}}_1,  \cdots, \widehat{\mathcal{E}}_{2N})
= diag( \lambda, \lambda,3  \lambda, I,  \frac{\mu_5}{\|\mathcal{E}_{3}\|}, \cdots,  \frac{\mu_{2N}}{\|\mathcal{E}_{3}\|}),
\end{equation}
\begin{equation}\label{Hessiandeg2}
\begin{array}{c}
  (\widehat{\mathcal{E}}_1,  \cdots, \widehat{\mathcal{E}}_{2N})^{\top}  (\lambda \mathfrak{M} + \mathfrak{B}-\frac{\lambda}{\mathfrak{m}^2}\mathfrak{M} \mathcal{E}_1\mathcal{E}_1^\top \mathfrak{M}-\frac{\lambda}{\mathfrak{m}^2}\mathfrak{M} \mathcal{E}_2\mathcal{E}_2^\top \mathfrak{M}-\frac{3\lambda}{I}\mathfrak{M} \mathcal{E}_3\mathcal{E}_3^\top \mathfrak{M})  (\widehat{\mathcal{E}}_1,  \cdots, \widehat{\mathcal{E}}_{2N}) \\
  = diag(0, 0,0, 0,  \frac{\mu_5}{\|\mathcal{E}_{3}\|}, \cdots,  \frac{\mu_{2N}}{\|\mathcal{E}_{3}\|}).
\end{array}
\end{equation}

These results will be useful in investigating the degeneracy of central configurations. For example, the nullity of $\lambda \mathfrak{M} + \mathfrak{B}+\mathfrak{M} \mathcal{E}_4\mathcal{E}_4^\top \mathfrak{M}$ in (\ref{Hessiandeg1}) is just the degree of degeneracy of $\mathbf{r}_0$; when investigating   central configurations of four-body with degree of degeneracy two,  the problem reduces to investigate that whether the matrix
\begin{equation*}
   \lambda \mathfrak{M} + \mathfrak{B}-\frac{\lambda}{\mathfrak{m}^2}\mathfrak{M} \mathcal{E}_1\mathcal{E}_1^\top \mathfrak{M}-\frac{\lambda}{\mathfrak{m}^2}\mathfrak{M} \mathcal{E}_2\mathcal{E}_2^\top \mathfrak{M}-\frac{3\lambda}{I}\mathfrak{M} \mathcal{E}_3\mathcal{E}_3^\top \mathfrak{M}
\end{equation*}
is a positive semi-definite matrix with rank $2$ by (\ref{Hessiandeg2}).

Of course, if one finds a basis of $\mathcal{P}^{\bot}_{\mathbf{r}_0}$, say $\mathcal{F}_1,\mathcal{F}_2,\cdots,\mathcal{F}_{2N-4}$, then it is simpler to consider the nullity of
\begin{equation*}
    ({\mathcal{F}}_1,  \cdots, {\mathcal{F}}_{2N-4})^{\top}  (\lambda \mathfrak{M} + \mathfrak{B})   ({\mathcal{F}}_1,  \cdots, {\mathcal{F}}_{2N-4}).
\end{equation*}
Based on these considerations, we will investigate the degeneracy of central configurations  with an axis of symmetry in the  four-body problem  in next section.

      \section{ Central Configurations of Four-body }\label{C.C.ofFour-body}
      \indent\par

Viewing as a preliminary study on the problem of degeneracy,   let us investigate the  degenerate central configurations of the planar four-body problem  with an axis of symmetry in this section. They consist of
systems of point particles in $\mathbb{R}^2$ whose configurations have the following geometric
properties:

1) There are 2 particles $m_3,m_4$ lying on a fixed line which is the axis of symmetry of problem, and the axis of symmetry is assumed as $y$-axis.

2) Another 2 particles $m_1,m_2$ are symmetric with regard to the $y$-axis.

This kind of configurations are usually called \textbf{kite configurations}. It is easy to see that $m_1=m_2$ is necessary  for a kite central configuration. Geometry of kite configurations may be seen in the following Figure \ref{equilateral central configuration}.

Without loss of generality, suppose
\begin{displaymath}
\begin{array}{c}
  m_1=m_2=1, \\
  \mathbf{r}_1 = (-s,-t)^\top,\\
  \mathbf{r}_2 = (s,-t)^\top, \\
  \mathbf{r}_3 = (0,u)^\top, \\
\mathbf{r}_4 = (0,u-1)^\top,
\end{array}
\end{displaymath}
where $s>0$ and $u= \frac{2t + m_4}{m_3 + m_4}$.

Then the equations (\ref{centralconfiguration}) of central configurations become:
\begin{equation}\label{centralconfiguration2}
\left\{
             \begin{array}{lr}
              \frac{1}{4s^3}+\frac{m_3}{[s^2+(u+t)^2]^\frac{3}{2}}+\frac{m_4}{[s^2+(u+t-1)^2]^\frac{3}{2}}=\lambda  &  \\
             \frac{m_3(u+t)}{[s^2+(u+t)^2]^\frac{3}{2}}+\frac{m_4(u+t-1)}{[s^2+(u+t-1)^2]^\frac{3}{2}}=\lambda t & \\
\frac{2(u+t)}{[s^2+(u+t)^2]^\frac{3}{2}}+ m_4=\lambda u & \\
\frac{2(u+t-1)}{[s^2+(u+t-1)^2]^\frac{3}{2}}-m_3=\lambda (u-1). &
             \end{array}
\right.
\end{equation}

Set \begin{displaymath}
\begin{array}{c}
  \mathcal{V}=(-1,0,1,0,0,0,0,0)^\top, \\
  \mathcal{P}=(0,-1,0,-1,0,\frac{2}{m_3 + m_4},0,\frac{2}{m_3 + m_4})^\top.
\end{array}
\end{displaymath}
Then
\begin{displaymath}
\begin{array}{c}
  \mathbf{i}\mathcal{V}=(0,-1,0,1,0,0,0,0)^\top, \\
  \mathbf{i}\mathcal{P}=(1,0,1,0,-\frac{2}{m_3 + m_4},0,-\frac{2}{m_3 + m_4},0)^\top.
\end{array}
\end{displaymath}
Note that
\begin{displaymath}
\begin{array}{c}
{\mathcal{E}}_3= \mathbf{r}=(-s,-t,s,-t,0,u,0,u-1)^\top, \\
  {\mathcal{E}}_4=\mathbf{i} \mathbf{r}=(t,-s,t,s,-u,0,1-u,0)^\top,
\end{array}
\end{displaymath}
\begin{displaymath}
\begin{array}{c}
   \mathfrak{B}  {\mathcal{E}}_3 = 2 \lambda \mathfrak{M} {\mathcal{E}}_3,\\
  \mathfrak{B}  {\mathcal{E}}_4 = -\lambda \mathfrak{M} {\mathcal{E}}_4,
\end{array}
\end{displaymath}
and
\begin{displaymath}
\begin{array}{c}
\langle\mathbf{i}\mathcal{V},{\mathcal{E}}_3\rangle= \langle\mathbf{i}\mathcal{P},{\mathcal{E}}_3\rangle=0, \\
\langle \mathcal{V}, {\mathcal{E}}_4\rangle= \langle \mathcal{P}, {\mathcal{E}}_4\rangle=0,\\
\langle \mathcal{V}, \mathcal{P}\rangle=\langle \mathbf{i}\mathcal{V},\mathbf{i}\mathcal{P}\rangle=0,\\
\langle \mathcal{V},\mathbf{i}\mathcal{V}\rangle=\langle \mathcal{V},\mathbf{i}\mathcal{P}\rangle=0,\\
\langle \mathcal{P},\mathbf{i}\mathcal{V}\rangle=\langle \mathcal{P},\mathbf{i}\mathcal{P}\rangle=0,
\end{array}
\end{displaymath}
\begin{displaymath}
\begin{array}{c}
I= \langle{\mathcal{E}}_3,{\mathcal{E}}_3\rangle= \langle{\mathcal{E}}_4,{\mathcal{E}}_4\rangle=2 (s^2+t^2)+\frac{4t^2+ m_3 m_4 }{ m_3 + m_4}, \\
\langle \mathcal{V}, \mathcal{V}\rangle= \langle \mathbf{i}\mathcal{V}, \mathbf{i}\mathcal{V}\rangle=2,\\
\langle \mathcal{P}, \mathcal{P}\rangle=\langle \mathbf{i}\mathcal{P},\mathbf{i}\mathcal{P}\rangle=2+\frac{4}{m_3 + m_4}.
\end{array}
\end{displaymath}
Then a straightforward computation shows that:
\begin{displaymath}
\begin{array}{c}
 \mathcal{V}^{\top}  \mathfrak{B}  \mathcal{V}  =\frac{2m_3}{[s^2+(u+t)^2]^\frac{3}{2}}(\frac{3s^2}{s^2+(u+t)^2}-1)+\frac{2m_4}{[s^2+(u+t-1)^2]^\frac{3}{2}}(\frac{3s^2}{s^2+(u+t-1)^2}-1)+\frac{1}{s^3}\\
\mathcal{V}^{\top}  \mathfrak{B}  \mathcal{P} =\mathcal{P}^{\top}  \mathfrak{B}  \mathcal{V} =\frac{6m_3 s(u+t)}{[s^2+(u+t)^2]^\frac{5}{2}}(\frac{2}{m_3 + m_4}+1)+\frac{6m_4 s(u+t-1)}{[s^2+(u+t-1)^2]^\frac{5}{2}}(\frac{2}{m_3 + m_4}+1)\\
\mathcal{P}^{\top}  \mathfrak{B}  \mathcal{P} =[\frac{2m_3}{[s^2+(u+t)^2]^\frac{3}{2}}(\frac{3(u+t)^2}{s^2+(u+t)^2}-1)+\frac{2m_4}{[s^2+(u+t-1)^2]^\frac{3}{2}}(\frac{3(u+t-1)^2}{s^2+(u+t-1)^2}-1)](\frac{2}{m_3 + m_4}+1)^2
\end{array}
\end{displaymath}

\begin{displaymath}
\begin{array}{c}
 (\mathbf{i}\mathcal{V})^{\top}  \mathfrak{B}  (\mathbf{i}\mathcal{V})= \frac{2m_3}{[s^2+(u+t)^2]^\frac{3}{2}}(\frac{3(u+t)^2}{s^2+(u+t)^2}-1)+\frac{2m_4}{[s^2+(u+t-1)^2]^\frac{3}{2}}(\frac{3(u+t-1)^2}{s^2+(u+t-1)^2}-1)-\frac{1}{2s^3}\\
(\mathbf{i}\mathcal{V})^{\top}  \mathfrak{B}  (\mathbf{i}\mathcal{P})= (\mathbf{i}\mathcal{P})^{\top}  \mathfrak{B}  (\mathbf{i}\mathcal{V})= \frac{-6m_3 s(u+t)}{[s^2+(u+t)^2]^\frac{5}{2}}(\frac{2}{m_3 + m_4}+1)+\frac{-6m_4 s(u+t-1)}{[s^2+(u+t-1)^2]^\frac{5}{2}}(\frac{2}{m_3 + m_4}+1)\\
(\mathbf{i}\mathcal{P})^{\top}  \mathfrak{B}  (\mathbf{i}\mathcal{P})= [\frac{2m_3}{[s^2+(u+t)^2]^\frac{3}{2}}(\frac{3s^2}{s^2+(u+t)^2}-1)+\frac{2m_4}{[s^2+(u+t-1)^2]^\frac{3}{2}}(\frac{3s^2}{s^2+(u+t-1)^2}-1)](\frac{2}{m_3 + m_4}+1)^2
\end{array}
\end{displaymath}

\begin{displaymath}
\begin{array}{c}
 \mathcal{V}^{\top}  \mathfrak{B}  (\mathbf{i}\mathcal{V})= 0\\
\mathcal{V}^{\top}  \mathfrak{B}  \mathbf{i}\mathcal{P}=0\\
\mathcal{P}^{\top}  \mathfrak{B}  (\mathbf{i}\mathcal{V})= 0\\
\mathcal{P}^{\top}  \mathfrak{B}  (\mathbf{i}\mathcal{P})= 0.
\end{array}
\end{displaymath}

Set \begin{displaymath}
\begin{array}{c}
  \tilde{\mathcal{V}}= \mathcal{V} - \frac{\langle\mathcal{V},{\mathcal{E}}_3\rangle}{\langle{\mathcal{E}}_3,{\mathcal{E}}_3\rangle}{\mathcal{E}}_3= \mathcal{V} - \frac{2s}{2 (s^2+t^2)+\frac{4t^2+ m_3 m_4 }{ m_3 + m_4}}{\mathcal{E}}_3,\\
 \tilde{ \mathcal{P}}=\mathcal{P} - \frac{\langle\mathcal{P},{\mathcal{E}}_3\rangle}{\langle{\mathcal{E}}_3,{\mathcal{E}}_3\rangle}{\mathcal{E}}_3= \mathcal{P} - \frac{2t+2u-\frac{2m_4}{m_3+ m_4}}{2 (s^2+t^2)+\frac{4t^2+ m_3 m_4 }{ m_3 + m_4}}{\mathcal{E}}_3.
\end{array}
\end{displaymath}

Then, by (\ref{Hessian}),  the central configuration $\mathcal{E}_3$ is {degenerate} if and only if the $4\times 4$ matrix \\
$(\tilde{\mathcal{V}},   \tilde{ \mathcal{P}},\mathbf{i}\tilde{\mathcal{V}},   \mathbf{i}\tilde{\mathcal{P}})^{\top} (\lambda \mathfrak{M} + \mathfrak{B})  (\tilde{\mathcal{V}},   \tilde{ \mathcal{P}},\mathbf{i}\tilde{\mathcal{V}},   \mathbf{i}\tilde{\mathcal{P}})$ is degenerate.

Due to
\begin{displaymath}
\begin{array}{c}
 (\tilde{\mathcal{V}},   \tilde{ \mathcal{P}},\mathbf{i}\tilde{\mathcal{V}},   \mathbf{i}\tilde{\mathcal{P}})^{\top} (\lambda \mathfrak{M} + \mathfrak{B})  (\tilde{\mathcal{V}},   \tilde{ \mathcal{P}},\mathbf{i}\tilde{\mathcal{V}},   \mathbf{i}\tilde{\mathcal{P}}) =  \\
 \left(
   \begin{array}{cc}
     (\mathcal{V},   \mathcal{P})^{\top} (\lambda \mathfrak{M} + \mathfrak{B})  (\mathcal{V},   \mathcal{P}) &   \\
       & (\mathbf{i}\mathcal{V},    \mathbf{i}\mathcal{P})^{\top} (\lambda \mathfrak{M} + \mathfrak{B})  (\mathbf{i}\mathcal{V},    \mathbf{i}\mathcal{P}) \\
   \end{array}
 \right)
 -\\

 \left(
   \begin{array}{cccc}
    \frac{3 \lambda |\langle\mathcal{V},{\mathcal{E}}_3\rangle|^2}{\langle{\mathcal{E}}_3,{\mathcal{E}}_3\rangle} &  \frac{ 3\lambda \langle\mathcal{V},{\mathcal{E}}_3\rangle\langle\mathcal{P},{\mathcal{E}}_3\rangle}{\langle{\mathcal{E}}_3,{\mathcal{E}}_3\rangle} &   &   \\
      \frac{ 3\lambda \langle\mathcal{V},{\mathcal{E}}_3\rangle\langle\mathcal{P},{\mathcal{E}}_3\rangle}{\langle{\mathcal{E}}_3,{\mathcal{E}}_3\rangle} & \frac{3 \lambda |\langle\mathcal{P},{\mathcal{E}}_3\rangle|^2}{\langle{\mathcal{E}}_3,{\mathcal{E}}_3\rangle} &   &   \\
       &   & 0 & 0 \\
      &   & 0 & 0 \\
   \end{array}
 \right),
\end{array}
\end{displaymath}
it follows that the degeneracy of $\mathcal{E}_3$ attributes to that the $2\times 2$ matrix
\begin{displaymath}
\begin{array}{c}
 \left(
   \begin{array}{cc}
     \mathcal{V}^{\top} (\lambda \mathfrak{M} + \mathfrak{B})  \mathcal{V}- \frac{3 \lambda |\langle\mathcal{V},{\mathcal{E}}_3\rangle|^2}{\langle{\mathcal{E}}_3,{\mathcal{E}}_3\rangle} &  \mathcal{V}^{\top} (\lambda \mathfrak{M} + \mathfrak{B})  \mathcal{P}-\frac{ 3\lambda \langle\mathcal{V},{\mathcal{E}}_3\rangle\langle\mathcal{P},{\mathcal{E}}_3\rangle}{\langle{\mathcal{E}}_3,{\mathcal{E}}_3\rangle} \\
       \mathcal{P}^{\top} (\lambda \mathfrak{M} + \mathfrak{B})  \mathcal{V}-\frac{ 3\lambda \langle\mathcal{V},{\mathcal{E}}_3\rangle\langle\mathcal{P},{\mathcal{E}}_3\rangle}{\langle{\mathcal{E}}_3,{\mathcal{E}}_3\rangle}  &  \mathcal{P}^{\top} (\lambda \mathfrak{M} + \mathfrak{B}) \mathcal{P} - \frac{3 \lambda |\langle\mathcal{P},{\mathcal{E}}_3\rangle|^2}{\langle{\mathcal{E}}_3,{\mathcal{E}}_3\rangle}\\
   \end{array}
 \right)
\end{array}
\end{displaymath}
or $(\mathbf{i}\mathcal{V},   \mathbf{i}\mathcal{P})^{\top} (\lambda \mathfrak{M} + \mathfrak{B})  (\mathbf{i}\mathcal{V},   \mathbf{i}\mathcal{P})$ is degenerate.

By the equations (\ref{centralconfiguration2}), we have
\begin{equation}\label{massbiaoshi}
\left\{
             \begin{array}{lr}
              \left(\frac{1}{[s^2+(u+t)^2]^\frac{3}{2}}-1\right)m_3= 2\left(u+t-1\right)\left(\frac{1}{8s^3}-\frac{1}{[s^2+(u+t-1)^2]^\frac{3}{2}}\right),&\\
              \left(1-\frac{1}{[s^2+(u+t-1)^2]^\frac{3}{2}}\right)m_4= 2\left(u+t\right)\left(\frac{1}{8s^3}-\frac{1}{[s^2+(u+t)^2]^\frac{3}{2}}\right).&
             \end{array}
\right.\nonumber
\end{equation}
Thus,  except  the  cases of $s=\frac{\sqrt{3}}{2}, t=\frac{1}{2}, u=1$ or $s=\frac{\sqrt{3}}{2}, t=-\frac{1}{2}, u=0$, the masses can be presented by geometric elements of the central configuration:
\begin{equation}\label{massbiaoshi1}
\left\{
             \begin{array}{lr}
              m_3= \frac{2\left(u+t-1\right)\left(\frac{1}{8s^3}-\frac{1}{[s^2+(u+t-1)^2]^\frac{3}{2}}\right)}{\frac{1}{[s^2+(u+t)^2]^\frac{3}{2}}-1}&\\
              m_4= \frac{2\left(u+t\right)\left(\frac{1}{8s^3}-\frac{1}{[s^2+(u+t)^2]^\frac{3}{2}}\right)}{1-\frac{1}{[s^2+(u+t-1)^2]^\frac{3}{2}}}&
             \end{array}
\right.
\end{equation}

The central configurations in the  cases of $s=\frac{\sqrt{3}}{2}, t=\frac{1}{2}, u=1$ or $s=\frac{\sqrt{3}}{2}, t=-\frac{1}{2}, u=0$ are formed by the three particles  at the vertices of an equilateral
triangle and a fourth particle at the centroid.  We shall call them the equilateral central configurations.
\subsection{Equilateral Central Configurations}\label{Equilateral Central Configurations}
\indent\par

\begin{figure}
  \centering
  \includegraphics[width=6cm]{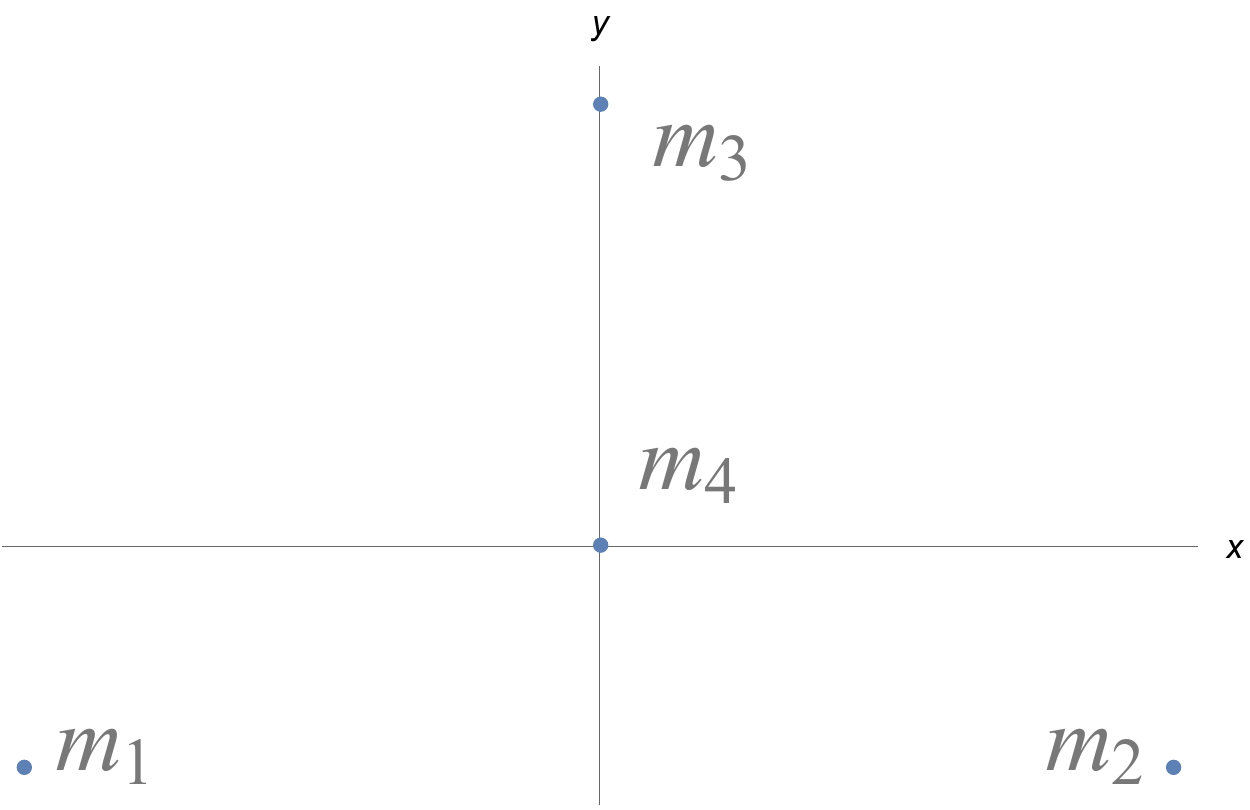}\\
  \caption{equilateral central configuration}\label{equilateral central configuration}
\end{figure}

In this subsection, let us investigate the degeneracy of equilateral central configurations, i.e., $s=\frac{\sqrt{3}}{2}, t=\frac{1}{2}, u=1$ or $s=\frac{\sqrt{3}}{2},
t=-\frac{1}{2}, u=0$. Without loss of generality, we only investigate the case of $s=\frac{\sqrt{3}}{2}, t=\frac{1}{2}, u=1$ as illustrated in Figure \ref{equilateral central configuration}.

Then $m_3=1$, and the equations (\ref{centralconfiguration2}) of central configurations become:
\begin{equation}\label{centralconfiguration4}
\lambda=\frac{1}{\sqrt{3}}+m_4.\nonumber
\end{equation}

It is a classical result that $m_4=(64\sqrt{3}+81)/249$ is the unique value of the mass parameter $m_4$ corresponding to central configurations
with  degree of degeneracy two by Palmore \cite{palmore1975classifying,Palmore1976Measure}.  We shall  reproduce this result and further pursue $\{ \mathcal{E}_5, \mathcal{E}_6,\mathcal{E}_7, \mathcal{E}_{8}\}$ in the following.

A straightforward computation shows that the $2\times 2$ matrices $(\tilde{\mathcal{V}},   \tilde{ \mathcal{P}})^{\top} (\lambda \mathfrak{M} + \mathfrak{B})  (\tilde{\mathcal{V}},   \tilde{ \mathcal{P}})$ and $ (\mathbf{i}\tilde{\mathcal{V}},   \mathbf{i}\tilde{\mathcal{P}})^{\top} (\lambda \mathfrak{M} + \mathfrak{B})  (\mathbf{i}\tilde{\mathcal{V}},   \mathbf{i}\tilde{\mathcal{P}})$ become respectively
\begin{displaymath}
\left(
   \begin{array}{cc}
     \frac{1}{2} \left(3 m_4+\sqrt{3}\right) &  \frac{(\sqrt{3}m_4-1)(3+m_4)}{2+2m_4}       \\
     \frac{(\sqrt{3}m_4-1)(3+m_4)}{2+2m_4} & \frac{\left(m_4+3\right) \left(9 m_4^2+\left(11 \sqrt{3}-45\right) m_4+9 \sqrt{3}\right)}{18 \left(m_4+1\right){}^2}       \\

   \end{array}
 \right)
\end{displaymath}
and
\begin{displaymath}
\left(
   \begin{array}{cc}
    \frac{1}{2} \left(3 m_4+\sqrt{3}\right) & -\frac{\left(m_4+3\right) \left(3 \sqrt{3} m_4+1\right)}{2 \left(m_4+1\right)} \\
    -\frac{\left(m_4+3\right) \left(3 \sqrt{3} m_4+1\right)}{2 \left(m_4+1\right)} & \frac{(m_4+3) \left(81 m_4^2+\left(11 \sqrt{3}+171\right) m_4+9 \sqrt{3}\right)}{18 (m_4+1)^2} \\
   \end{array}
 \right)
\end{displaymath}
Since both of the determinants of the above two matrices are
\begin{displaymath}
\frac{m_4 \left(m_4+3\right) \left(\left(5 \sqrt{3}-18\right) m_4+3 \sqrt{3}+2\right)}{3 \left(m_4+1\right){}^2},
\end{displaymath}
the central configuration $\mathcal{E}_3$ is {degenerate} if and only if
\begin{equation}\label{degeneracy of rhombic central configurations2}
(5 \sqrt{3}-18) m_4+3 \sqrt{3}+2=0.\nonumber
\end{equation}
Obviously, $m_4=\frac{3 \sqrt{3}+2}{18-5 \sqrt{3}}=\frac{81+64\sqrt{3}}{249}$ is the unique solution of above equation.

Furthermore, it is easy to see that the $4\times 4$ matrix
$(\tilde{\mathcal{V}},   \tilde{ \mathcal{P}},\mathbf{i}\tilde{\mathcal{V}},   \mathbf{i}\tilde{\mathcal{P}})^{\top} (\lambda \mathfrak{M} + \mathfrak{B})  (\tilde{\mathcal{V}},   \tilde{ \mathcal{P}},\mathbf{i}\tilde{\mathcal{V}},   \mathbf{i}\tilde{\mathcal{P}})$ is positive definite for $0<m_4<\frac{81+64\sqrt{3}}{249}$ and indefinite for $m_4>\frac{81+64\sqrt{3}}{249}$. So the equilateral central configuration is a local
minimum of the function $I^{\frac{1}{2}}\mathcal{U}$ for $0<m_4<\frac{81+64\sqrt{3}}{249}$, and the equilateral central configuration is a saddle point of the function $I^{\frac{1}{2}}\mathcal{U}$ for $m_4>\frac{81+64\sqrt{3}}{249}$.

As a result of (\ref{Hessian}), $\{ \mathcal{E}_5, \mathcal{E}_6,\mathcal{E}_7, \mathcal{E}_{8}\}$ can be obtained by calculating eigenvectors of the matrix $\lambda \mathbb{I} + \mathfrak{M}^{-1} \mathfrak{B}$.

A straightforward computation shows that, for the case $m_4=\frac{81+64\sqrt{3}}{249}$,
\begin{equation}\label{eigenvectors}
\begin{array}{c}
  \mathcal{E}_5=\left(\frac{64 \sqrt{3}+81}{498} ,-\frac{741 \sqrt{3}+908}{1494},\frac{64 \sqrt{3}+81}{498} ,\frac{741 \sqrt{3}+908}{1494},0,0,-1,0\right)^\top \\
  \mathcal{E}_6 = \left(\frac{165 \sqrt{3}+179}{747} ,-\frac{371 \sqrt{3}+738}{2241},-\frac{165 \sqrt{3}+179}{747} ,-\frac{371 \sqrt{3}+738}{2241},0,\frac{2 \sqrt{3}+9}{27} ,0,1\right)^\top\\
  \mathcal{E}_7 = \left(\frac{275 \sqrt{3}+243}{1494},\frac{9 \sqrt{3}+49}{166} ,\frac{275 \sqrt{3}+243}{1494},-\frac{9 \sqrt{3}+49}{166} ,-\frac{1}{3 \sqrt{3}},0,-1,0\right)^\top\\
  \mathcal{E}_8 = \left(\frac{9 \sqrt{3}+49}{166},\frac{81-19 \sqrt{3}}{1494},-\frac{9 \sqrt{3}+49}{166} ,\frac{81-19 \sqrt{3}}{1494},0,\frac{211 \sqrt{3}+162}{747} ,0,-1\right)^\top
\end{array}
\end{equation}
The corresponding eigenvalues of the matrix $\lambda \mathbb{I} + \mathfrak{M}^{-1} \mathfrak{B}$ are
\begin{equation}\label{eigenvalues}
  \frac{\mu_5}{\sqrt{3}}=0, ~~~~
  \frac{\mu_6}{\sqrt{3}}=0, ~~~~
  \frac{\mu_7}{\sqrt{3}}=\frac{799 \sqrt{3}+1233}{498}, ~~~~
 \frac{\mu_8}{\sqrt{3}}=\frac{799 \sqrt{3}+1233}{498}\nonumber
\end{equation}

\subsection{Rhombic Central Configurations}
\indent\par

Let us  investigate the degeneracy of rhombic central configurations in this subsection, that is, $t=0, u=\frac{1}{2}$ for (\ref{centralconfiguration2}). Then $m_3=m_4$, and the equations (\ref{centralconfiguration2}) of central configurations become:
\begin{equation}\label{centralconfiguration3}
\left\{
             \begin{array}{lr}
              \frac{1}{4s^3}+\frac{2\tilde{m}}{\tilde{r}^3}=\lambda  &  \\
\frac{2}{\tilde{r}^3}+ 2\tilde{m}=\lambda  &
             \end{array}
\right.
\end{equation}
where $\tilde{m}= m_3 = m_4, \tilde{r}=\sqrt{s^2+\frac{1}{4}}$.

Aa a result,  the $4\times 4$ matrix
$(\tilde{\mathcal{V}},   \tilde{ \mathcal{P}},\mathbf{i}\tilde{\mathcal{V}},   \mathbf{i}\tilde{\mathcal{P}})^{\top} (\lambda \mathfrak{M} + \mathfrak{B})  (\tilde{\mathcal{V}},   \tilde{ \mathcal{P}},\mathbf{i}\tilde{\mathcal{V}},   \mathbf{i}\tilde{\mathcal{P}})$ becomes
\begin{displaymath}
\left(
   \begin{array}{cccc}
     \frac{12\tilde{m}s^2}{\tilde{r}^5}+ \frac{3}{2s^3}-\frac{12 \lambda s^2}{2s^2+\frac{\tilde{m}}{2}} &  0 &   &   \\
     0 & (\frac{3\tilde{m}}{\tilde{r}^5}+\frac{1}{2s^3}-\frac{2\lambda}{\tilde{m}+1})(\frac{1}{\tilde{m}}+1)^2 &   &   \\
       &   & \frac{3\tilde{m}}{\tilde{r}^5} & 0 \\
      &   & 0 & (\frac{12\tilde{m}s^2}{\tilde{r}^5}+\frac{1}{2s^3}-\frac{2\lambda}{\tilde{m}+1})(\frac{1}{\tilde{m}}+1)^2 \\
   \end{array}
 \right)
\end{displaymath}
so the central configuration $\mathcal{E}_3$ is {degenerate} if and only if
\begin{equation}\label{degeneracy of rhombic central configurations}
(\frac{12\tilde{m}s^2}{\tilde{r}^5}+ \frac{3}{2s^3}-\frac{12 \lambda s^2}{2s^2+\frac{\tilde{m}}{2}})(\frac{3\tilde{m}}{\tilde{r}^5}+\frac{1}{2s^3}-\frac{2\lambda}{\tilde{m}+1})(\frac{12\tilde{m}s^2}{\tilde{r}^5}+\frac{1}{2s^3}-\frac{2\lambda}{\tilde{m}+1})=0
\end{equation}

Note that $\tilde{r}=1$ is impossible for the equations (\ref{centralconfiguration3}). Then, by  the equations (\ref{centralconfiguration3}), we have
\begin{equation}\label{mlambda}
\left\{
             \begin{array}{lr}
             \tilde{m}=\frac{\frac{\tilde{r}^3}{8s^3}-1}{\tilde{r}^3-1}, &  \\
             \lambda=\frac{2}{\tilde{r}^3}+ 2\tilde{m}= \frac{\frac{2\tilde{r}^3}{8s^3}-\frac{2}{\tilde{r}^3}}{\tilde{r}^3-1}.  &
             \end{array}
\right.
\end{equation}
To solve the equation (\ref{degeneracy of rhombic central configurations}), let us introduce the following rational transformations:
\begin{equation}\label{rational transformations}
\left\{
             \begin{array}{lr}
             \tilde{r}=\frac{\zeta^2+1}{4\zeta} &  \\
             s=\frac{\zeta^2-1}{4\zeta}  &
             \end{array}
\right.\nonumber
\end{equation}
According to $s>0$, we can think that $\zeta>1$. Then the equations (\ref{mlambda}) become
\begin{equation}\label{mlambda1}
\left\{
             \begin{array}{lr}
             \tilde{m}=\frac{-8 \zeta ^3 \left(\zeta ^2-3\right) \left(7 \zeta ^4-6 \zeta ^2+3\right)}{\left(\zeta ^2-1\right)^3 \left(\zeta ^2-4 \zeta +1\right) \left(\zeta ^4+4 \zeta ^3+18 \zeta ^2+4 \zeta +1\right)}, &  \\
             \lambda=\frac{16 \zeta ^3 \left(\zeta ^{12}+6 \zeta ^{10}-512 \zeta ^9+15 \zeta ^8+1536 \zeta ^7+20 \zeta ^6-1536 \zeta ^5+15 \zeta ^4+512 \zeta ^3+6 \zeta ^2+1\right)}{\left(\zeta ^4-1\right)^3 \left(\zeta ^6+3 \zeta ^4-64 \zeta ^3+3 \zeta ^2+1\right)},  &
             \end{array}
\right.
\end{equation}
and the equation (\ref{degeneracy of rhombic central configurations}) attributes to the following equation:
\begin{equation}\label{degeneracy of rhombic central configurations1}
\begin{aligned}
&[(\zeta ^2-3) (7 \zeta ^{10}-45 \zeta ^8+70 \zeta ^6+256 \zeta ^5-90 \zeta ^4+35 \zeta ^2-9) \\
&(7 \zeta ^{16}-88 \zeta ^{14}-448 \zeta ^{13}-44 \zeta ^{12}+12352 \zeta ^{11}+184 \zeta ^{10}-37504 \zeta ^9-70 \zeta ^8 \\
&+34176 \zeta ^7-296 \zeta ^6-13248 \zeta ^5-12 \zeta ^4+576 \zeta ^3+72 \zeta ^2-9) \\
&(17 \zeta ^{16}-56 \zeta ^{14}-2432 \zeta ^{13}-4 \zeta ^{12}+14720 \zeta ^{11}+248 \zeta ^{10}-32768 \zeta ^9+70 \zeta ^8\\
&+30720 \zeta ^7-136 \zeta ^6-14976 \zeta ^5+60 \zeta ^4+2688 \zeta ^3+72 \zeta ^2-15)]/[\left(\zeta ^2-4 \zeta +1\right)\\
&\left(\zeta ^{12}-120 \zeta ^9-3 \zeta ^8+408 \zeta ^7-360 \zeta ^5+3 \zeta ^4+136 \zeta ^3-1\right)\\
&\left(\zeta ^{12}-4 \zeta ^{10}-64 \zeta ^9+5 \zeta ^8+224 \zeta ^7-160 \zeta ^5-5 \zeta ^4+64 \zeta ^3+4 \zeta ^2-1\right)]
=0
\end{aligned}\nonumber
\end{equation}
By the  software Mathematica, all the solutions of above equation are the following£»
\begin{displaymath}
\zeta = \sqrt{3}, \zeta \approx 1.4943350238941826,\zeta \approx 5.0458391643884815.
\end{displaymath}
However, due to the equations (\ref{mlambda1}), the corresponding value of $\tilde{m}$ are all nonpositive:
\begin{displaymath}
\tilde{m} = 0, \tilde{m} \approx -1.4888656305411752,\tilde{m} \approx -0.6716522831120226.
\end{displaymath}

So the following theorem holds:
\begin{theorem}\label{Rhombic Central Configurations}
All the rhombic central configurations are nondegenerate.
\end{theorem}

Indeed, thanks to (\ref{mlambda1}), a routine computation gives rise to
\begin{displaymath}
\tilde{m} > 0 \Leftrightarrow \sqrt{3}<\zeta <\sqrt{3}+2.
\end{displaymath}
Therefore, all the admissible configurations to be central configurations are just
\begin{displaymath}
\begin{array}{lr}
            \left(-\frac{\zeta ^2-1}{4 \zeta },0,\frac{\zeta ^2-1}{4 \zeta },0,0,\frac{1}{2},0,-\frac{1}{2}\right), & ~~~for~ \sqrt{3}<\zeta <\sqrt{3}+2.
             \end{array}
\end{displaymath}

As a result of (\ref{Hessian}) and some tedious computation, the  eigenvalues of the matrix $\lambda \mathbb{I} + \mathfrak{M}^{-1} \mathfrak{B}$ are
\begin{equation}
  \frac{\mu_1}{\sqrt{I}}=\lambda, ~~~~
  \frac{\mu_2}{\sqrt{I}}=\lambda, ~~~~
  3\lambda, ~~~~
 0 \nonumber
\end{equation}
and
\begin{equation}\label{Rhombiceigenvalues}
\begin{array}{lr}
             \frac{\mu_5}{\sqrt{I}}=-\frac{48 \zeta ^3 \left(7 \zeta ^{10}-45 \zeta ^8+70 \zeta ^6+256 \zeta ^5-90 \zeta ^4+35 \zeta ^2-9\right)}{\left(\zeta ^2-1\right)^3 \left(\zeta ^2+1\right)^2 \left(\zeta ^6+3 \zeta ^4-64 \zeta ^3+3 \zeta ^2+1\right)},  & \\
             \frac{\mu_6}{\sqrt{I}}=\frac{384 \zeta ^3 \left(\zeta ^{12}-4 \zeta ^{10}-64 \zeta ^9+5 \zeta ^8+224 \zeta ^7-160 \zeta ^5-5 \zeta ^4+64 \zeta ^3+4 \zeta ^2-1\right)}{\left(\zeta ^2-1\right)^3 \left(\zeta ^2+1\right)^3 \left(\zeta ^6+3 \zeta ^4-64 \zeta ^3+3 \zeta ^2+1\right)},  & \\
             \begin{aligned}\frac{\mu_7}{\sqrt{I}}=
&16 \zeta ^3 (7 \zeta ^{16}-88 \zeta ^{14}-448 \zeta ^{13}-44 \zeta ^{12}+12352 \zeta ^{11}+184 \zeta ^{10}-37504 \zeta ^9\\
&-70 \zeta ^8+34176 \zeta ^7-296 \zeta ^6-13248 \zeta ^5-12 \zeta ^4+576 \zeta ^3+72 \zeta ^2-9)\\
&/\left(1-\zeta ^2\right)^3 \left(\zeta ^2+1\right)^5 \left(\zeta ^6+3 \zeta ^4-64 \zeta ^3+3 \zeta ^2+1\right),
\end{aligned}
               & \\
             \begin{aligned}\frac{\mu_8}{\sqrt{I}}=
&16 \zeta ^3 (17 \zeta ^{16}-56 \zeta ^{14}-2432 \zeta ^{13}-4 \zeta ^{12}+14720 \zeta ^{11}+248 \zeta ^{10}-32768 \zeta ^9\\
&+70 \zeta ^8+30720 \zeta ^7-136 \zeta ^6-14976 \zeta ^5+60 \zeta ^4+2688 \zeta ^3+72 \zeta ^2-15)\\
&/\left(\zeta ^2-1\right)^3 \left(\zeta ^2+1\right)^5 \left(\zeta ^6+3 \zeta ^4-64 \zeta ^3+3 \zeta ^2+1\right).
\end{aligned}  &
             \end{array}
\end{equation}
where\begin{equation}
I=\frac{\left(\zeta ^2+1\right)^2 \left(\zeta ^{12}-4 \zeta ^{10}-64 \zeta ^9+5 \zeta ^8+224 \zeta ^7-160 \zeta ^5-5 \zeta ^4+64 \zeta ^3+4 \zeta ^2-1\right)}{8 \zeta ^2 \left(\zeta ^2-1\right)^3 \left(\zeta ^6+3 \zeta ^4-64 \zeta ^3+3 \zeta ^2+1\right)}.\nonumber
\end{equation}

A straightforward computation shows that
\begin{equation}
\begin{array}{lr}
             \frac{\kappa^{\frac{1}{2}}}{\sqrt{I}}=\frac{\sqrt{2\lambda}}{\sqrt{I}}=  & \\
16 \sqrt{\frac{\zeta ^5 \left(\zeta ^{12}+6 \zeta ^{10}-512 \zeta ^9+15 \zeta ^8+1536 \zeta ^7+20 \zeta ^6-1536 \zeta ^5+15 \zeta ^4+512 \zeta ^3+6 \zeta ^2+1\right)}{\left(\zeta ^2+1\right)^5 \left(\zeta ^{12}-4 \zeta ^{10}-64 \zeta ^9+5 \zeta ^8+224 \zeta ^7-160 \zeta ^5-5 \zeta ^4+64 \zeta ^3+4 \zeta ^2-1\right)}} ,&
             \end{array}\nonumber
\end{equation}
and
\begin{equation}\label{Rhombiceigenvalues1}
\begin{array}{lr}
              0< \kappa^{\frac{1}{2}},\mu_5<\mu_6 <\mu_7,\mu_8, & ~~for~~ \sqrt{3}<\zeta <\sqrt{3}+2\\
               \lambda_5<\sqrt{\kappa },  & ~~for~~ \zeta_1<\zeta <\zeta_2\\
               \lambda_5=\sqrt{\kappa },  & ~~for~~ \zeta =\zeta_1, \zeta_2\\
               \lambda_5>\sqrt{\kappa },  & ~~for~~ \sqrt{3}<\zeta<\zeta_1 ~or ~<\zeta_2<\zeta<\sqrt{3}+2\\
              \mu_8 <\mu_7,  & ~~for~~ \sqrt{3}<\zeta <\sqrt{2}+1\\
              \mu_7 <\mu_8,  & ~~for~~ \sqrt{2}+1<\zeta <\sqrt{3}+2\\
              \mu_7 =\mu_8,  & ~~for~~ \zeta =\sqrt{2}+1
             \end{array}
\end{equation}
where$\zeta_1 \approx 1.7889580612081344, \zeta_2 \approx 3.705602221466667$.

It is clear that all the $\kappa^{\frac{1}{2}},\mu_5, \mu_6, \mu_7, \mu_8$ and their ratios are irrational for almost all $\zeta \in(\sqrt{3},\sqrt{3}+2)$, since each of them is rational only for countable $\zeta$.

\subsection{Central Configurations With Two Degrees of Degeneracy}
\indent\par
Let us  investigate the  central configurations with degree of degeneracy two  in this subsection. We only consider the cases of except equilateral central configurations.

Thus the equations (\ref{centralconfiguration2}) can be written as
\begin{equation}\label{massbiaoshi}
\left\{
             \begin{array}{lr}
              m_3= \frac{2\left(u+t-1\right)\left(\frac{1}{8s^3}-\frac{1}{[s^2+(u+t-1)^2]^\frac{3}{2}}\right)}{\frac{1}{[s^2+(u+t)^2]^\frac{3}{2}}-1}&\\
              m_4= \frac{2\left(u+t\right)\left(\frac{1}{8s^3}-\frac{1}{[s^2+(u+t)^2]^\frac{3}{2}}\right)}{1-\frac{1}{[s^2+(u+t-1)^2]^\frac{3}{2}}}&\\
              \lambda = \frac{1}{4s^3}+\frac{m_3}{[s^2+(u+t)^2]^\frac{3}{2}}+\frac{m_4}{[s^2+(u+t-1)^2]^\frac{3}{2}}  &
             \end{array}
\right.\nonumber
\end{equation}

To study the degeneracy of central configurations, let us introduce the following rational transformations:
\begin{equation}\label{rational transformations1}
\left\{
             \begin{array}{lr}
             u+t =s \frac{\xi^2-1}{2\xi} &  \\
             u+t-1=s \frac{\eta^2-1}{2\eta}  &
             \end{array}
\right.\nonumber
\end{equation}
Without losing generality, we may assume that $\xi>\eta>0$.
Then
\begin{equation}\label{rational transformations2}
\left\{
             \begin{array}{lr}
             \sqrt{s^2+(u+t)^2} =s \frac{\xi^2+1}{2\xi} &  \\
             \sqrt{s^2+(u+t-1)^2}=s \frac{\eta^2+1}{2\eta}  &\\
             s= \frac{1}{\frac{\xi^2-1}{2\xi}-\frac{\eta^2-1}{2\eta}}  &
             \end{array}
\right.\nonumber
\end{equation}
and
\begin{equation}\label{rational transformations3}
\left\{
             \begin{array}{lr}
             m_3= \frac{-(\eta -1) (\eta +1) \left(\eta ^2-4 \eta +1\right) \left(\eta ^4+4 \eta ^3+18 \eta ^2+4 \eta +1\right) \left(\xi ^2+1\right)^3 (\eta -\xi )^2 (\eta  \xi +1)^2}{32 \left(\eta ^2+1\right)^3 \xi ^2 \left(\eta ^2 \xi +2 \eta -\xi \right) \left(\eta ^4 \xi ^2-3 \eta ^3 \xi ^3+\eta ^3 \xi +3 \eta ^2 \xi ^4-2 \eta ^2 \xi ^2+\eta ^2+3 \eta  \xi ^3-\eta  \xi +\xi ^2\right)} &\\
              m_4=\frac{\left(\eta ^2+1\right)^3 (\xi -1) (\xi +1) \left(\xi ^2-4 \xi +1\right) \left(\xi ^4+4 \xi ^3+18 \xi ^2+4 \xi +1\right) (\eta -\xi )^2 (\eta  \xi +1)^2}{32 \eta ^3 \left(\xi ^2+1\right)^3 \left(2 \eta  \xi -\xi ^2+1\right) \left(\eta ^4 \xi ^2-\eta ^3 \xi ^3+\eta ^3 \xi +\eta ^2 \xi ^4-2 \eta ^2 \xi ^2+\eta ^2+3 \eta  \xi ^3-3 \eta  \xi +3 \xi ^2\right)}&\\
              \lambda = \frac{1}{4s^3}+\frac{8\xi^3m_3}{(\xi^2+1)^3s^3}+\frac{8\eta^3m_4}{(\eta^2+1)^3s^3}  &
             \end{array}
\right.\nonumber
\end{equation}

If $\xi=1$ or $\eta=1$, then the corresponding configurations are formed by the three particles  on a common straight line but the fourth particle not on the straight line, by the well-known \emph{perpendicular bisector theorem} \cite{moeckel1990central}, this kind of configurations cannot be central configurations.

If $\xi \eta=1$, then the corresponding configurations are rhombus.

If $\xi=2\pm\sqrt{3}$ , then $m_4=0$ or the corresponding configurations are equilateral  configurations; similarly,  if $\eta=2\pm\sqrt{3}$, then $m_3=0$ or the corresponding configurations are equilateral  configurations.

Therefore, we investigate only the cases of $\xi, \eta \neq 1, 2\pm\sqrt{3}$ and $\xi \eta \neq 1$ in the following.

Note that
\begin{displaymath}
(\mathbf{i}\tilde{\mathcal{V}})^{\top} (\lambda \mathfrak{M} + \mathfrak{B})  \mathbf{i}\tilde{\mathcal{V}}=\frac{6m_3(u+t)^2}{[s^2+(u+t)^2]^\frac{5}{2}}+\frac{6m_4(u+t-1)^2}{[s^2+(u+t-1)^2]^\frac{5}{2}}\rangle0,
\end{displaymath}
hence the central configurations with  degree of degeneracy two  satisfy the following equations
\begin{equation}\label{quanweiling}
(\tilde{\mathcal{V}},   \tilde{ \mathcal{P}})^{\top} (\lambda \mathfrak{M} + \mathfrak{B})  (\tilde{\mathcal{V}},   \tilde{ \mathcal{P}})=\left(
                                                                                                                                          \begin{array}{cc}
                                                                                                                                            0 & 0 \\
                                                                                                                                            0 & 0 \\
                                                                                                                                          \end{array}
                                                                                                                                        \right)
\end{equation}
or
\begin{equation}\label{dierzhong}
\left\{
             \begin{array}{lr}
             Det \left[(\tilde{\mathcal{V}},   \tilde{ \mathcal{P}})^{\top} (\lambda \mathfrak{M} + \mathfrak{B})  (\tilde{\mathcal{V}},   \tilde{ \mathcal{P}}) \right]=0&  \\
             Det\left[(\mathbf{i}\tilde{\mathcal{V}},   \mathbf{i}\tilde{\mathcal{P}})^{\top} (\lambda \mathfrak{M} + \mathfrak{B})  (\mathbf{i}\tilde{\mathcal{V}},   \mathbf{i}\tilde{\mathcal{P}}) \right]=0  &
             \end{array}
\right.
\end{equation}

By the  software Mathematica, there is no  solution of  the equations (\ref{quanweiling}) (\ref{dierzhong})  such that $m_3>0$ and $m_4>0$. It is necessary to contain computer-aided proofs for this  part, due to the size of some of the polynomials we worked with. We will not write all
the explicit expressions of the polynomials. Instead, we shall provide the
steps followed to calculate all the important polynomials.

By the above rational transformations, the equations (\ref{quanweiling}) (\ref{dierzhong}) can respectively be written as
\begin{equation}\label{rational1}
\left\{
             \begin{array}{lr}
             Q_{11}(\xi, \eta)=0&  \\
             Q_{12}(\xi, \eta)=0  &\\
             Q_{22}(\xi, \eta)=0 &
             \end{array}
\right.
\end{equation}
and
\begin{equation}\label{rational2}
\left\{
             \begin{array}{lr}
             Q_1(\xi, \eta)=Q_{11}(\xi, \eta)Q_{22}(\xi, \eta)-Q_{12}^2(\xi, \eta)=0&  \\
             Q_2(\xi, \eta)=0  &
             \end{array}
\right.
\end{equation}
where $Q_{11}, Q_{22}, Q_{12}, Q_1, Q_2$ are rational functions of $\xi, \eta$.

We mainly consider the numerators of above  rational functions. Note that
\begin{equation}\label{biggerzero}
\left\{
             \begin{array}{lr}
             m_3 >0 \Leftrightarrow (1 - \eta)\left(\eta ^2-4 \eta +1\right)\left(\eta ^2 \xi +2 \eta -\xi \right)>0 &\\
              m_4>0 \Leftrightarrow (\xi -1) \left(\xi ^2-4 \xi +1\right)\left(2 \eta  \xi -\xi ^2+1\right)>0 &
             \end{array}
\right.
\end{equation}
and
\begin{displaymath}
0<\frac{2}{m_3+m_4}+1=\frac{P_{numer}(\xi, \eta)}{P_{denom}(\xi, \eta)},
\end{displaymath}
where $P_{numer}$ is a polynomial which is the numerator of $\frac{2}{m_3+m_4}+1$, and $P_{denom}$ is a polynomial which is the denominator of $\frac{2}{m_3+m_4}+1$; please pay particular attention to the polynomial
\begin{equation}
\begin{aligned}
&P_{numer}=\\
&(\xi ^6 \eta ^{10}-\xi ^4 \eta ^{10}-2 \xi ^7 \eta ^9+4 \xi ^5 \eta ^9-2 \xi ^3 \eta ^9-\xi ^8 \eta ^8-8 \xi ^6 \eta ^8-4 \xi ^4 \eta ^8-3 \xi ^2 \eta ^8+2 \xi ^9 \eta ^7\\
&+8 \xi ^5 \eta ^7-8 \xi ^3 \eta ^7-2 \xi  \eta ^7-\xi ^{10} \eta ^6+2 \xi ^8 \eta ^6-9 \xi ^6 \eta ^6+13 \xi ^4 \eta ^6-4 \xi ^2 \eta ^6-\eta ^6-4 \xi ^9 \eta ^5-\\
&8 \xi ^7 \eta ^5+8 \xi ^3 \eta ^5+4 \xi  \eta ^5+\xi ^{10} \eta ^4-2 \xi ^8 \eta ^4-31 \xi ^6 \eta ^4-9 \xi ^4 \eta ^4-8 \xi ^2 \eta ^4+\eta ^4+2 \xi ^9 \eta ^3+\\
&8 \xi ^7 \eta ^3-8 \xi ^5 \eta ^3-2 \xi  \eta ^3+\xi ^8 \eta ^2-2 \xi ^6 \eta ^2+2 \xi ^4 \eta ^2-\xi ^2 \eta ^2+2 \xi ^7 \eta -4 \xi ^5 \eta +2 \xi ^3 \eta +\xi ^6-\\
&\xi ^4)(\xi ^9 \eta ^{12}+3 \xi ^7 \eta ^{12}-64 \xi ^6 \eta ^{12}+3 \xi ^5 \eta ^{12}+\xi ^3 \eta ^{12}-3 \xi ^{10} \eta ^{11}-6 \xi ^8 \eta ^{11}+192 \xi ^7 \eta ^{11}-\\
&192 \xi ^5 \eta ^{11}+6 \xi ^4 \eta ^{11}+3 \xi ^2 \eta ^{11}+3 \xi ^{11} \eta ^{10}+3 \xi ^9 \eta ^{10}-192 \xi ^8 \eta ^{10}-6 \xi ^7 \eta ^{10}+384 \xi ^6 \eta ^{10}-
\end{aligned}\nonumber
\end{equation}
\begin{equation}
\begin{aligned}
&6 \xi ^5 \eta ^{10}-192 \xi ^4 \eta ^{10}+3 \xi ^3 \eta ^{10}+3 \xi  \eta ^{10}-\xi ^{12} \eta ^9-3 \xi ^{10} \eta ^9-64 \xi ^9 \eta ^9-3 \xi ^8 \eta ^9-384 \xi ^7 \eta ^9\\
&-384 \xi ^5 \eta ^9+3 \xi ^4 \eta ^9-192 \xi ^3 \eta ^9+3 \xi ^2 \eta ^9+\eta ^9+6 \xi ^{11} \eta ^8+192 \xi ^{10} \eta ^8+3 \xi ^9 \eta ^8-21 \xi ^7 \eta ^8\\
&+960 \xi ^6 \eta ^8-21 \xi ^5 \eta ^8-768 \xi ^4 \eta ^8+3 \xi ^3 \eta ^8-192 \xi ^2 \eta ^8+6 \xi  \eta ^8-3 \xi ^{12} \eta ^7-192 \xi ^{11} \eta ^7+\\
&6 \xi ^{10} \eta ^7+21 \xi ^8 \eta ^7-576 \xi ^7 \eta ^7+1344 \xi ^5 \eta ^7-21 \xi ^4 \eta ^7-384 \xi ^3 \eta ^7-6 \xi ^2 \eta ^7-192 \xi  \eta ^7+3 \eta ^7\\
&+64 \xi ^{12} \eta ^6-384 \xi ^{10} \eta ^6-960 \xi ^8 \eta ^6+960 \xi ^4 \eta ^6+384 \xi ^2 \eta ^6-64 \eta ^6-3 \xi ^{12} \eta ^5+192 \xi ^{11} \eta ^5\\
&+6 \xi ^{10} \eta ^5+21 \xi ^8 \eta ^5-2496 \xi ^7 \eta ^5-576 \xi ^5 \eta ^5-21 \xi ^4 \eta ^5-384 \xi ^3 \eta ^5-6 \xi ^2 \eta ^5+192 \xi  \eta ^5+\\
&3 \eta ^5-6 \xi ^{11} \eta ^4+192 \xi ^{10} \eta ^4-3 \xi ^9 \eta ^4+768 \xi ^8 \eta ^4+21 \xi ^7 \eta ^4-960 \xi ^6 \eta ^4+21 \xi ^5 \eta ^4-3 \xi ^3 \eta ^4\\
&-192 \xi ^2 \eta ^4-6 \xi  \eta ^4-\xi ^{12} \eta ^3-3 \xi ^{10} \eta ^3+64 \xi ^9 \eta ^3-3 \xi ^8 \eta ^3+3 \xi ^4 \eta ^3-64 \xi ^3 \eta ^3+3 \xi ^2 \eta ^3+\eta ^3\\
&-3 \xi ^{11} \eta ^2-3 \xi ^9 \eta ^2+192 \xi ^8 \eta ^2+6 \xi ^7 \eta ^2-384 \xi ^6 \eta ^2+6 \xi ^5 \eta ^2+192 \xi ^4 \eta ^2-3 \xi ^3 \eta ^2-3 \xi  \eta ^2-\\
&3 \xi ^{10} \eta -6 \xi ^8 \eta +192 \xi ^7 \eta -192 \xi ^5 \eta +6 \xi ^4 \eta +3 \xi ^2 \eta -\xi ^9-3 \xi ^7+64 \xi ^6-3 \xi ^5-\xi ^3).
\end{aligned}\nonumber
\end{equation}
It is easy to see that $P_{numer}$ and $P_{denom}$ cannot simultaneously be zero
except the cases abandoned by us. So  the polynomial $P_{numer}$ cannot  be zero in the scope of our interest.

Then, by using factorization and combining what has been said above, especially, one should observe that $P_{numer}$ is a factor of numerators of $ Q_{22}, Q_{12}, Q_1, Q_2$,  we can respectively reduce the equations (\ref{rational1}) and (\ref{rational2}) to
\begin{equation}\label{polynomial1}
\left\{
             \begin{array}{lr}
             P_{11}(\xi, \eta)=0&  \\
             P_{12}(\xi, \eta)=0  &\\
             P_{22}(\xi, \eta)=0 &
             \end{array}
\right.\nonumber
\end{equation}
and
\begin{equation}\label{polynomial2}
\left\{
             \begin{array}{lr}
             P_1(\xi, \eta)=0&  \\
             P_2(\xi, \eta)=0  &
             \end{array}
\right.\nonumber
\end{equation}
where $P_{11}, P_{22}, P_{12}, P_1, P_2$ are polynomial functions of $\xi, \eta$.

Then one can easily seek the solutions of the above equations one by one combining the conditions (\ref{biggerzero}). It is noteworthy that using the resultant of two polynomials can evidently economize calculating time.   We finally see that there is not a solution of  the equations (\ref{quanweiling}) (\ref{dierzhong})  such that $m_3>0$ and $m_4>0$.

So the following theorem holds:
\begin{theorem}\label{Rhombic Central Configurations}
If a kite central configuration except equilateral central configurations is degenerate, then the degree of degeneracy is one.
\end{theorem}

Following this result and a crude dimension count, we can venturesomely conjecture that:

\emph{ All the four-body central configurations except equilateral central configurations have zero or one degree of degeneracy.}

That is to say, we believe that the equilateral central configurations  founded by Palmore \cite{palmore1975classifying,Palmore1976Measure} are  the only degenerate central configurations with  degree of degeneracy two  for the four-body problem, although we cannot prove it  now. This conjecture is important for resolving $PISPW$ of the four-body problem. Indeed, if
this conjecture is true, due to Corollary \ref{lastcorollary}, one can completely solve  $PISPW$ of the four-body problem.

\begin{remark}
The results and method in this section are novel according to what I know,  even though for the kind of kite central configurations, Leandro \cite{Leandro2003Finiteness} has investigated a more general case,  included works on a kind of degenerate central configurations. Indeed the above rational transformations are inspired by his work. However, the definition of degeneracy in his work is different from ours, thus the results in \cite{Leandro2003Finiteness} cannot be directly applied. The point is that the degenerate central configuration by his definition is still degenerate by our definition; the opposite, however, isn't necessarily true. The problem can be explained in this way:

Given a function $f$ defined on a manifold $\Omega$, $x\in \Omega$ is a critical point of the function $f$. Let $\Omega_1$ be a submanifold of $\Omega$ and $x\in \Omega_1$,  although  the point $x$ may  still be a critical point of the function $f$ restricted to $\Omega_1$, the degeneracy of the critical point $x$ may change: $x$ is degenerate on $\Omega_1$ certainly  implies  that $x$ is degenerate on $\Omega$, but the converse isn't necessarily true.
\end{remark}


\section{Equations of Motion in Various Coordinate Systems}\label{Lagrangiansystems}

\subsection{Lagrangian Dynamical Systems}\label{LagrangianDynamicalsystems}
      \indent\par
For the sake of readability, we sketchily give the theory of Lagrangian dynamical systems used in deducing the general equations of motion. The exposition follows \cite{zbMATH03601149,zbMATH05031968}, to which we refer the reader for proofs and details.

\begin{definition}
Let $\mathcal{M}$ be a differentiable manifold, $\mathrm{T}\mathcal{M}$ its tangent bundle, and $\mathcal{L}: \mathrm{T}\mathcal{M} \rightarrow \mathbb{R}$
a differentiable function. A map $\gamma: \mathbb{R} \rightarrow \mathcal{M}$ is called a motion in the Lagrangian
system with configuration manifold $\mathcal{M}$ and Lagrangian function $\mathcal{L}$ if $\gamma$ is an
extremal of the Lagrangian action functional
\begin{equation}
\mathcal{A}(\gamma) = \int^{t_2}_{t_1}{ \mathcal{L}(\dot{\gamma}(t)) dt}.\nonumber
\end{equation}
where $\dot{\gamma}$ is the velocity vector $\dot{\gamma}(t) \in \mathrm{T}_{\gamma(t)}\mathcal{M}$.
\end{definition}

\begin{theorem}\label{Euler-Lagrange
equations}
The evolution of the local coordinates $q = (q_1, \cdots, q_n)$ of a point $\gamma(t)$
under motion in a Lagrangian system on a manifold satisfies the \emph{Euler-Lagrange
equations}
\begin{equation}
\frac{d}{dt}\frac{\partial \mathcal{L}}{\partial \dot{q}} = \frac{\partial \mathcal{L}}{\partial {q}} ,\nonumber
\end{equation}
where $\mathcal{L}(q, \dot{q})$ is the expression for the function $\mathcal{L}: \mathrm{T}\mathcal{M} \rightarrow \mathbb{R}$ in the coordinates
$q$ and $\dot{q}$ on $\mathrm{T}\mathcal{M}$.
\end{theorem}

Theorem \ref{Euler-Lagrange
equations} yields a quick method
for writing equations of motion in various coordinate systems, even in larger class of coordinate transformations which  contain time. Indeed, to write the equations of motion in a new coordinate system, it is sufficient to
express the Lagrangian function in the new coordinates. In fact, we have
\begin{theorem}
If the orbit $\gamma: q = \varphi(t)$ of Euler-Lagrange equations $\frac{d}{dt}\frac{\partial \mathcal{L}}{\partial \dot{q}} = \frac{\partial \mathcal{L}}{\partial {q}}$ is written as $\gamma: Q = \Phi(t)$ in the local coordinates $Q, t$ (where $Q =
Q(q, t)$), then the function $\Phi(t)$ satisfies Euler-Lagrange equations $\frac{d}{dt}\frac{\partial \tilde{\mathcal{L}}}{\partial \dot{Q}} = \frac{\partial \tilde{\mathcal{L}}}{\partial {Q}}$, where $\tilde{\mathcal{L}}(Q, \dot{Q}, t) = \mathcal{L}(q,\dot{q}, t)$.

\end{theorem}

\begin{remark}
By the additional dependence of the
Lagrangian function on time:
\begin{equation}
\mathcal{L}: \mathrm{T}\mathcal{M}\times \mathbb{R} \rightarrow \mathbb{R}   ~~~~~~\mathcal{L}=\mathcal{L}(q,\dot{q},t),\nonumber
\end{equation}
one can consider a Lagrangian nonautonomous system and the results above are also valid.
  \end{remark}

\begin{definition}
In mechanics,  $\frac{\partial \mathcal{L}}{\partial \dot{q}}$ are called generalized
momenta, $\frac{\partial \mathcal{L}}{\partial {q}}$ are called generalized forces.
\end{definition}

\begin{definition}
Given a Lagrangian function $\mathcal{L}(q,\dot{q},t)$,  a coordinate $q_j$ is called ignorable (or cyclic ) if it does not enter into the
Lagrangian:  $\frac{\partial \mathcal{L}}{\partial q_j}=0$.
\end{definition}

\begin{theorem}
The generalized momentum corresponding to an ignorable coordinate is
conserved: $p_j = \frac{\partial \mathcal{L}}{\partial \dot{q}_j}=const$.
\end{theorem}

The following  content is \emph{Routh's method for eliminating
ignorable coordinates}.\\
Suppose that the Lagrangian $\mathcal{L}(q,\dot{q},\dot{\xi})$ does not involve the coordinate $\xi$, i.e., $\xi$ is ignorable. Using the equality $\frac{\partial \mathcal{L}}{\partial \dot{\xi}}=c$ we represent the  velocity
$\dot{\xi}$ as a function of $q,\dot{q}$ and $c$. Following Routh we introduce the function
\begin{equation}
\mathcal{R}_c(q,\dot{q})=\mathcal{L}(q,\dot{q},\dot{\xi}) - c \dot{\xi}|_{q,\dot{q},c} (= \mathcal{L}(q,\dot{q},\dot{\xi}) - \frac{\partial \mathcal{L}}{\partial \dot{\xi}}  \dot{\xi}|_{q,\dot{q},c} ).\nonumber
\end{equation}

\begin{theorem}
A vector-function $(q(t),\xi(t))$ with the constant value of generalized momentum $\frac{\partial \mathcal{L}}{\partial \dot{\xi}}=c$ satisfies the Euler-Lagrange
equations $
\frac{d}{dt}\frac{\partial \mathcal{L}}{\partial \dot{q}} = \frac{\partial \mathcal{L}}{\partial {q}}$
 if and only
if $q(t)$ satisfies
the Euler-Lagrange
equations $\frac{d}{dt}\frac{\partial \mathcal{R}_c}{\partial \dot{q}} = \frac{\partial \mathcal{R}_c}{\partial {q}}. $
\end{theorem}

{\bf Example. } Compare Newton's equations (\ref{eq:Newton's equation1})
\begin{equation*}
\frac{d{m_k \dot{\mathbf{r}}_k}}{dt} - \frac{\partial \mathcal{U}}{\partial \mathbf{r}_k}=0,
\end{equation*}
with the Euler-Lagrange equations
\begin{equation}
\frac{d}{dt}\frac{\partial \mathcal{L}}{\partial \dot{q}} - \frac{\partial \mathcal{L}}{\partial {q}}=0 \nonumber
\end{equation}
of  Lagrangian
system with configuration manifold $\mathcal{X}$ and Lagrangian function
\begin{displaymath}
\mathcal{L}  =  \mathcal{K} + \mathcal{U}.
\end{displaymath}
\begin{theorem}
Motions of the mechanical system (\ref{eq:Newton's equation1}) coincide with extremals of
the functional $\mathcal{A}(\mathbf{r})= \int^{t_2}_{t_1}{ \mathcal{L}dt}.$

\end{theorem}

\subsection{Equations of Motion in the Moving Coordinate System}\label{EquationsMotionMovingCoordinateSystem}
      \indent\par
The equations (\ref{Euler-Lagrangeequations0}) of motion in the coordinates $ r, \theta, z_5, \cdots, z_{2N}$ are:
\begin{equation}\label{Euler-Lagrangeequations0detail}
\left\{
             \begin{array}{l}
             {r}^2 [ \ddot{z}_i + \ddot{\theta} \sum_{j=5}^{2N}q_{ij} z_j+ 2 \dot{\theta} \sum_{j=5}^{2N}q_{ij} \dot{z}_j]
             +{r}^2 [\frac{z_i \sum_{j=5}^{2N}(\ddot{z}_j z_j + \dot{z}^2_j)}{z^2_3} + \frac{{z}_i (\sum_{j=5}^{2N}{z}_j \dot{z}_j)^2}{z^4_3}]\\
 + 2r\dot{r}[\frac{z_i \sum_{j=5}^{2N}\dot{z}_j z_j}{z^2_3} + \dot{z}_i + \dot{\theta} \sum_{j=5}^{2N}q_{ij} z_j] - \frac{1}{r}\frac{\partial U(z)}{\partial z_i}=0,  ~~~~~ i=5,\cdots, 2N, \\
             \ddot{r}- {r} [\frac{ (\sum_{j=5}^{2N}{z}_j \dot{z}_j)^2}{z^2_3}+\sum_{j=5}^{2N}\dot{z}^2_j + 2\dot{\theta}\sum_{j,k=5}^{2N}q_{jk}\dot{z}_j {z}_k+ \dot{\theta}^2]+ \frac{ U(z)}{r^2}=0,   \\
             2r\dot{r} (\sum_{j,k=5}^{2N}q_{jk}\dot{z}_j {z}_k+ \dot{\theta}) + r^2 (\sum_{j,k=5}^{2N}q_{kj}\ddot{z}_k {z}_j+ \ddot{\theta})=0.
             \end{array}
\right.
\end{equation}
The equations (\ref{Euler-Lagrangeequations1}) of motion on the level set of $\mathcal{J}$ are:
\begin{equation}\label{Euler-Lagrangeequations1detail}
\left\{
             \begin{array}{lr}
              {r}^2 [\ddot{z}_i - \sum_{j=5}^{2N}q_{ij} z_j\sum_{j,k=5}^{2N}q_{jk}\ddot{z}_j {z}_k- 2  \sum_{j=5}^{2N}q_{ij} \dot{z}_j\sum_{j,k=5}^{2N}q_{jk}\dot{z}_j {z}_k]\\
               +{r}^2 [\frac{z_i \sum_{j=5}^{2N}(\ddot{z}_j z_j + \dot{z}^2_j)}{z^2_3} + \frac{{z}_i (\sum_{j=5}^{2N}{z}_j \dot{z}_j)^2}{z^4_3}]+ 2\mathcal{J}\sum_{j=5}^{2N}q_{ij} \dot{z}_j\\
 +2r\dot{r}[\frac{z_i \sum_{j=5}^{2N}\dot{z}_j z_j}{z^2_3} + \dot{z}_i -  \sum_{j=5}^{2N}q_{ij} z_j\sum_{j,k=5}^{2N}q_{jk}\dot{z}_j {z}_k] - \frac{1}{r}\frac{\partial U(z)}{\partial z_i}=0,  & i=5,\cdots, 2N, \\
             \ddot{r}- {r} [\frac{ (\sum_{j=5}^{2N}{z}_j \dot{z}_j)^2}{z^2_3}+\sum_{j=5}^{2N}\dot{z}^2_j - (\sum_{j,k=5}^{2N}q_{jk}\dot{z}_j {z}_k)^2+\frac{\mathcal{J}^2}{r^4}]+ \frac{ U(z)}{r^2}=0.  &
             \end{array}
\right.
\end{equation}
Note that $z_3 = \sqrt{1 - \sum_{j = 5}^{2N} z^2_j}$.

\begin{remark}
By Legendre Transforms, the systems above can respectively be reduced to the corresponding Hamiltonian  systems, although  the amount of calculation  is large in doing transformations. The benefits are that the theory of KAM is successfully applied to study  the  $N$-body problem near relative equilibrium solutions.
\end{remark}

\section{Some Notes on Equations of Motion for Collision Orbits}\label{MotionforCollisionOrbits}

\subsection{On Gradient-like Flow}\label{Gradient-likeFlow}
 \indent\par

The main aim of the subsection is to show the proposition below.

\begin{proposition}\label{gradient-likeflowr0}
Consider the system (\ref{equationzr3}),  we claim that the flow of the system (\ref{equationzr3}) restricted to the invariant manifold
 \begin{center}
 $\{r=0\}\bigcap\mathcal{N}_{\mathcal{H}}=\{(z, Z, 0, \Upsilon)\in \mathcal{N}_{\mathcal{H}}|{\Upsilon}^2+ [\frac{ (\sum_{j=5}^{2N}{z}_j Z_j)^2}{z^2_3}+\sum_{j=5}^{2N}Z^2_j - (\sum_{j,k=5}^{2N}q_{jk}Z_j {z}_k)^2]= 2U(z)\}$
 \end{center}
 is gradient-like with respect
to
\begin{equation}
    \pi_{\Upsilon}: \{r=0\}\bigcap\mathcal{N}_{\mathcal{H}}\rightarrow \mathbb{R}: (z, Z, 0, \Upsilon)\mapsto -\Upsilon.\nonumber
\end{equation}
\end{proposition}

Before  proving Proposition \ref{gradient-likeflowr0},
let us recall the concept of ``gradient-like" \cite{Mcgehee1974Triple}:
\begin{definition}
Let $\phi$ be a flow on a complete metric space $\mathcal{M}$. Suppose
there is a continuous function $f: \mathcal{M}\rightarrow \mathbb{R}$ such that
\begin{equation*}
    f(\phi(t,p))<f(p)  ~~~~~  ~~~~~  ~~~~~  if  ~~~ t>0
\end{equation*}
unless $p$ is an equilibrium point of $\phi$. Suppose further that the equilibrium points of $\phi$ are
isolated. Then $\phi$ is called gradient-like (with respect to $f$).
\end{definition}

{\bf Proof.} 
Proposition \ref{gradient-likeflowr0} will be proved if we can show that $\Upsilon$ is strictly monotone increasing with respect
to $\tau$. To this end,  it suffices to prove that
\begin{equation}
 \Upsilon'=   \frac{1}{2}\Upsilon^2+ \frac{ (\sum_{j=5}^{2N}{z}_j  {Z_j})^2}{z^2_3}+\sum_{j=5}^{2N}{Z_j}^2 - (\sum_{j,k=5}^{2N}q_{jk}{Z_j} {z}_k)^2 -   U(z)>0\nonumber
\end{equation}
or
\begin{equation}
\Upsilon'= 0,   ~~~~~
  \Upsilon''=0,  ~~~~~
  \Upsilon'''>0\nonumber
\end{equation}
 unless $(z, Z, 0, \Upsilon)$ is an equilibrium point of the system (\ref{equationzr3}).

First, it is easy to see that equilibrium points of the system (\ref{equationzr3}) are exactly points in the set
\begin{equation*}
   CP:= \{(z, Z, r, \Upsilon)|~~Z=0,~~r=0,~~\Upsilon=\pm\sqrt{2U(z)}, ~~z ~~is ~~a~~ critical~~ point~~of~~ U~~ in~~ \mathbf{B}^{2N-4} \}.
\end{equation*}
According to Proposition \ref{isolatedcriticalpoint1}, the set $CP$ consists of finitely isolated points.

By the relation of total energy
\begin{equation*}
    {\Upsilon}^2+ [\frac{ (\sum_{j=5}^{2N}{z}_j Z_j)^2}{z^2_3}+\sum_{j=5}^{2N}Z^2_j - (\sum_{j,k=5}^{2N}q_{jk}Z_j {z}_k)^2]- 2U(z)=0,
\end{equation*}
it follows that
\begin{equation*}
   \Upsilon'=  U(z)- \frac{1}{2}\Upsilon^2=\frac{ (\sum_{j=5}^{2N}{z}_j  {Z_j})^2}{z^2_3}+\sum_{j=5}^{2N}{Z_j}^2 - (\sum_{j,k=5}^{2N}q_{jk}{Z_j} {z}_k)^2.
\end{equation*}
Note that
\begin{equation*}
    \begin{array}{l}
      \frac{ (\sum_{j=5}^{2N}{z}_j  {Z_j})^2}{z^2_3}+\sum_{j=5}^{2N}{Z_j}^2 - (\sum_{j,k=5}^{2N}q_{jk}{Z_j} {z}_k)^2 \\
      \geq \frac{ (\sum_{j=5}^{2N}{z}_j  {Z_j})^2}{z^2_3}+\sum_{j=5}^{2N}{Z_j}^2 - \sum_{j=5}^{2N}{Z_j}^2 \sum_{j=5}^{2N}(\sum_{k=5}^{2N}q_{jk} {z}_k)^2\\
      = \frac{ (\sum_{j=5}^{2N}{z}_j  {Z_j})^2}{z^2_3}+\sum_{j=5}^{2N}{Z_j}^2 - \sum_{j=5}^{2N}{Z_j}^2 \sum_{j=5}^{2N}{z}_j^2 \\
      = \frac{ (\sum_{j=5}^{2N}{z}_j  {Z_j})^2}{z^2_3}+\sum_{j=5}^{2N}{Z_j}^2 z^2_3\geq 0,
    \end{array}
\end{equation*}
and  equality above hold if and only if $Z=0$, here we used the Cauchy-Schwarz inequality and the fact that $Q=(q_{jk})$ is an orthogonal matrix to deduce the inequality above.

Therefore $\Upsilon'>0$ holds unless
\begin{equation}
   Z=0,  ~~~~~U(z)= \frac{1}{2}\Upsilon^2.\nonumber
\end{equation}

By
\begin{equation*}
   \Upsilon'=  U(z)- \frac{1}{2}\Upsilon^2,
\end{equation*}
it follows that
\begin{equation*}
   \Upsilon''= \sum_{i=5}^{2N} \frac{\partial U(z)}{\partial z_i}z_i'- \Upsilon \Upsilon'
\end{equation*}
and
\begin{equation*}
   \Upsilon'''= \sum_{i=5}^{2N} \frac{\partial U(z)}{\partial z_i}z_i''+\sum_{i,j=5}^{2N} \frac{\partial^2 U(z)}{\partial z_i\partial z_j}z_i'z_j'- \Upsilon \Upsilon''-\Upsilon'^2,
\end{equation*}
where
\begin{equation}\label{zZd}
    \begin{array}{c}
       z'_i  =  Z_i,   \\
             Z'_i  =  \frac{\partial U(z)}{\partial z_i} - \frac{z_i \sum_{j=5}^{2N}({Z_j'} z_j + {Z_j}^2)}{z^2_3} - \frac{{z}_i (\sum_{j=5}^{2N}{z}_j  {Z_j})^2}{z^4_3}-  {\sum_{j,k=5}^{2N}q_{kj}{Z_j}' {z}_k} \sum_{j=5}^{2N}q_{ij} z_j\\
- 2 \sum_{j,k=5}^{2N}q_{kj}{Z_j} {z}_k \sum_{j=5}^{2N}q_{ij} {Z_j}-\frac{\Upsilon}{2} [ {Z_i} +\frac{z_i \sum_{j=5}^{2N}{Z_j} z_j}{z^2_3} +  \sum_{j,k=5}^{2N}q_{kj}{Z_j} {z}_k \sum_{j=5}^{2N}q_{ij} z_j]
    \end{array}
\end{equation}
for $i=5,\cdots,2N.$

As a result, when $\Upsilon'=0$, we have
\begin{equation*}
   \Upsilon''= 0
\end{equation*}
and
\begin{equation*}
    \Upsilon'''= \sum_{i=5}^{2N} \frac{\partial U(z)}{\partial z_i}Z_i'.
\end{equation*}
By (\ref{zZd}), it follows that
\begin{equation*}
   Z'_i  =  \frac{\partial U(z)}{\partial z_i} - \frac{z_i \sum_{j=5}^{2N}{Z_j'} z_j }{z^2_3} -  {\sum_{j,k=5}^{2N}q_{kj}{Z_j}' {z}_k} \sum_{j=5}^{2N}q_{ij} z_j.
\end{equation*}
Consequently,
\begin{equation*}
   \begin{array}{c}
     \sum_{i=5}^{2N}  Z'^2_i  =  \sum_{i=5}^{2N} \frac{\partial U(z)}{\partial z_i}Z'_i  -\sum_{i=5}^{2N} \frac{Z'_i z_i \sum_{j=5}^{2N}{Z_j'} z_j }{z^2_3} -  \sum_{i=5}^{2N} Z'_i{\sum_{j,k=5}^{2N}q_{kj}{Z_j}' {z}_k} \sum_{j=5}^{2N}q_{ij} z_j \\\\
     =  \Upsilon'''  - \frac{ (\sum_{j=5}^{2N}{Z_j'} z_j)^2 }{z^2_3} +  ({\sum_{j,k=5}^{2N}q_{jk}{Z_j}' {z}_k} )^2,
   \end{array}
\end{equation*}
here we used the fact that $Q=(q_{jk})$ is an  anti-symmetric matrix.

So
\begin{equation*}
    \Upsilon'''=\sum_{i=5}^{2N}  Z'^2_i  + \frac{ (\sum_{j=5}^{2N}{Z_j'} z_j)^2 }{z^2_3} -  ({\sum_{j,k=5}^{2N}q_{jk}{Z_j}' {z}_k} )^2.
\end{equation*}
An argument similar to the one used above shows that $\Upsilon'''>0$ holds unless
\begin{equation}
   Z'=0.\nonumber
\end{equation}
This yields that
\begin{equation*}
    \frac{\partial U(z)}{\partial z_i}=0, ~~~\forall i\in\{5,\cdots,2N\}
\end{equation*}

To summarize what we have proved, $\Upsilon$ is strictly monotone increasing with respect
to $\tau$, unless
  $(z, Z, 0, \Upsilon)$ is a point of $CP$.

We have thus proved that  the flow of the system (\ref{equationzr3}) restricted to the invariant manifold
 $\{r=0\}\bigcap\mathcal{N}_{\mathcal{H}}$
 is gradient-like with respect
to $\pi_{\Upsilon}$.

$~~~~~~~~~~~~~~~~~~~~~~~~~~~~~~~~~~~~~~~~~~~~~~~~~~~~~~~~~~~~~~~~~~~~~~~~~~~~~~~~~~~~~~~~~~~~~~~~~~~~~~~~~~~~~~~~~~~~~~~~~~~~~~~~~~~~~~~\Box$\\

\begin{remark}
By Proposition \ref{isolatedcriticalpoint1} and the following lemma
\begin{lemma} \emph{( \cite{Mcgehee1974Triple})}
Let $\phi$  be a flow on a locally compact metric space $\mathcal{M}$. Let
$p\in \mathcal{M}$ be such that $\omega(p)$ is a non-empty compact set. Suppose $\phi$ restricted
to $\omega(p)$ is gradient-like. Then $\omega(p)$ is a single point.
\end{lemma}
it is easy to see that a solution $( z(\tau), Z(\tau),r(\tau), \Upsilon(\tau), \theta(\tau))$ of the system (\ref{equationzr3}) corresponding to a total collision orbit $\mathbf{r}(t)$ of equations (\ref{eq:Newton's equation1})  asymptotically approaches an isolated point  in $CP$.  Consequently, the following asymptotic estimations on the total collision orbit $\mathbf{r}(t)$ hold:
\begin{itemize}
  \item $I(\mathbf{r}(t)) \sim (\frac{3}{2})^{\frac{4}{3}} \kappa^{\frac{2}{3}} t^{\frac{4}{3}}, \dot{I}(\mathbf{r}(t)) \sim (12)^{\frac{1}{3}} \kappa^{\frac{2}{3}} t^{\frac{1}{3}}, \ddot{I}(\mathbf{r}(t)) \sim (\frac{2}{3})^{\frac{2}{3}} \kappa^{\frac{2}{3}} t^{-\frac{2}{3}}$ as $t \rightarrow 0+$.
  \item $\mathcal{U}(\mathbf{r}(t))\sim  (\frac{1}{18})^{\frac{1}{3}} \kappa^{\frac{2}{3}} t^{-\frac{2}{3}}, \mathcal{K}(\dot{\mathbf{r}}(t))\sim (\frac{1}{18})^{\frac{1}{3}} \kappa^{\frac{2}{3}} t^{-\frac{2}{3}}$ as $t \rightarrow 0+$.
  \item $\hat{\mathbf{r}}(t) \rightarrow \textbf{CC}_{\lambda}$ as $t \rightarrow 0+$.
\end{itemize}
So the moving coordinates  provide  a setting to retrieve asymptotic estimation on total collision orbits.
\end{remark}

\subsection{Diagonalization of the Linear Part}\label{DiagonalizationoftheLinearPart}
\indent\par

The aim of the subsection is to diagonalize the linear part of equations (\ref{equationz4}) in detail.

Without loss of generality, suppose
\begin{displaymath}
\mu_j =0 ~~~~for ~~~~j\in \{5, \cdots, n_0+4 \},
\end{displaymath}
\begin{displaymath}
\mu_j >0 ~~~~for ~~~~j\in \{n_0+5, \cdots, n_0+n_p+4 \},
\end{displaymath}
\begin{displaymath}
0> \mu_j > - \frac{\kappa}{16} ~~~~for ~~~~j\in \{n_0+n_p+5, \cdots, n_0+n_p+n_1+4 \},
\end{displaymath}
\begin{displaymath}
\mu_j < - \frac{\kappa}{16} ~~~~for ~~~~j\in \{n_0+n_p+n_1+5, \cdots, n_0+n_p+n_1+n_2+4 \},
\end{displaymath}
\begin{displaymath}
\mu_j = - \frac{\kappa}{16} ~~~~for ~~~~j\in \{n_0+n_p+n_1+n_2+5, \cdots, 2N \}.
\end{displaymath}

Set
\begin{displaymath}
\Lambda_0 = diag( \mu_5, \cdots,  \mu_{n_0+4}),
\end{displaymath}
\begin{displaymath}
\Lambda_p = diag( \mu_{n_0+5}, \cdots, \mu_{n_0+n_p+4}),
\end{displaymath}
\begin{displaymath}
\Lambda_1 = diag( \mu_{n_0+n_p+5}, \cdots,  \mu_{n_0+n_p+n_1+4}),
\end{displaymath}
\begin{displaymath}
\Lambda_2 = diag( \mu_{n_0+n_p+n_1+5}, \cdots, \mu_{n_0+n_p+n_1+n_2+4}),
\end{displaymath}
\begin{displaymath}
\Lambda_3 = diag( \mu_{n_0+n_p+n_1+n_2+5}, \cdots,  \mu_{2N}),
\end{displaymath}
\begin{displaymath}
\tilde{\Lambda}_0 = diag( \tilde{\mu}_5, \cdots,  \tilde{\mu}_{n_0+4}),
\end{displaymath}
\begin{displaymath}
\tilde{\Lambda}_p = diag( \tilde{\mu}_{n_0+5}, \cdots, \tilde{\mu}_{n_0+n_p+4}),
\end{displaymath}
\begin{displaymath}
\tilde{\Lambda}_1 = diag( \tilde{\mu}_{n_0+n_p+5}, \cdots,  \tilde{\mu}_{n_0+n_p+n_1+4}),
\end{displaymath}
\begin{displaymath}
\tilde{\Lambda}_2 = diag( \tilde{\mu}_{n_0+n_p+n_1+5}, \cdots, \tilde{\mu}_{n_0+n_p+n_1+n_2+4}),
\end{displaymath}
\begin{displaymath}
\tilde{\Lambda}_3 = diag( \tilde{\mu}_{n_0+n_p+n_1+n_2+5}, \cdots,  \tilde{\mu}_{2N}),
\end{displaymath}
where
\begin{displaymath}
\tilde{\mu}_j = -\frac{\kappa^{\frac{1}{2}}}{4}+\sqrt{ \mu_j + \frac{\kappa}{16}} ~~~~for ~~~~j\in \{5, \cdots, 2N \}.
\end{displaymath}

Note that
$\tilde{\Lambda}_2$ is consisting of conjugate complex pair numbers. Without loss of generality, suppose that
\begin{displaymath}
\tilde{\mu}_{j} = \overline{\tilde{\mu}_{j+1}} ~~~~for ~~~~\lfloor \frac{ j+1}{2} \rfloor \in \{\lfloor \frac{ n_0+n_p+n_1+6}{2}\rfloor +k: k=1,  \cdots, \frac{n_2}{2}\}
\end{displaymath}
where $\overline{{\tilde{\mu}}}$ denotes the complex conjugation of a complex number $\tilde{\mu}$ and $\lfloor \cdot \rfloor$ denotes  the greatest integer function.

Set
\begin{equation}
\mathfrak{P}=\left(
    \begin{array}{cccc}
      \mathbb{I} &   & \mathbb{I} &   \\
        & \mathbb{I} &   & \mathbb{I} \\
      diag(\tilde{\Lambda}_0,\tilde{\Lambda}_p, \tilde{\Lambda}_1,\tilde{\Lambda}_2) &   & - \frac{\kappa^{\frac{1}{2}}}{2}\mathbb{I} - diag(\tilde{\Lambda}_0,\tilde{\Lambda}_p, \tilde{\Lambda}_1,\tilde{\Lambda}_2) &   \\
        & \tilde{\Lambda}_3 &   & \epsilon \mathbb{I} +\tilde{\Lambda}_3 \\
    \end{array}
  \right),\nonumber
\end{equation}
then one can demonstrate that $\mathfrak{P}$ is invertible for any $\epsilon >0$ and
\begin{displaymath}
\mathfrak{P}^{-1} =\left(
    \begin{array}{cccc}
      \frac{\mathbb{I}}{2}- \frac{\kappa^{\frac{1}{2}}}{8}D_{0} &   & -\frac{1}{2}D_{0}  &   \\
        & \mathbb{I} + \frac{\tilde{\Lambda}_3}{\epsilon} &   & -\frac{\mathbb{I}}{\epsilon} \\
      \frac{\mathbb{I}}{2} + \frac{\kappa^{\frac{1}{2}}}{8}D_{0}  &   &\frac{1}{2} D_{0} &   \\
         &-\frac{\tilde{\Lambda}_3}{\epsilon}    &  &\frac{\mathbb{I}}{\epsilon} \\
    \end{array}
  \right),
\end{displaymath}
\begin{displaymath}
\mathfrak{P}^{-1} \left(
           \begin{array}{cc}
              0 & \mathbb{I}    \\
             \Lambda & - \frac{\kappa^{\frac{1}{2}}}{2}\mathbb{I}   \\
           \end{array}
         \right) \mathfrak{P}=\left(
    \begin{array}{cccc}
      diag(\tilde{\Lambda}_0,\tilde{\Lambda}_p, \tilde{\Lambda}_1,\tilde{\Lambda}_2) &   &   &   \\
        & \tilde{\Lambda}_3 &   & \epsilon \mathbb{I} \\
       &   & - \frac{\kappa^{\frac{1}{2}}}{2}\mathbb{I} - diag(\tilde{\Lambda}_0,\tilde{\Lambda}_p, \tilde{\Lambda}_1,\tilde{\Lambda}_2) &   \\
         &    &  &\tilde{\Lambda}_3 \\
    \end{array}
  \right),
\end{displaymath}
where $D_{0}=diag(\frac{1}{ \sqrt{ \mu_5 + \frac{\kappa}{16}}},\cdots,\frac{1}{ \sqrt{ \mu_{2N-n_3} + \frac{\kappa}{16}}})$.

As a result, after applying the linear substitution
\begin{displaymath}
\left(
  \begin{array}{c}
    z \\
    Z \\
  \end{array}
\right)
=\mathfrak{P} q,
\end{displaymath}
the equations (\ref{equationz4}) are transformed to the following equations
\begin{equation}\label{equationz5}
\left\{
             \begin{array}{lr}
             q'_k=  - \varphi_k(q, \gamma),  & k\in \{5, \cdots, n_0+4\} \\
             q'_{2N+k}= - \frac{\kappa^{\frac{1}{2}}}{2} {q}_{2N+k} + \varphi_k(q, \gamma),  & k\in \{5, \cdots, n_0+4\}
             \end{array}
\right.
\end{equation}
\begin{equation}\label{equationz6}
\left\{
             \begin{array}{lr}
             q'_k= \tilde{{\mu}}_k {q}_k - \varphi_k(q, \gamma),  & k\in \{ n_0+5, \cdots, 2N-n_3\} \\
             q'_{2N+k}= (- \frac{\kappa^{\frac{1}{2}}}{2}-\tilde{\mu}_k) {q}_{2N+k} + \varphi_k(q, \gamma),  & k\in \{ n_0+5, \cdots, 2N-n_3\}
             \end{array}
\right.
\end{equation}
\begin{equation}\label{equationz7}
\left\{
             \begin{array}{lr}
             q'_k= - \frac{\kappa^{\frac{1}{2}}}{4} {q}_k + \epsilon {q}_{2N+k} - \varphi_k(q, \gamma),  & k\in \{2N-n_3 +1, \cdots, 2N\} \\
             q'_{2N+k}= - \frac{\kappa^{\frac{1}{2}}}{4} {q}_{2N+k} + \varphi_k(q, \gamma),  & k\in \{2N-n_3 +1, \cdots, 2N\}
             \end{array}
\right.
\end{equation}
\begin{equation}\label{equationz8}
\gamma'  = \kappa^{\frac{1}{2}}\gamma - \varphi_0(q, \gamma),
\end{equation} and  (\ref{equationtheta4}) is transformed to (\ref{equationtheta5}).
Furthermore, note the following facts:
\begin{equation}\label{conjugation}
\left\{
             \begin{array}{lr}
             q_k \in \mathbb{R},  & k ~or~ k-2N\in \{ 5, \cdots,2N-n_3-n_2, 2N-n_3+1,\cdots, 2N\}\\
             q_{k}=\overline{q_{k+1}},  & \lfloor \frac{ k+1}{2} \rfloor~or~ \lfloor \frac{ k+1}{2} \rfloor-N \in \{\lfloor \frac{ n_0+n_p+n_1+6}{2}\rfloor +j: j=1,  \cdots, \frac{n_2}{2}\}\\
             \overline{\varphi_k}(\overline{q}, \overline{\gamma})= \varphi_k(q, \gamma),  & k ~or~ k-2N\in \{ 5, \cdots,2N-n_3-n_2, 2N-n_3+1,\cdots, 2N\} \\
             \overline{\varphi_k}(\overline{q}, \overline{\gamma})= \varphi_{k+1}(q, \gamma),  & \lfloor \frac{ k+1}{2} \rfloor~or~ \lfloor \frac{ k+1}{2} \rfloor-N \in \{\lfloor \frac{ n_0+n_p+n_1+6}{2}\rfloor +j: j=1,  \cdots, \frac{n_2}{2}\}\\
\overline{\varphi_0}(\overline{q}, \overline{\gamma})= \varphi_0(q, \gamma), &
             \end{array}
\right.
\end{equation}
where $\overline{\varphi_k}=\overline{\varphi_k}(q,\gamma)$ denotes the power-series obtained by replacing the coefficients in the power-series $\varphi_k(q, \gamma)$ by their complex conjugates.

Note that
\begin{equation*}
 \begin{array}{l}
   \mathfrak{C}=\begin{pmatrix}\mathfrak{P}^{-1}&\\&1\end{pmatrix} \mathfrak{A} \begin{pmatrix}\mathfrak{P}&\\&1\end{pmatrix} \\
    =\begin{pmatrix}diag(\tilde{\Lambda}_0,\tilde{\Lambda}_p, \tilde{\Lambda}_1,\tilde{\Lambda}_2) &   &   & &  \\
        & \tilde{\Lambda}_3 &   & \epsilon \mathbb{I} &\\
       &   & - \frac{\kappa^{\frac{1}{2}}}{2}\mathbb{I} - diag(\tilde{\Lambda}_0,\tilde{\Lambda}_p, \tilde{\Lambda}_1,\tilde{\Lambda}_2) &   &\\
         &    &  &- \frac{\kappa^{\frac{1}{2}}}{2}\mathbb{I}-\tilde{\Lambda}_3 &\\
       &&&  & \kappa^{\frac{1}{2}}\end{pmatrix}\\
       =\begin{pmatrix}diag(0,\tilde{\Lambda}_p, \tilde{\Lambda}_1,\tilde{\Lambda}_2) &   &   & &  \\
        & -\frac{\kappa^{\frac{1}{2}}}{4}\mathbb{I} &   & \epsilon \mathbb{I} &\\
       &   & - \frac{\kappa^{\frac{1}{2}}}{2}\mathbb{I} - diag(0,\tilde{\Lambda}_p, \tilde{\Lambda}_1,\tilde{\Lambda}_2) &   &\\
         &    &  &-\frac{\kappa^{\frac{1}{2}}}{4}\mathbb{I} &\\
       &&&  & \kappa^{\frac{1}{2}}\end{pmatrix}.
 \end{array}
\end{equation*}
\\

      \section{Normal Forms}\label{NormalFormsth}
      \indent\par

The key ideas of resolving $PISPW$ focus on estimating the rate of tending to zero in (\ref{Asymptotic condition}) to ensure that
$\theta(\tau)$ approaches a fixed limit as $\tau \rightarrow -\infty$.
For this purpose, we need further simplify equations (\ref{equationz4}) by the theory of normal forms (or reduction theorem).

Note that equations in (\ref{equationz5}) are degenerate (may be too degenerate) corresponding to the degeneration of central configuration. In general, degenerate equations are very difficult to handle.

Based on some original ideas of Siegel  \cite{C1967Lectures}, we will naturally absorb more general  results of normal forms than that of in \cite{C1967Lectures} to explore the problem contact with  degenerate equations.

As a matter of notational convenience, in the following we will abuse some notations which have slightly different  meaning from the previous content.

Suppose the origin $0$ is an equilibrium point of (\ref{differential system}), then it follows from Taylor expansion near the origin that
\begin{equation}
v(q) = \frac{\partial v(0)}{\partial q}q + o(q),\nonumber
\end{equation}
the system (\ref{differential system}) becomes
\begin{equation}
\dot{q} = \mathfrak{C}q+h(q),\nonumber
\end{equation}
where $\mathfrak{C}=\frac{\partial v(0)}{\partial q}$, and $h(q)=o(q)$, i.e., $h(0)=0, \frac{\partial h(0)}{\partial q}=0$.

If the function $v$ is $C^m$-smooth,  then the expansion
\begin{equation}
h(q) = h_2(q) + \cdots +h_m(q) + o_m(q)\nonumber
\end{equation}
is valid, where $h_k(q)$ is a homogeneous polynomial of the power
$k$; hereafter $o_m(q)$ stands for the terms which vanish at the origin along with
the first $m$ derivatives. Furthermore, if
$v$ is analytic, the expansion
\begin{equation}
h(q) = h_2(q) + \cdots +h_m(q) + \cdots\nonumber
\end{equation}
is valid, which is power-series starting with quadratic terms.

Let $\mu_1, \cdots, \mu_n$ denote the eigenvalues of the matrix $\mathfrak{C}$. The set $\{\mu_1, \cdots, \mu_n\}$ of the eigenvalues $\mu_1, \cdots, \mu_n$  is called
a \textbf{resonant} set if there exists a linear relationship
\begin{equation}
\mu_k = (\alpha, \mu) = \alpha_1 \mu_1 + \cdots +\alpha_n \mu_n\nonumber
\end{equation}
where $\alpha=(\alpha_1,   \cdots,  \alpha_n) \in \mathbb{N}^n$ is a multiindex (i.e.,the row of non-negative integers) such that $|\alpha|=\alpha_1 + \cdots +\alpha_n \geq 2$. The relation itself is called a resonance and $|\alpha|$ is called the order
of the resonance. And any constant multiple of the  monomial $ q^\alpha e_k= q^{\alpha_1}_1\cdots q^{\alpha_n}_n e_k$ is called the resonant monomial, where $e_k = (\underbrace{0,\cdots,0,1}_k,0,\cdots,0)^\top$.

\begin{definition}
A collection $\mu = \{\mu_1, \cdots, \mu_n\}$ of the eigenvalues belongs to the Poincar\'{e} domain if the
convex hull of the $n$ points $\mu_1, \cdots, \mu_n$
 in the complex plane does not
contain zero. Otherwise, the collection $\mu = \{\mu_1, \cdots, \mu_n\}$  belongs to the Siegel
region.
\end{definition}
Geometrically, the condition for Poincar\'{e} domain means that there exists a line  in the complex plane which separate the eigenvalues $\mu_1, \cdots,
\mu_n$  from
the origin (i.e. the eigenvalues are on one and the same side of the line while 0 is on
the other side).  Furthermore, one can prove that
\begin{proposition}\label{Poincare}
A collection $\mu = \{\mu_1, \cdots, \mu_n\}$ belongs to the Poincar\'{e} domain if and only if there exists a positive number $c$ such that
\begin{equation}
|(\alpha, \mu)|\geq c |\alpha|\nonumber
\end{equation}
for any  multiindex  $\alpha=(\alpha_1,   \cdots,  \alpha_n)$, $|\alpha| \geq 1$.
\end{proposition}
{\bf Proof.} 
We can assume that $\mu_k \neq 0$ for
any $k \in \{1,\cdots,n\}$.

If the collection $\mu = \{\mu_1, \cdots, \mu_n\}$ belongs to the Poincar\'{e} domain, after proper rotation of the complex plane, we can suppose that there exists a $\delta >0$ such that $Re \mu_k \geq \delta$ for
any $k \in \{1,\cdots,n\}$. Thus
\begin{equation}
|(\alpha, \mu)|\geq |\alpha_1 Re \mu_1 + \cdots +\alpha_n Re \mu_n| \geq \delta |\alpha|.\nonumber
\end{equation}

If $|(\alpha, \mu)|\geq c |\alpha|$ for some $c>0$, after proper rotation of the complex plane, we can suppose that
\begin{equation}
\alpha_1 Im \mu_1 + \cdots +\alpha_n Im \mu_n=0.\nonumber
\end{equation}
It follows that
\begin{equation}
|(\alpha, \mu)| = |\alpha_1 Re \mu_1 + \cdots +\alpha_n Re \mu_n| \geq c |\alpha|.\nonumber
\end{equation}
We claim that all of the real parts $Re \mu_1, \cdots , Re \mu_n$ are simultaneously greater or less than 0. Otherwise, by using reduction to absurdity, assume
\begin{description}
  \item[i.] $Re \mu_1=0$;
  \item[ii.] $Re \mu_1 >0 >Re \mu_2$.
\end{description}

First, it can easily be checked  that the case \textbf{i} is impossible.

Next, for the case \textbf{ii}, let
 \begin{equation}
\begin{array}{c}
  \alpha_3 = \cdots =\alpha_n=0,\nonumber \\
  \alpha_1>0,  ~~~\alpha_2>0,\nonumber
\end{array}
\end{equation}
then
\begin{equation}
 |\alpha_1 Re \mu_1 + \alpha_2 Re \mu_2| \geq c (\alpha_1  + \alpha_2)\nonumber
\end{equation}
or
\begin{equation}
 |\frac{\alpha_1}{\alpha_2}  +  \frac{Re \mu_2}{Re \mu_1}| \geq \frac{c (\alpha_1  + \alpha_2)}{\alpha_2 Re \mu_1}\geq \frac{c }{Re \mu_1}.\nonumber
\end{equation}
However, this contradicts with $|\frac{\alpha_1}{\alpha_2}  +  \frac{Re \mu_2}{Re \mu_1}|$ can be arbitrarily small by appropriately selecting positive integers $\alpha_1,\alpha_2$.

$~~~~~~~~~~~~~~~~~~~~~~~~~~~~~~~~~~~~~~~~~~~~~~~~~~~~~~~~~~~~~~~~~~~~~~~~~~~~~~~~~~~~~~~~~~~~~~~~~~~~~~~~~~~~~~~~~~~~~~~~~~~~~~~~~~~~~~~~~~~~~~~~~~~~~~~~~~~~~~~~~~~\Box$
Thus
each collection in the Poincar\'{e} region satisfies at most a finite number of resonances.

\begin{definition}
A collection $\mu = \{\mu_1, \cdots, \mu_n\}$ of the eigenvalues is said to be of type $(c,\upsilon)$ if there exists a positive number $c$ such that, for
any $k \in \{1,\cdots,n\}$ and all  multiindexes  $\alpha=(\alpha_1,   \cdots,  \alpha_n)$, $|\alpha| \geq 2$, we have
\begin{equation}
|(\alpha, \mu)-\mu_k|\geq \frac{c}{|\alpha|^\upsilon}. \nonumber
\end{equation}
\end{definition}

\subsection{Analytical Case}
\indent\par

Consider a  system
\begin{equation}\label{differential system1}
q' = \mathfrak{C}q + \varphi(q),
\end{equation}
where $q=(q_1,\cdots,q_{n})^\top$,
$\mathfrak{C}$ denotes the constant square matrix of order $n$,
$\varphi$ is a column vector whose components $\varphi_k$ are power-series in
the $n$ independent variables $q_1,\cdots,q_{n}$ with real coefficients and starting with quadratic terms.

Suppose $\varphi$ is convergent for sufficiently small $|q|$, and $\mu_1, \cdots, \mu_n$ denote the eigenvalues of the matrix $\mathfrak{C}$. Then for the system (\ref{differential system1}),  it is well known that:
\begin{theorem}\label{normalformps}\emph{(\textbf{Poincar\'{e}-Siegel})}
If  the collection $\{\mu_1, \cdots, \mu_n\}$ of the eigenvalues belongs to the Poincar\'{e} domain and is nonresonant,  or is  of type $(c,\upsilon)$, then one can
introduce a nonlinear substitution of the form
\begin{equation}\label{nonlinear substitutionps}
u_k = q_k - F_k(q_1, \cdots,q_{n}),  ~~~ k \in \{1,\cdots,n\},\nonumber
\end{equation}
so that the system (\ref{differential system1}) can be reduced to the following simple form
\begin{equation}\label{simple formsps}
u' = \mathfrak{C} u,\nonumber
\end{equation}
where  the $F_k$ are power-series starting with quadratic terms and convergent
for small $|q|$.
\end{theorem}
\begin{theorem}\label{normalformpd}\emph{(\textbf{Poincar\'{e}-Dulac})}
If  the collection $\{\mu_1, \cdots, \mu_n\}$ of the eigenvalues belongs to the Poincar\'{e} domain, then one can
introduce a nonlinear substitution of the form
\begin{equation}\label{nonlinear substitutionpd}
u_k =  q_k - F_k(q_1, \cdots,q_{n}),  ~~~ k \in \{1,\cdots,n\},\nonumber
\end{equation}
 so that the system (\ref{differential system1}) can be reduced to the following simple form
\begin{equation}\label{simple formspd}
u' = \mathfrak{C} u+ R(u) ,\nonumber
\end{equation}
where the $F_k$ are power-series starting with quadratic terms and convergent
for small $|q|$, and
the $R(u)$ is a finite-order polynomial composed by resonant monomials.
\end{theorem}

Even for the hyperbolic equilibrium point, the collection of the eigenvalues is frequently resonant, thus the above classical results can not be directly utilized.
It is natural to consider other available forms. Suppose  the matrix $\mathfrak{C}$ in the system (\ref{differential system1}) is block-diagonal form $\mathfrak{C}=\left(
                                                                                                                                                \begin{array}{cc}
                                                                                                                                                  \mathfrak{C}^+ &   \\
                                                                                                                                                    & \mathfrak{C}^- \\
                                                                                                                                                \end{array}
                                                                                                                                              \right)
$, and $\mu_1, \cdots, \mu_m$ denote the eigenvalues of the matrix $\mathfrak{C}^+$, $\mu_{m+1}, \cdots, \mu_n$ denote the eigenvalues of the matrix $\mathfrak{C}^-$. Set
\begin{displaymath}
\begin{array}{c}
  q^+ = (q_1, \cdots, q_m)^\top \\
  q^- = (q_{m+1}, \cdots, q_n)^\top
\end{array}
\end{displaymath}
Essentially Siegel gave the following simpler form  of the system (\ref{differential system1}) in \cite{C1967Lectures}.

\begin{theorem}\label{normalform1}\emph{(\textbf{Siegel})}
If  the collection $\{\mu_1, \cdots, \mu_m\}$ of the eigenvalues belongs to the Poincar\'{e} domain, then we can
introduce a nonlinear substitution of the form
\begin{equation}\label{nonlinear substitution1}
\left\{
             \begin{array}{lr}
             u^+ =  q^+ - F^+(q_1, \cdots,q_{m}), \\
             u^- =  q^- - F^-(q_1, \cdots,q_{m}),
             \end{array}
\right.\nonumber
\end{equation}
where  the vector-valued functions $F^+, F^-$ are power-series in the $m$ independent
 variables $q_1, \cdots,q_{m}$ only, starting with quadratic terms and convergent
for small $|q_1|, \cdots,|q_{m}|$. So that we can write the system (\ref{differential system1}) in the simpler form
\begin{equation}\label{simple forms1}
\left\{
             \begin{array}{lr}
             {u'^+} = \mathfrak{C}^+ u^+ + R^+(u) +  \psi^+(u), \\
             {u'^-} = \mathfrak{C}^- u^-  +  \psi^-(u),
             \end{array}
\right.\nonumber
\end{equation}
where
the vector-valued function $R^+(u)$ is a finite-order polynomial composed by resonant monomials with regard to $\{\mu_1, \cdots, \mu_m\}$, and the vector-valued functions $\psi^+, \psi^-$ are power-series in the $n$ independent
 variables $u_1, \cdots,u_{n}$ starting with quadratic terms and convergent
for small $|u_1|, \cdots,|u_{n}|$, furthermore, for $u_{m+1} = \cdots = u_{n} =0$,
\begin{equation}\label{psi}
\psi^+(u_1, \cdots,u_{n}) = \psi^-(u_1, \cdots,u_{n}) \equiv 0\nonumber
\end{equation}
\end{theorem}

Based upon the ideas and methods of Siegel in proving the above form  in \cite{C1967Lectures},  we prove
\begin{theorem}\label{normalform2}\emph{(\textbf{Siegel})}
If  the collection $\{\mu_1, \cdots, \mu_m\}$ of the eigenvalues belongs to the Poincar\'{e} domain, then we can
introduce a nonlinear substitution of the form
\begin{equation}\label{nonlinear substitution2}
\left\{
             \begin{array}{lr}
             u^+ =  q^+ , \\
             u^- =  q^- - F^-(q_1, \cdots,q_{m}),
             \end{array}
\right.\nonumber
\end{equation}
where  the vector-valued function $F^-$ are power-series in the $m$ independent
 variables $q_1, \cdots,q_{m}$ only, starting with quadratic terms and convergent
for small $|q_1|, \cdots,|q_{m}|$. So that we can write the system (\ref{differential system1}) in the simpler form
\begin{equation}\label{simple forms2}
\left\{
             \begin{array}{lr}
             {u'^+} = \mathfrak{C}^+ u^+  +  \psi^+(u), \\
             {u'^-} = \mathfrak{C}^- u^-  +  \psi^-(u),
             \end{array}
\right.\nonumber
\end{equation}
where
the vector-valued functions $\psi^+, \psi^-$ are power-series in the $n$ independent
 variables $u_1, \cdots,u_{n}$ starting with quadratic terms and convergent
for small $|u_1|, \cdots,|u_{n}|$, furthermore, for $u_{m+1} = \cdots = u_{n} =0$,
\begin{equation}\label{psi}
\psi^-(u_1, \cdots,u_{n}) \equiv 0  .\nonumber
\end{equation}
\end{theorem}

Since the proof of Theorem \ref{normalform2} is long and quite similar to that given  by Siegel in \cite{C1967Lectures}, no proof will be given here.

\begin{remark}
Note the difference between \emph{Theorem \ref{normalform1}} and \emph{Theorem \ref{normalform2}}.

In particular, note that the following  two statements  are equivalent:
\begin{enumerate}
  \item For $u_{m+1} = \cdots = u_{n} =0$,
\begin{equation}
\psi^+(u_1, \cdots,u_{n}) = \psi^-(u_1, \cdots,u_{n}) \equiv 0\nonumber
\end{equation}
where
the vector-valued functions $\psi^+, \psi^-$ are power-series in the $n$ independent
 variables $u_1, \cdots,u_{n}$ starting with quadratic terms and convergent
for small $|u_1|, \cdots,|u_{n}|$.
  \item The vector-valued functions $\psi^+, \psi^-$ have the  following forms:
\begin{equation}
\left\{
             \begin{array}{lr}
             \psi^+(u) =  \tilde{\psi}^+(u)u^-, \\
             \psi^-(u) =   \tilde{\psi}^-(u)u^-,
             \end{array}
\right.\nonumber
\end{equation}
where the matrix-valued functions $\tilde{\psi}^+, \tilde{\psi}^-$ are power-series in the $n$ independent
 variables $u_1, \cdots,u_{n}$ starting with linear terms and convergent
for small $|u_1|, \cdots,|u_{n}|$.
\end{enumerate}

\end{remark}
\begin{remark}
It is noteworthy that the estimate $|\alpha_1\mu_1 + \cdots +\alpha_m \mu_m| \geq c |\alpha_1 + \cdots + \alpha_m|$ is necessary for proving the convergence of the nonlinear substitution in \emph{Theorem \ref{normalform1}} and \emph{Theorem \ref{normalform2}}. That is, according to \emph{Proposition \ref{Poincare}}, the condition of Poincar\'{e} is necessary for \emph{Theorem \ref{normalform1}} and \emph{Theorem \ref{normalform2}}. It is essentially different from the estimation of type $(c,\upsilon)$ in \emph{Theorem \ref{normalformps}}, since  the convergence of the nonlinear substitution in  \emph{Theorem \ref{normalform1}} and \emph{Theorem \ref{normalform2}} with the condition of type $(c,\upsilon)$ can not be guaranteed by the well-known scheme of KAM.

Indeed, even the estimation $|\alpha_1\mu_1 + \cdots +\alpha_m \mu_m| \geq c$, which is stronger than that of type $(c,\upsilon)$, is not sufficient for proving the convergence of the nonlinear substitution in \emph{Theorem \ref{normalform1}} or \emph{Theorem \ref{normalform2}:}

Let us consider the two-dimensional system
\begin{equation}
\left\{
             \begin{array}{lr}
             p'  = p^2 &  \\
             q'  = p+ q &
             \end{array}
\right.\nonumber
\end{equation}
There is a unique formal diffeomorphism $p=u,q=v+F(u)$ that transforms the
previous systems into its normal form
\begin{equation}
\left\{
             \begin{array}{lr}
             u'  = u^2 &  \\
             v'  = v &
             \end{array}
\right.\nonumber
\end{equation}
It is easy to see that the formal function $F$ is
\begin{equation}
F(u)= - \sum_{k\geq 1}(k-1)! u^k.\nonumber
\end{equation}
This does not converge in a neighborhood of the origin, although we have the  estimate $|\alpha_1 1| \geq 1$.
\end{remark}

It is easy to show that it follows from Theorem \ref{normalform2} that the following reduction theorem is  true.
\begin{corollary}\label{normalform3}\emph{(\textbf{Reduction Theorem})}
If the collections $\{\mu_1, \cdots, \mu_m\}$ and $\{\mu_{m+1}, \cdots, \mu_n\}$ of the eigenvalues belong to the Poincar\'{e} domain, then we can
introduce a nonlinear substitution of the form
\begin{equation}\label{nonlinear substitution3}
\left\{
             \begin{array}{lr}
             u^+ =  q^+  - F^+ (q_{m+1}, \cdots,q_{n}),  & \\
             u^-  =  q^-  - F^- (q_1, \cdots,q_{m}),  &
             \end{array}
\right.\nonumber
\end{equation}
where  the vector-valued functions $F^+$ and $F^-$ are respectively power-series starting with quadratic terms and convergent
for small $|q_{m+1}|, \cdots,|q_{n}|$ and $|q_1|, \cdots,|q_{m}|$. So that we can write the system (\ref{differential system1}) in the simpler form
\begin{equation}\label{simple forms3}
\left\{
             \begin{array}{lr}
             {u'^+} = \mathfrak{C}^+ u^+  +  \psi^+(u)u^+, \\
             {u'^-} = \mathfrak{C}^- u^-  +  \psi^-(u)u^-,
             \end{array}
\right.
\end{equation}
where
the matrix-valued functions $\psi^+, \psi^-$ are power-series in the $n$ independent
 variables $u_1, \cdots,u_{n}$ starting with linear terms and convergent
for small $|u_1|, \cdots,|u_{n}|$.
\end{corollary}

In particular, for a hyperbolic equilibrium
point $0$, the collection of  the eigenvalues of the matrix $\mathfrak{C}$ having positive real part  or  the eigenvalues of the matrix $\mathfrak{C}$ having negative real part
belongs to the Poincar\'{e} domain.
So it follows from the above Reduction Theorem that  the following celebrated result holds.
\begin{theorem}\label{analyticalinvariantmanifolds}\emph{(\textbf{Poincar\'{e}-
Lyapunov})}
Suppose the eigenvalues of the matrix $\mathfrak{C}^+$ have positive real part, the eigenvalues of the matrix $\mathfrak{C}^-$ have negative real part, that is,
the  equilibrium
point $0$ of (\ref{differential system1}) is hyperbolic. Then $0$  has  local analytic invariant manifolds $\mathcal{W}^s_{loc}$ and $\mathcal{W}^u_{loc}$ whose equations are

\begin{equation}\label{analytic invariant manifolds}
\begin{array}{c}
  \mathcal{W}^s_{loc}:  ~~~ q^+  =  F^+(q^-)\\
  \mathcal{W}^u_{loc}:  ~~~ q^-  =  F^-(q^+)
\end{array}\nonumber
\end{equation}
where  the  vector-valued functions $F^+$ and $F^-$ are respectively power-series starting with quadratic terms and convergent
for small $|q_{m+1}|, \cdots,|q_{n}|$ and $|q_1|, \cdots,|q_{m}|$.
\end{theorem}
\begin{remark}
Of course, similar to an argument in the following subsection, the existence of the local analytic invariant manifolds $\mathcal{W}^s_{loc}$ and $\mathcal{W}^u_{loc}$ can be used to transform the system (\ref{differential system1}) in the simpler form (\ref{simple forms3}).
\end{remark}

We can solve $PISPW$ corresponding to the nondegenerate central configuration by using the above Reduction Theorem or Theorem \ref{normalform1}.

\subsection{Smooth Case}
\indent\par

As before, consider the  system
\begin{equation}\label{differential system2}
q' = \mathfrak{C}q + \varphi(q).
\end{equation}

When one of the eigenvalues $\mu_1, \cdots, \mu_n$ of the matrix $\mathfrak{C}$ is on the imaginary axis, we have to consider smooth change of variables to simplify the system (\ref{differential system2}),   even though the right side of the system (\ref{differential system2}) is  analytic power-series in
the $n$ independent variables $q_1,\cdots,q_{n}$. Suppose $\varphi(q)$ is $C^l$-smooth ($1\leq l \leq \infty$) in
the $n$ independent variables $q_1,\cdots,q_{n}$ in the following.

Suppose  the matrix $\mathfrak{C}$ in the system (\ref{differential system2}) is block-diagonal form
\begin{displaymath}
\mathfrak{C}=\left(
                                                                                                                                                \begin{array}{ccc}
                                                                                                                                                  \mathfrak{C}^0 &  & \\
                                                                                                                                                    & \mathfrak{C}^+ & \\
                                                                                                                                                   &  & \mathfrak{C}^-\\
                                                                                                                                                \end{array}
                                                                                                                                              \right).
\end{displaymath}
Set
\begin{displaymath}
\begin{array}{c}
  q^0 = (q_1, \cdots, q_{n_0})^\top ,\\
  q^+ = (q_{n_0+1}, \cdots, q_{n_0+m})^\top, \\
  q^- = (q_{n_0+m+1}, \cdots, q_n)^\top.
\end{array}
\end{displaymath}
Suppose the eigenvalues of the matrix $\mathfrak{C}^0$ are all on the imaginary axis, the eigenvalues of the matrix $\mathfrak{C}^+$ lie
to the right of the imaginary axis, and the eigenvalues of the matrix $\mathfrak{C}^-$ lie
to the left of the imaginary axis.


Similar to Theorem \ref{analyticalinvariantmanifolds} of {Poincar\'{e}-Lyapunov}, it is well known that the following results hold (our main reference on this issue is \cite{Shilnikov}).
\begin{theorem}\label{Center-Stable Manifold}\emph{\textbf{(Center-Stable Manifold)}}
In a small neighborhood of the origin
there exists an $n-m$-dimensional invariant center-stable manifold $\mathcal{W}^{cs}_{loc}| q^+ = F^{cs}(q^0,q^-)$ of class $C^l$ ($l< \infty$), which contains the origin and which is tangent to the subspace
$q^+=0$ at the origin. The manifold $\mathcal{W}^{cs}_{loc}$ contains all orbits which stay in a small
neighborhood of the origin for all positive times. Though the center-stable manifold
is not defined uniquely, for any two manifolds $\mathcal{W}^{cs}_{1}$ and $\mathcal{W}^{cs}_{2}$
the functions
$F^{cs}_{1}$ and $F^{cs}_{2}$
 have the same Taylor expansion at  the origin (and at each point whose
 positive semiorbit stays in a small neighborhood of the origin).
\end{theorem}

\begin{remark}
Note that even if the system is $C^\infty$-smooth,
the center-stable manifold has, in general, only finite smoothness. Of course,
if the original system is $C^\infty$-smooth, it is $C^l$-smooth
for any finite $l$. Therefore, in this case
one may apply the center-stable manifold theorem with any given $l$ which implies that: for any finite
$l$ there exists a neighborhood $\mathcal{N}_l$ of $0$ where $\mathcal{W}^{cs}_{loc}$ is $C^l$-smooth.
In principle, however, these neighborhoods may shrink to zero as $l\rightarrow \infty$.
\end{remark}

\begin{theorem}\label{Center-Unstable Manifold}\emph{\textbf{(Center-Unstable Manifold)}}
In a small neighborhood of the origin
there exists an $n_0 + m$-dimensional invariant center-unstable manifold $\mathcal{W}^{cu}_{loc}| q^- = F^{cu}(q^0,q^+)$ of class $C^l$ ($l< \infty$), which contains the origin and which is tangent to the subspace
$q^-=0$ at the origin. The manifold $\mathcal{W}^{cu}_{loc}$ contains all orbits which stay in a small
neighborhood of the origin for all negative times. Though the center-unstable manifold
is not defined uniquely, for any two manifolds $\mathcal{W}^{cu}_{1}$ and $\mathcal{W}^{cu}_{2}$
the functions
$F^{cu}_{1}$ and $F^{cu}_{2}$
 have the same Taylor expansion at the origin (and at each point whose
 negative semiorbit stays in a small neighborhood of the origin). In
the case where the system is $C^\infty$-smooth, the center-unstable manifold has, in
general, only finite smoothness.
\end{theorem}

Note that the condition of invariance of the manifolds $\mathcal{W}^{cs}_{loc}$ and $\mathcal{W}^{cu}_{loc}$ may be
expressed as
\begin{equation}\label{condition of invariance}
\begin{array}{c}
  q'^+ = \frac{\partial F^{cs}(q^0,q^-)}{\partial q^0}q'^0 +\frac{\partial F^{cs}(q^0,q^-)}{\partial q^-}q'^-   ~~~ ~~when ~~q^+  =  F^{cs}(q^0,q^-),\\
  q'^- = \frac{\partial F^{cu}(q^0,q^+)}{\partial q^0}q'^0 +\frac{\partial F^{cu}(q^0,q^+)}{\partial q^+}q'^+   ~~~ ~~when ~~q^-  =  F^{cu}(q^0,q^+),
\end{array}\nonumber
\end{equation}
or
\begin{equation}\label{condition of invariance1}
\begin{array}{c}
  \mathfrak{C}^+ F^{cs}(q^0,q^-) + \varphi^+(q^0,F^{cs},q^-) = \\
\frac{\partial F^{cs}(q^0,q^-)}{\partial q^0}[\mathfrak{C}^0 q^0 + \varphi^0(q^0,F^{cs},q^-)] +\frac{\partial F^{cs}(q^0,q^-)}{\partial q^-}[\mathfrak{C}^- q^- + \varphi^-(q^0,F^{cs},q^-)],\\
   \mathfrak{C}^- F^{cu}(q^0,q^+) + \varphi^-(q^0,q^+, F^{cu}) = \\
\frac{\partial F^{cu}(q^0,q^+)}{\partial q^0}[\mathfrak{C}^0 q^0 + \varphi^0(q^0,q^+, F^{cu})] +\frac{\partial F^{cu}(q^0,q^+)}{\partial q^+}[\mathfrak{C}^+ q^+ + \varphi^+(q^0,q^+, F^{cu})].
\end{array}\nonumber
\end{equation}
The above relations  yield an algorithm for computing the invariant
manifolds which will be used in the following.

Let us introduce new variables:
\begin{displaymath}
\left\{
             \begin{array}{lr}
             u^0  =  q^0    &\\
             u^+ =  q^+  -  F^{cs}(q^0,q^-)  & \\
             u^-  =  q^-    &
             \end{array}
\right.
\end{displaymath}
Then  the system (\ref{differential system1}) becomes:
\begin{equation}\label{simple forms}
\left\{
             \begin{array}{lr}
             {u'^0} = \mathfrak{C}^0 u^0 + \varphi^0(u^0,  F^{cs}(u^0,u^-),u^-) + \psi^0(u)u^+, \\
             {u'^+} = \mathfrak{C}^+ u^+  + \psi^+(u)u^+, \\
             {u'^-} = \mathfrak{C}^- u^-  + \varphi^-(u^0,  F^{cs}(u^0,u^-),u^-)+ \psi^-(u)u^+
             \end{array}
\right.\nonumber
\end{equation}
As a matter of fact, from the new variables, it follows that
\begin{equation}\label{guodu}
\left\{
             \begin{array}{lr}
             {u'^0} = \mathfrak{C}^0 u^0 + \varphi^0(u^0,u^+ +  F^{cs}(u^0,u^-),u^-), \\
{u'^+} = q'^+-\frac{\partial F^{cs}(q^0,q^-)}{\partial q^0}q'^0 -\frac{\partial F^{cs}(q^0,q^-)}{\partial q^-}q'^-, \\
             {u'^-} = \mathfrak{C}^- u^- + \varphi^-(u^0,u^+ +  F^{cs}(u^0,u^-),u^-),
             \end{array}
\right.\nonumber
\end{equation}
By (\ref{condition of invariance1}), the second equation in (\ref{guodu}) may be rewritten as
\begin{eqnarray}
{u'^+} &=& \mathfrak{C}^+ u'^+ +  \mathfrak{C}^ + F^{cs} + \varphi^+(q^0,u^++F^{cs},q^-)-\frac{\partial F^{cs}}{\partial q^0}[\mathfrak{C}^0 q^0 + \varphi^0(q^0,u^++F^{cs},q^-)]\nonumber\\
&-&\frac{\partial F^{cs}}{\partial q^-}[\mathfrak{C}^- q^- + \varphi^-(q^0,u^++F^{cs},q^-)]\nonumber \\
&=& \mathfrak{C}^+ u'^+   + [\varphi^+(q^0,u^++F^{cs},q^-)-\varphi^+(q^0,+F^{cs},q^-)]\nonumber\\
&-&\frac{\partial F^{cs}}{\partial q^0}[ \varphi^0(q^0,u^++F^{cs},q^-)-\varphi^0(q^0,F^{cs},q^-)]\nonumber\\
&-&\frac{\partial F^{cs}}{\partial q^-}[ \varphi^-(q^0,u^++F^{cs},q^-)-\varphi^-(q^0,F^{cs},q^-)]\nonumber
\end{eqnarray}
Since
\begin{equation}\label{guodu1}
    \varphi(q^0,u^++F^{cs},q^-)-\varphi(q^0,F^{cs},q^-)= [\int^1_0 \frac{\partial\varphi(q^0,t u^++F^{cs},q^-)}{\partial q^+}dt] u^+,\nonumber
\end{equation}
we can write the system (\ref{differential system1}) in the simpler form
\begin{equation}\label{simple forms4}
\left\{
             \begin{array}{lr}
             {u'^0} = \mathfrak{C}^0 u^0 + \varphi^0(u^0,  F^{cs}(u^0,u^-),u^-) + \psi^0(u)u^+, \\
             {u'^+} = \mathfrak{C}^+ u^+  + \psi^+(u)u^+, \\
             {u'^-} = \mathfrak{C}^- u^-  + \varphi^-(u^0,  F^{cs}(u^0,u^-),u^-)+ \psi^-(u)u^+,
             \end{array}
\right.\nonumber
\end{equation}
where
\begin{displaymath}
\begin{array}{c}
  \psi^0(u)= \int^1_0 \frac{\partial\varphi^0(u^0,t u^++F^{cs},u^-)}{\partial q^+}dt, \\
  \psi^-(u)= \int^1_0 \frac{\partial\varphi^-(u^0,t u^++F^{cs},u^-)}{\partial q^+}dt
\end{array}
\end{displaymath}
are $C^{l}$-smooth, and
\begin{displaymath}
\begin{array}{c}
  \psi^+(u)= \int^1_0 \frac{\partial\varphi^+(u^0,t u^++F^{cs},u^-)}{\partial q^+}dt-\int^1_0\frac{\partial F^{cs}}{\partial u^0}  \frac{\partial\varphi^0(u^0,t u^++F^{cs},u^-)}{\partial q^+}dt-\frac{\partial F^{cs}}{\partial u^-}  \frac{\partial\varphi^-(u^0,t u^++F^{cs},u^-)}{\partial q^+}dt
\end{array}
\end{displaymath}
is $C^{l-1}$-smooth. Moreover, by virtue of $\varphi(0)=0$ and $ \frac{\partial\varphi(0)}{\partial q}=0$,  $\psi(0)=0$ holds.

Similarly, we have
\begin{corollary}\label{normalform5}\emph{(\textbf{Reduction Theorem})}
One can introduce new variables:
\begin{displaymath}
\left\{
             \begin{array}{lr}
             u^0  =  q^0    &\\
             u^+ =  q^+    & \\
             u^-  =  q^-  -  F^{cu}(q^0,q^+)  &
             \end{array}
\right.
\end{displaymath}
so that we can write the system (\ref{differential system2}) in the simpler form
\begin{equation}
\left\{
             \begin{array}{lr}
             {u'^0} = \mathfrak{C}^0 u^0 + \varphi^0(u^0,u^+,  F^{cu}(u^0,u^+)) + \psi^0(u)u^-, \\
             {u'^+} = \mathfrak{C}^+ u^+  + \varphi^+(u^0,u^+,  F^{cu}(u^0,u^+)) + \psi^+(u)u^-, \\
             {u'^-} = \mathfrak{C}^- u^-  + \psi^-(u)u^-,
             \end{array}
\right.\nonumber
\end{equation}
where
the functions $\psi^0, \psi^+$ are $C^{l}$-smooth and $ \psi^-$ is  $C^{l-1}$-smooth; in addition, all the functions $\psi^0, \psi^+, \psi^-$ are vanishing at the origin, i.e., $\psi(0)=0$.
\end{corollary}

In fact, a stronger version of reduction theorem holds \cite{Shilnikov}:
\begin{theorem}\label{Reduction Theorem}\emph{(\textbf{Reduction Theorem})}
By a $C^{l-1}$-smooth transformation the system (\ref{differential system2}) can be
locally reduced to the simpler form
\begin{equation}\label{simple forms6}
\left\{
             \begin{array}{lr}
             {u'^0} = \mathfrak{C}^0 u^0 + \psi^0_0(u)u^0 + \psi^0_+(u)u^+ + \psi^0_-(u)u^-, \\
             {u'^+} = \mathfrak{C}^+ u^+  + \psi^+(u)u^+, \\
             {u'^-} = \mathfrak{C}^- u^-  + \psi^-(u)u^-,
             \end{array}
\right.\nonumber
\end{equation}
where
the functions $\psi^0_0, \psi^0_+, \psi^0_-$ are $C^{l-1}$-smooth and vanishing at the origin;
the functions $ \psi^+, \psi^-$ are $C^{l}$-smooth and vanishing at the origin. Furthermore, $\psi^0_+$ vanishes identically at $u^-=0$, and   $\psi^0_-$ vanishes identically at $u^+=0$.
\end{theorem}

However, we will not utilize this stronger version of reduction theorem in this paper, because the concomitant computational capacity is large. Instead, we will directly utilize
a version of theorem of  center manifold.

When we investigate the orbits on the center-unstable manifold, that is, the orbits such that $q^-  =  F^{cu}(q^0,q^+)$ or $u^- = 0$, the system is reduced to a form
\begin{equation}\label{simple forms6}
\left\{
             \begin{array}{lr}
             {u'^0} = \mathfrak{C}^0 u^0 + \varphi^0(u^0,u^+,  F^{cu}(u^0,u^+)), \\
             {u'^+} = \mathfrak{C}^+ u^+  + \varphi^+(u^0,u^+,  F^{cu}(u^0,u^+)).
             \end{array}
\right.
\end{equation}

Note that this system is of class $C^l$.
For this kind of system, it is well known that the following theorems hold.
\begin{theorem}\label{Center Manifold}\emph{\textbf{(Center Manifold)}}
Consider the system (\ref{simple forms6}), in a small neighborhood of $0$
there exists an $n_0$-dimensional invariant center manifold $\mathcal{W}^{c}_{loc}| u^+ = F^{c}(u^0)$ of class $C^l$, which contains $0$ and which is tangent to the subspace
$u^+=0$ at $0$. The manifold $\mathcal{W}^{c}_{loc}$ contains all orbits which stay in a small
neighborhood of $0$ for all times. Though the center manifold
is not defined uniquely, for any two manifolds $\mathcal{W}^{c}_{1}$ and $\mathcal{W}^{c}_{2}$
the functions
$F^{c}_{1}$ and $F^{c}_{2}$
 have the same Taylor expansion at  $0$ (and at each point whose
orbit stays in a small neighborhood of $0$). In
the case where the system is $C^\infty$-smooth, the center manifold has, in
general, only finite smoothness.
\end{theorem}

Note also that the condition of invariance of the manifold $\mathcal{W}^{c}_{loc}$  may be
expressed as
\begin{equation}
u'^+ = \frac{\partial F^{c}(u^0)}{\partial u^0}u'^0   ~~~ ~~when ~~u^+  =  F^{c}(u^0)\nonumber
\end{equation}
or
\begin{equation}
\begin{array}{c}
  \mathfrak{C}^+ F^{c}(u^0) + \varphi^+(u^0,F^{cu}(u^0,F^{c}(u^0),F^{c}(u^0))) = \\
\frac{\partial F^{c}(u^0)}{\partial u^0}[\mathfrak{C}^0 u^0 + \varphi^0(u^0,F^{c}(u^0),F^{cu}(u^0,F^{c}(u^0)))].
\end{array}\nonumber
\end{equation}
The above relation  yields an algorithm for computing the center
manifolds which will be used in the following.

The existence of a center manifold allows some problems related to the nonhyperbolic
equilibrium to be reduced to the study of an $n_0$-dimensional system
\begin{equation}\label{center system}
{u'^0} = \mathfrak{C}^0 u^0 + \varphi^0 \left(u^0,F^{c}(u^0),  F^{cu}\left(u^0,F^{c}(u^0) \right) \right)
\end{equation}

In particular, when an orbit on the center-unstable manifold approaches the origin for $\tau\rightarrow -\infty$, we have the following result \cite{zbMATH03727247}:
\begin{theorem}\label{centerapproximatetheorem}
Let $(u^0(\tau), u^+(\tau))$ be a solution of the system (\ref{simple forms6}). Suppose that
\begin{displaymath}
(u^0(\tau), u^+(\tau))\rightarrow 0 ~~~~~~~~~as ~~~~~~\tau \rightarrow -\infty,
\end{displaymath}
then there exists a solution $\upsilon(\tau)$ of the system (\ref{center system}) such that as $\tau\rightarrow -\infty$
\begin{equation}\label{centerapproximate}
\left\{
             \begin{array}{lr}
             {u^0}(\tau) = \upsilon(\tau) + O(e^{\sigma \tau}), \\
             {u^+}(\tau) = F^{c}(\upsilon(\tau))  + O(e^{\sigma \tau}).
             \end{array}
\right.
\end{equation}
where $\sigma >0$ is a constant depending only on $\mathfrak{C}^+$.
\end{theorem}

As a conclusion, it is well known that an orbit on the
unstable manifold of a hyperbolic equilibrium point exponentially approaches the hyperbolic equilibrium point, i.e., the following result holds.
\begin{corollary}\label{hyperbolicapproximatetheorem}
Consider the  system
(\ref{differential system2}), assume that the origin
$q=0$ is a hyperbolic equilibrium point. Let $q(\tau)$ be a solution on the
unstable manifold of the origin, i.e.,
\begin{equation*}
    q(\tau)\rightarrow 0 ~~~~~~~~~as ~~~\tau \rightarrow -\infty,
\end{equation*}
then there are two positive constants $c$ and $\sigma$  such that
\begin{equation}
 dist(q(\tau),0) \leq c e^{\sigma \tau}\nonumber
\end{equation}
for sufficiently small $\tau$,
and $\sigma $  depends only on $\mathfrak{C}^+$.
\end{corollary}

      \section{Plane Equilibrium Points}\label{Plane Equilibrium Points}
      \indent\par

In this section, we will discuss some aspects of plane equilibrium points. In particular, for those orbits tending to the equilibrium points, we estimate their rate of tending to the equilibrium points.

Let us consider an autonomous system on the plane $\mathbb{R}^2$
\begin{equation}\label{system on the plane}
\left\{
             \begin{array}{lr}
             \zeta' = f(\zeta,\eta), \\
             \eta' = g(\zeta,\eta),
             \end{array}
\right.
\end{equation}
where $f,g$ are continuous for small $\zeta,\eta$ and
\begin{displaymath}
f(0,0)=g(0,0)= 0 .
\end{displaymath}

One can introduce polar coordinates
\begin{displaymath}
\zeta=\rho \cos \vartheta, \eta=\rho \sin \vartheta
\end{displaymath}
to transform the system (\ref{system on the plane})
into
\begin{equation}\label{system on the planepolar}
\left\{
             \begin{array}{lr}
             \rho' = f(\rho \cos \vartheta,\rho \sin \vartheta)\cos \vartheta +g(\rho \cos \vartheta,\rho \sin \vartheta)\sin \vartheta, \\
             \rho\vartheta' = g(\rho \cos \vartheta,\rho \sin \vartheta)\sin \vartheta -f(\rho \cos \vartheta,\rho \sin \vartheta)\cos \vartheta.
             \end{array}
\right.
\end{equation}

As in \cite{Hartman},
a direction $\vartheta=\vartheta_0$ at the origin is called \textbf{characteristic} for the system (\ref{system on the plane}), if
there exists a sequence $(\rho_1, \vartheta_1), (\rho_1, \vartheta_1), \cdots$ such that:
\begin{description}
  \item[1)] $
\rho_k \rightarrow 0+,  \vartheta_k \rightarrow \vartheta_0 ~~~~~~as ~~~k\rightarrow\infty;
$
  \item[2)] $(f_k,g_k)$  is not $(0,0)$, and the angle ($mod ~\pi$) between the vectors $(f_k,g_k)$ and $(\cos \vartheta_k,  \sin \vartheta_k)$ tends to zero as $k\rightarrow\infty$, i.e.,
\begin{equation}\label{characteristic}
    \begin{array}{lr}
             \frac{g_k\cos \vartheta_k-f_k\sin \vartheta_k}{\sqrt{f^2_k+g^2_k}} \rightarrow 0 ~~~~~~as ~~~k\rightarrow\infty.
             \end{array}
\end{equation}
\end{description}
Where $\rho_k \rightarrow 0+$ denotes $\rho_k \rightarrow 0$ and $\rho_k > 0$;
$(f_k,g_k)$ is the vector field $(f,g)$ evaluated at $(\zeta_k, \eta_k)=(\rho_k \cos \vartheta_k, \rho_k \sin \vartheta_k)$.

The following two lemmas in \cite{Hartman} are important to the investigation of plane equilibrium points.
\begin{lemma}\label{near w }
Let $f,g$ be continuous for small $\zeta,\eta$ and $f^2(\zeta,\eta)+g^2(\zeta,\eta) >0 $ except at the origin, that is, the origin is an isolated equilibrium point of the system (\ref{system on the plane}). Let the system (\ref{system on the plane}) possess a solution $(\zeta(\tau),\eta(\tau))$ for $-\infty < \tau \leq 0$ such that
\begin{displaymath}
\zeta^2(\tau)+\eta^2(\tau)  \rightarrow 0+ ~~~~~~as ~~~\tau\rightarrow -\infty,
\end{displaymath}
Let $\rho(\tau)=\sqrt{f^2(\zeta,\eta)+g^2(\zeta,\eta)}>0$  and $\vartheta(\tau)$ a continuous determination of $\arctan \frac{\eta(\tau)}{\zeta(\tau)}$. Let $\vartheta= \vartheta_0$ be a noncharacteristic direction. Then either $\vartheta'(\tau)>0$ or $\vartheta'(\tau)<0$ for all $\tau$ near $-\infty$ for which $\vartheta(\tau)= \vartheta_0$ $mod ~2\pi$.
\end{lemma}
\begin{lemma}\label{vartheta}
Let $f,g$ and $(\zeta(\tau),\eta(\tau))$ be as in Lemma \ref{near w }. Suppose that every $\vartheta$-interval, $\alpha< \vartheta< \beta$, contains a noncharacteristic direction. Then either
\begin{equation}\label{vartheta1}
    \begin{array}{lr}
              \vartheta_0= \lim_{\tau\rightarrow -\infty} \vartheta(\tau)~~~~~~exists ~(and~ is~ finite)
             \end{array}
\end{equation}
or $(\zeta(\tau),\eta(\tau))$ is a spiral; i.e.,
\begin{equation}\label{vartheta2}
    \begin{array}{lr}
     |\vartheta(\tau)|\rightarrow \infty~~~~~~as ~~~\tau\rightarrow -\infty.
             \end{array}
\end{equation}
In the case (\ref{vartheta1}), $\vartheta(\tau)= \vartheta_0$ is a characteristic direction.
\end{lemma}
\begin{figure}
  \center
  \includegraphics[width=6cm]{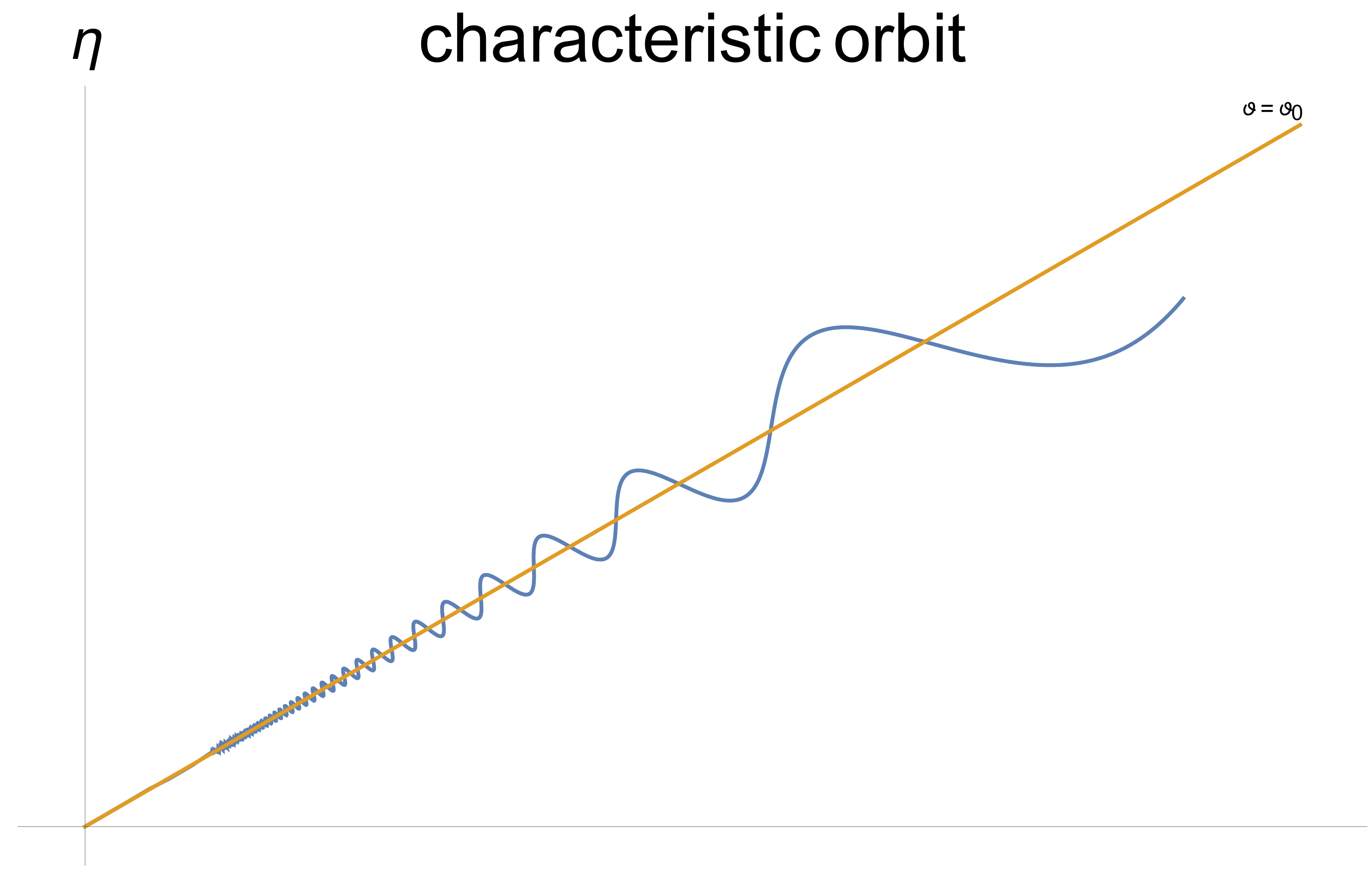}~~~~~~~~~~\includegraphics[width=6cm]{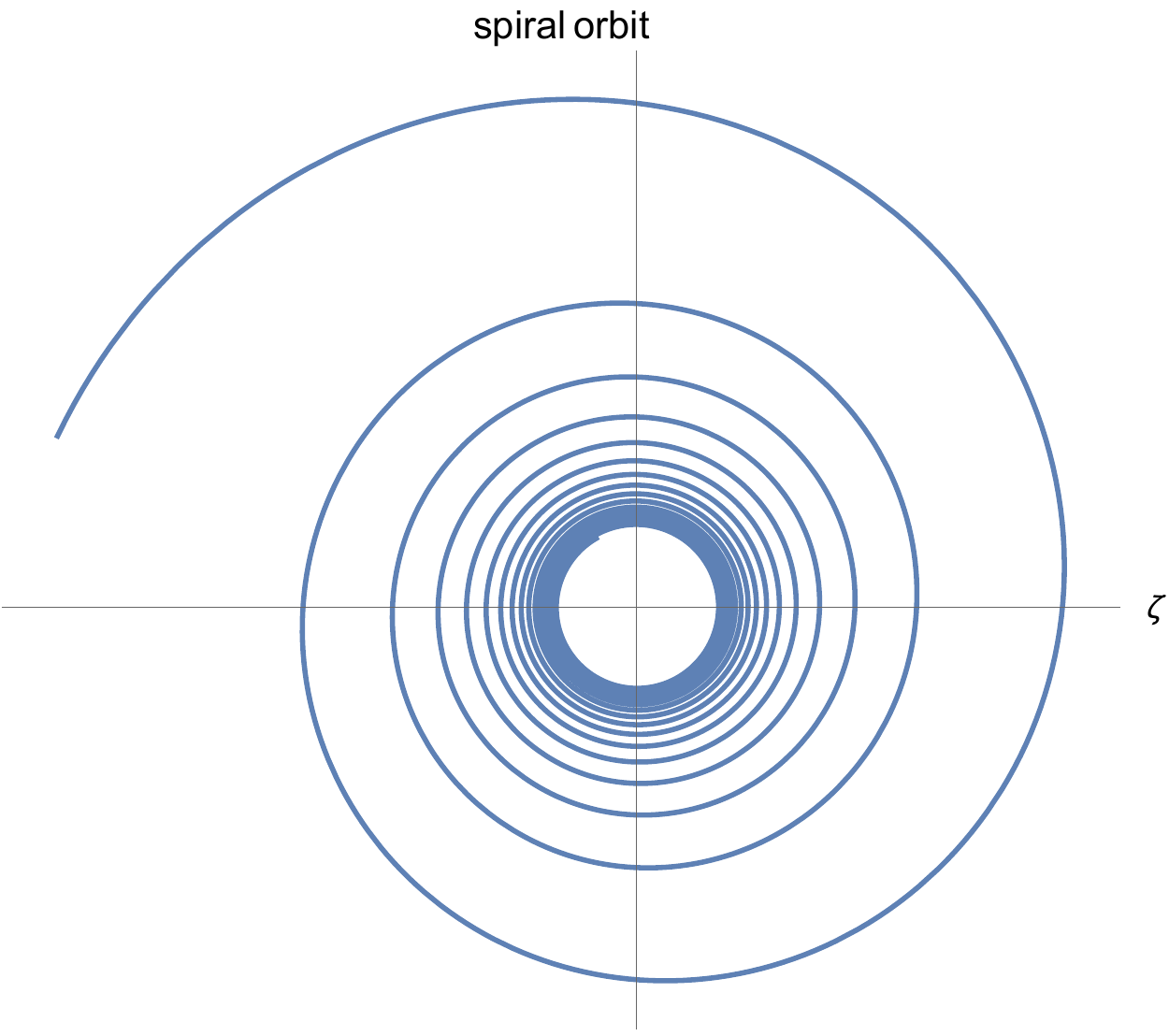}\\
  \caption{{Lemma} \ref{vartheta}}
\end{figure}

In the following, let us further assume that the functions $f,g$ are 
\begin{equation}\label{functioncondition}
\left\{
             \begin{array}{lr}
             f = P_m(\zeta,\eta) + p_m(\zeta,\eta), \\
             g = Q_m(\zeta,\eta) + q_m(\zeta,\eta),
             \end{array}
\right.\nonumber
\end{equation}
where $P_m, Q_m$ are homogeneous polynomials of degree $m>1$ and
\begin{displaymath}
p^2(\tau)+q^2(\tau)  = o (\rho^{2m}) ~~~~~~as ~~~\rho\rightarrow 0.
\end{displaymath}

In terms of polar coordinates, define
\begin{equation}\label{functioncondition1}
\left\{
             \begin{array}{lr}
             \Phi(\vartheta) = \rho^{-m}\left(P_m(\rho \cos \vartheta,\rho \sin \vartheta)\cos \vartheta +Q_m(\rho \cos \vartheta,\rho \sin \vartheta)\sin \vartheta\right), \\
             \Psi(\vartheta) = \rho^{-m}\left(Q_m(\rho \cos \vartheta,\rho \sin \vartheta)\cos \vartheta -P_m(\rho \cos \vartheta,\rho \sin \vartheta)\sin \vartheta\right);
             \end{array}
\right.\nonumber
\end{equation}
then $\Phi,\Psi$ are homogeneous polynomials of $\cos \vartheta, \sin \vartheta$ of degree $m+1$.

In terms of polar coordinates, (\ref{system on the planepolar}) can be written as
\begin{equation}\label{system on the planepolarmain}
\left\{
             \begin{array}{lr}
             \rho' = \rho^{m}\Phi(\vartheta) + o (\rho^{m}), \\
             \vartheta' = \rho^{m-1}\Psi(\vartheta) + o (\rho^{m-1}).
             \end{array}
\right.
\end{equation}

\begin{theorem}\label{varthetamain}
Assume $\Psi(\vartheta)\not\equiv 0$ and $m$ is  an even number. Let $f,g$ and $(\zeta(\tau),\eta(\tau))$ be as in Lemma \ref{near w }.
Then
\begin{equation}\label{vartheta11}
    \begin{array}{lr}
              \vartheta_0= \lim_{\tau\rightarrow -\infty} \vartheta(\tau)~~~~~~exists ~(and~ is~ finite)
             \end{array}
\end{equation}
and $\Psi(\vartheta_0)= 0$.
\end{theorem}
{\bf Proof.} 
Our first goal is to show  (\ref{vartheta11}).
By Lemma \ref{vartheta}, it suffices to get rid of the possibility of (\ref{vartheta2}).

By using reduction to absurdity, suppose that $(\zeta(\tau),\eta(\tau))$ is a spiral, i.e., (\ref{vartheta2}) holds.

Without loss of generality, suppose
\begin{equation}\label{vartheta21}
    \begin{array}{lr}
     \vartheta(\tau)\rightarrow \infty~~~~~~as ~~~\tau\rightarrow -\infty.
             \end{array}
\end{equation}

According to $\Psi(\vartheta)\not\equiv 0$, it has only a
finite number of zeros ($mod ~2\pi$).
Because  $\Psi$ are homogeneous polynomials of $\cos \vartheta, \sin \vartheta$ of degree $m+1$, we have
\begin{equation}\label{Psi}
         \begin{array}{lr}
             \Psi(\vartheta + \pi) = - \Psi(\vartheta).
             \end{array}
\end{equation}
As a result, the plane is split into several $\vartheta$-intervals, $\vartheta_k< \vartheta< \vartheta_{k+1}$ ($k=1,\cdots,n$), such that
\begin{displaymath}
 n\geq 1,\vartheta_{n+1}=\vartheta_1+2\pi, \Psi(\vartheta_k)=0
\end{displaymath}
and
\begin{displaymath}
\Psi(\vartheta)\neq 0, ~~~~~~for~~ \vartheta_k< \vartheta< \vartheta_{k+1}.
\end{displaymath}

It follows from (\ref{Psi}) that we can assume $
\Psi(\vartheta)> 0$ for $\vartheta_1< \vartheta< \vartheta_{2}.
$
Then there exists a $\vartheta$-interval, $\alpha_1< \vartheta< \alpha_{2}$ included in $\vartheta_1< \vartheta< \vartheta_{2}$, such that
\begin{displaymath}
\Psi(\vartheta) >\sigma, ~~~~~~for~~ \alpha_1< \vartheta< \alpha_{2},
\end{displaymath}
where $\sigma >0$ is a constant.

Note that $\frac{o (\rho^{m-1})}{\rho^{m-1}}$ tend to zero as $\rho\rightarrow 0$ uniformly in $\vartheta$. It follows that,  there exists a  sufficiently  small $\rho_0>0$ such that
\begin{equation}
\rho^{m-1}\sigma + o (\rho^{m-1}) >0 ~~~~~~for~any~0<\rho<\rho_0.\nonumber
\end{equation}

Since there exists a  real number $\tau_0<0$ such that
\begin{equation}
0<\rho(\tau)<\rho_0 ~~~~~~for~any~\tau \leq \tau_0.\nonumber
\end{equation}
Taking  into consideration (\ref{vartheta21}), we know that there exists a sequence of $\tau$-intervals $(\beta_1, \gamma_1), (\beta_2, \gamma_2), \cdots $  such that
\begin{equation}\label{contradiction}
\left\{
             \begin{array}{lr}
               \gamma_{k+1} <\beta_k < \gamma_1< \tau_0,    & \\
           \beta_k,  \gamma_k \rightarrow -\infty   &as ~~~k\rightarrow\infty; \\
              \vartheta'(\tau) >0  &for~~ \tau\in [\beta_k,  \gamma_{k}]\\
           \left[\vartheta(\beta_k),  \vartheta(\gamma_k)\right] \subset (\vartheta_1 - 2 n_k \pi, \vartheta_{2}- 2 n_k \pi)  &for~some~ n_k\in \mathbb{N}.
             \end{array}
\right.
\end{equation}

We claim that
\begin{equation}
\vartheta(\tau) \leq \vartheta(\frac{\beta_1+\gamma_1}{2}) <\vartheta(\gamma_1)\nonumber
\end{equation}
for any $\tau \leq \frac{\beta_1+\gamma_1}{2}$. Obviously, this is a contradiction to (\ref{contradiction}). So (\ref{vartheta11}) will be proved if we can show the claim.

Set
\begin{equation}
\Omega= \{\tilde{\tau}| \vartheta(\tau) \leq \vartheta(\frac{\beta_1+\gamma_1}{2}) ~~for~any~\tau \in (\tilde{\tau}, \gamma_1]\} .\nonumber
\end{equation}
Let $\tau_i=\inf\Omega$ be  the infimum of above set.
Following from (\ref{contradiction}), it is clear that $\tau_i<\beta_1$.

The above claim will be proved by showing that $\tau_i=-\infty$. If otherwise, then $\tau_i$ is a certain negative number.
It follows that
\begin{equation}
\vartheta(\beta_1) <\vartheta(\tau_i) = \vartheta(\frac{\beta_1+\gamma_1}{2}) <\vartheta(\gamma_1)\nonumber
\end{equation}
However, it is easy to prove that  $\tau_i$ is not the infimum of $\Omega$ by above inequality. This leads to a contradiction.

Our task now is to show $\Psi(\vartheta_0)=0$. Following from Lemma \ref{vartheta}, $\vartheta(\tau)= \vartheta_0$ is a characteristic direction. If $\Psi(\vartheta_0)\neq0$, It is straightforward to show that  (\ref{characteristic}) attributes to the following
\begin{equation}
    \begin{array}{lr}
             \frac{\Psi(\vartheta_0)}{\sqrt{\Phi^2(\vartheta_0)+\Psi^2(\vartheta_0)}} =0.\nonumber
             \end{array}
\end{equation}
This leads to a contradiction.

The theorem is now evident from what we have proved.

$~~~~~~~~~~~~~~~~~~~~~~~~~~~~~~~~~~~~~~~~~~~~~~~~~~~~~~~~~~~~~~~~~~~~~~~~~~~~~~~~~~~~~~~~~~~~~~~~~~~~~~~~~~~~~~~~~~~~~~~~~~~~~~~~~~~~~~~~~~~~~~~~~~~~~~~~~~~~~~~~~~~\Box$

Using the same argument as in the proof of above theorem, we can prove the following more general result£»
\begin{corollary}
Let $f,g$ and $(\zeta(\tau),\eta(\tau))$ be as in Lemma \ref{near w }. Assume $\Psi(\vartheta)$ has both positive and negative values.
Then
\begin{equation}
    \begin{array}{lr}
              \vartheta_0= \lim_{\tau\rightarrow -\infty} \vartheta(\tau)~~~~~~exists ~(and~ is~ finite)\nonumber
             \end{array}
\end{equation}
and $\Psi(\vartheta_0)= 0$.
\end{corollary}

So a spiral of the system (\ref{system on the plane}) can occur  for $\Psi$ be invariably nonnegative or nonpositive and $m$ is  an odd number.

\begin{theorem}\label{rhoamain}
Under the  conditions in Theorem \ref{varthetamain}, if $\Phi(\vartheta_0) \neq 0$,
then
\begin{displaymath}
 \rho = (\frac{1}{(m-1)\Phi(\vartheta_0)})^\frac{1}{m-1}(\frac{1}{-\tau})^\frac{1}{m-1} +o((\frac{1}{-\tau})^\frac{1}{m-1}).
\end{displaymath}
\end{theorem}
{\bf Proof.} 
Let us consider the first equation of (\ref{system on the planepolarmain}):
\begin{equation}
\rho' = \rho^{m}\Phi(\vartheta) + o (\rho^{m}).\nonumber
\end{equation}
According to L'H\'{o}pital's rule, it follows  that
\begin{displaymath}
\lim_{\tau \rightarrow -\infty} \frac{1}{\tau \rho^{m-1}}= (1-m)\Phi(\vartheta_0),
\end{displaymath}
or
\begin{displaymath}
 \rho= \left(\frac{1}{ (1-m)c_{m}\Phi(\vartheta_0)\tau}\right)^\frac{1}{m-1} +o((\frac{1}{-\tau})^\frac{1}{m-1}).
\end{displaymath}

The proof of the theorem is now complete.

$~~~~~~~~~~~~~~~~~~~~~~~~~~~~~~~~~~~~~~~~~~~~~~~~~~~~~~~~~~~~~~~~~~~~~~~~~~~~~~~~~~~~~~~~~~~~~~~~~~~~~~~~~~~~~~~~~~~~~~~~~~~~~~~~~~~~~~~~~~~~~~~~~~~~~~~~~~~~~~~~~~~\Box$

Similarly, when $m$ is  an odd number, one can prove the following theorem:
\begin{theorem}\label{rhoamain1}
Let $f,g$ and $(\zeta(\tau),\eta(\tau))$ be as in Lemma \ref{near w }.
If $\Phi(\vartheta) \neq 0$ for any $\vartheta$,
then there exists a positive number $c$  such that
\begin{displaymath}
 \rho \leq c(\frac{1}{-\tau})^\frac{1}{m-1}.
\end{displaymath}
\end{theorem}
{\bf Proof.} 

By $\Phi(\vartheta) \neq 0$ for any $\vartheta$, it follows that  there exists a positive number $\sigma$  such that
\begin{displaymath}
|\Phi(\vartheta)|\geq \sigma.
\end{displaymath}
According to
\begin{displaymath}
\rho(\tau)  \rightarrow 0+ ~~~~~~as ~~~\tau\rightarrow -\infty,
\end{displaymath}
it follows that
\begin{displaymath}
\Phi(\vartheta)\geq \sigma.
\end{displaymath}

Let us consider the first equation of (\ref{system on the planepolarmain}):
\begin{equation}
\rho' = \rho^{m}\Phi(\vartheta) + o (\rho^{m}).\nonumber
\end{equation}

It is clear that there exists a  real number $\tau_0<0$ such that
\begin{equation}
\rho^{m}\Phi(\vartheta) + o (\rho^{m}) \geq \frac{\sigma}{2} \rho^{m} ~~~~~~for~any~\tau \leq \tau_0.\nonumber
\end{equation}
Then
\begin{equation}
\rho^{1-m}(\tau)  \geq \rho^{1-m}(\tau_0) +\frac{\sigma}{2} (m-1)(\tau_0-\tau) ~~~~~~for~any~\tau \leq \tau_0.\nonumber
\end{equation}

As a result, it is evident to see that  the theorem holds.

$~~~~~~~~~~~~~~~~~~~~~~~~~~~~~~~~~~~~~~~~~~~~~~~~~~~~~~~~~~~~~~~~~~~~~~~~~~~~~~~~~~~~~~~~~~~~~~~~~~~~~~~~~~~~~~~~~~~~~~~~~~~~~~~~~~~~~~~~~~~~~~~~~~~~~~~~~~~~~~~~~~~\Box$

  \end{appendices}


\section*{Acknowledgements}
\indent\par
I wish to express my gratitude to  Shanzhen Chen, Rick Moeckel, Shiqing Zhang and Shuqiang Zhu for their  very helpful
discussions. I am also deeply indebted to the  reviewer  for carefully reading the original manuscript and suggesting very helpful improvements with paying enormous energy. I would like to acknowledge the support from  NSFC (No.11701464).

\newpage

\bibliographystyle{plain}


\begin{thebibliography}{10}

\bibitem{Albouy2012Some}
A. Albouy, H. Cabral, and A. Santos.
\newblock Some problems on the classical n-body problem.
\newblock {\em Celestial Mechanics \& Dynamical Astronomy}, 113(4):369--375,
  2012.

\bibitem{zbMATH06074021}
A. Albouy and V. Kaloshin.
\newblock Finiteness of central configurations of five bodies in the plane.
\newblock {\em Annals of Mathematics}, 176(1):535--588, 2012.

\bibitem{zbMATH03601149}
V.~I. {Arnold}.
\newblock {\em {Mathematical methods of classical mechanics. Translated by K.
  Vogtman and A. Weinstein}}, volume~60.
\newblock Springer, New York, NY, 1978.

\bibitem{Geometricalmethods}
V.~I. Arnold.
\newblock {\em Geometrical methods in the theory of ordinary differential
  equations}.
\newblock Springer-Verlag: New York, 1988.

\bibitem{zbMATH05031968}
V.~I. {Arnold}, V.~V. {Kozlov}, and A.~I. {Neishtadt}.
\newblock {\em {Mathematical aspects of classical and celestial mechanics.
  Transl. from the Russian by E. Khukhro. 3rd revised}}.
\newblock Berlin: Springer, 3rd revised ed. edition, 2006.

\bibitem{zbMATH03727247}
J. {Carr}.
\newblock {\em {Applications of centre manifold theory}}, volume~35.
\newblock Springer, New York, NY, 1981.

\bibitem{Chazy}
J. Chazy.
\newblock Sur certaines trajectoires du probleme des n corps.
\newblock {\em Bull. Astron.}, 35:321--389, 1918.

\bibitem{Hampton2006Finiteness}
M. Hampton and R. Moeckel.
\newblock Finiteness of relative equilibria of the four-body problem.
\newblock {\em Inventiones mathematicae}, 163(2):289--312, 2006.

\bibitem{Hartman}
P. Hartman.
\newblock {\em Ordinary differential equations , Classics in applied
  mathematics 38}.
\newblock Society for Industrial and Applied Mathematics, 2002.

\bibitem{zbMATH00579429}
C. {Hassell} and E. {Rees}.
\newblock {The index of a constrained critical point}.
\newblock {\em {American Mathematical Monthly}}, 100(8):772--778, 1993.

\bibitem{Hirsch1970Invariant}
M.~W. Hirsch, C.~C. Pugh, and M. Shub.
\newblock Invariant manifolds.
\newblock {\em Bulletin of the American Mathematical Society},
  76(5):1015--1019, 1970.

\bibitem{Leandro2003Finiteness}
E.~S.~G. Leandro.
\newblock Finiteness and bifurcations of some symmetrical classes of central
  configurations.
\newblock {\em Archive for Rational Mechanics \& Analysis}, 167(2):147--177,
  2003.

\bibitem{Mcgehee1974Triple}
R. Mcgehee.
\newblock Triple collision in the collinear three-body problem.
\newblock {\em Inventiones Mathematicae}, 27(3):191--227, 1974.

\bibitem{moeckel1990central}
R. Moeckel.
\newblock On central configurations.
\newblock {\em Mathematische Zeitschrift}, 205(1):499--517, 1990.

\bibitem{moeckel2014lectures}
R. Moeckel.
\newblock Lectures on central configurations.
\newblock http://www.math.umn.edu/~rmoeckel, 2014.

\bibitem{Moser2010Regularization}
J.~Moser.
\newblock Regularization of kepler's problem and the averaging method on a
  manifold.
\newblock {\em Communications on Pure and Applied Mathematics}, 23(4):609--636,
  2010.

\bibitem{palmore1975classifying}
J.~I. Palmore.
\newblock Classifying relative equilibria. ii.
\newblock {\em Bulletin of the American Mathematical Society}, 81(2):489--491,
  1975.

\bibitem{Palmore1976Measure}
J.~I. Palmore.
\newblock Measure of degenerate relative equilibria. i.
\newblock {\em Annals of Mathematics}, 104(3):421--429, 1976.

\bibitem{zbMATH03657180}
D.~G. {Saari}.
\newblock {On the role and the properties of n body central configurations}.
\newblock {\em {Celestial Mechanics}}, 21:9--20, 1980.

\bibitem{Saari1984The}
D.~G. Saari.
\newblock The manifold structure for collision and for hyperbolic-parabolic
  orbits in the n-body problem.
\newblock {\em Journal of Differential Equations}, 55(3):300--329, 1984.

\bibitem{Saari1981On}
D.~G. Saari and Neal~D Hulkower.
\newblock On the manifolds of total collapse orbits and of completely parabolic
  orbits for the n-body problem.
\newblock {\em Journal of Differential Equations}, 41(1):27--43, 1981.

\bibitem{Shilnikov}
L. P. Shilnikov, A.~Shilnikov, D.~Turaev, and L.~Chua.
\newblock {\em Methods of Qualitative Theory in Nonlinear Dynamics. Part I.}
\newblock World Scientific Press, 1998.

\bibitem{C1967Lectures}
C. L. Siegel.
\newblock {\em Lectures On The Singularities Of The Three-Body Problem}.
\newblock Tata Institute of Fundamental Research, Lectures on Mathematics,
  1967.

\bibitem{zbMATH01516323}
S. {Smale}.
\newblock {Mathematical problems for the next century}.
\newblock {\em The Mathematical Intelligencer}, 20(2):7--15, 1998.

\bibitem{smale1970topology}
S. Smale.
\newblock Topology and mechanics. ii.
\newblock {\em Inventiones mathematicae}, 11(1):45--64, 1970.

\bibitem{Sperling1970On}
H.~J. Sperling.
\newblock On the real singularities of the n-body problem.
\newblock {\em Journal F\"{u}r Die Reine Und Angewandte Mathematik},
  1970(245):15--40, 1970.

\bibitem{Sundman}
K.~F. Sundman.
\newblock Memoire sur le probleme des trois corps.
\newblock {\em Acta Mathematica}, 36:105--179, 1912.

\bibitem{wintner1941analytical}
A. Wintner.
\newblock {\em The analytical foundations of celestial mechanics}.
\newblock Princeton University Press, 1941.

\end{thebibliography}


\end{document}